

\input amstex
\expandafter\ifx\csname mathdefs.tex\endcsname\relax
  \expandafter\gdef\csname mathdefs.tex\endcsname{}
\else \message{Hey!  Apparently you were trying to
  \string twice.   This does not make sense.} 
\errmessage{Please edit your file (probably \jobname.tex) and remove
any duplicate ``\string\input'' lines} \fi




\catcode`\X=12\catcode`\@=11

\def\n@wcount{\alloc@0\count\countdef\insc@unt}
\def\n@wwrite{\alloc@7\write\chardef\sixt@@n}
\def\n@wread{\alloc@6\read\chardef\sixt@@n}
\def\r@s@t{\relax}\def\v@idline{\par}\def\@mputate#1/{#1}
\def\l@c@l#1X{\firstpart.#1}\def\gl@b@l#1X{#1}\def\t@d@l#1X{{}}

\def\crossrefs#1{\ifx\all#1\let\tr@ce=\all\else\def\tr@ce{#1,}\fi
   \n@wwrite\cit@tionsout\openout\cit@tionsout=\jobname.cit 
   \write\cit@tionsout{\tr@ce}\expandafter\setfl@gs\tr@ce,}
\def\setfl@gs#1,{\def\@{#1}\ifx\@\empty\let\next=\relax
   \else\let\next=\setfl@gs\expandafter\xdef
   \csname#1tr@cetrue\endcsname{}\fi\next}
\def\m@ketag#1#2{\expandafter\n@wcount\csname#2tagno\endcsname
     \csname#2tagno\endcsname=0\let\tail=\all\xdef\all{\tail#2,}
   \ifx#1\l@c@l\let\tail=\r@s@t\xdef\r@s@t{\csname#2tagno\endcsname=0\tail}\fi
   \expandafter\gdef\csname#2cite\endcsname##1{\expandafter
     \ifx\csname#2tag##1\endcsname\relax?\else\csname#2tag##1\endcsname\fi
     \expandafter\ifx\csname#2tr@cetrue\endcsname\relax\else
     \write\cit@tionsout{#2tag ##1 cited on page \folio.}\fi}
   \expandafter\gdef\csname#2page\endcsname##1{\expandafter
     \ifx\csname#2page##1\endcsname\relax?\else\csname#2page##1\endcsname\fi
     \expandafter\ifx\csname#2tr@cetrue\endcsname\relax\else
     \write\cit@tionsout{#2tag ##1 cited on page \folio.}\fi}
   \expandafter\gdef\csname#2tag\endcsname##1{\expandafter
      \ifx\csname#2check##1\endcsname\relax
      \expandafter\xdef\csname#2check##1\endcsname{}%
      \else\immediate\write16{Warning: #2tag ##1 used more than once.}\fi
      \multit@g{#1}{#2}##1/X%
      \write\t@gsout{#2tag ##1 assigned number \csname#2tag##1\endcsname\space
      on page \number\count0.}%
   \csname#2tag##1\endcsname}}
\def\multit@g#1#2#3/#4X{\def\t@mp{#4}\ifx\t@mp\empty%
      \global\advance\csname#2tagno\endcsname by 1 
      \expandafter\xdef\csname#2tag#3\endcsname
      {#1\number\csname#2tagno\endcsnameX}%
   \else\expandafter\ifx\csname#2last#3\endcsname\relax
      \expandafter\n@wcount\csname#2last#3\endcsname
      \global\advance\csname#2tagno\endcsname by 1 
      \expandafter\xdef\csname#2tag#3\endcsname
      {#1\number\csname#2tagno\endcsnameX}
      \write\t@gsout{#2tag #3 assigned number \csname#2tag#3\endcsname\space
      on page \number\count0.}\fi
   \global\advance\csname#2last#3\endcsname by 1
   \def\t@mp{\expandafter\xdef\csname#2tag#3/}%
   \expandafter\t@mp\@mputate#4\endcsname
   {\csname#2tag#3\endcsname\lastpart{\csname#2last#3\endcsname}}\fi}
\def\t@gs#1{\def\all{}\m@ketag#1e\m@ketag#1s\m@ketag\t@d@l p
   \m@ketag\gl@b@l r \n@wread\t@gsin
   \openin\t@gsin=\jobname.tgs \re@der \closein\t@gsin
   \n@wwrite\t@gsout\openout\t@gsout=\jobname.tgs }
\outer\def\localtags{\t@gs\l@c@l}
\outer\def\globaltags{\t@gs\gl@b@l}
\outer\def\newlocaltag#1{\m@ketag\l@c@l{#1}}
\outer\def\newglobaltag#1{\m@ketag\gl@b@l{#1}}

\newif\ifpr@ 
\def\m@kecs #1tag #2 assigned number #3 on page #4.%
   {\expandafter\gdef\csname#1tag#2\endcsname{#3}
   \expandafter\gdef\csname#1page#2\endcsname{#4}
   \ifpr@\expandafter\xdef\csname#1check#2\endcsname{}\fi}
\def\re@der{\ifeof\t@gsin\let\next=\relax\else
   \read\t@gsin to\t@gline\ifx\t@gline\v@idline\else
   \expandafter\m@kecs \t@gline\fi\let \next=\re@der\fi\next}
\def\pretags#1{\pr@true\pret@gs#1,,}
\def\pret@gs#1,{\def\@{#1}\ifx\@\empty\let\n@xtfile=\relax
   \else\let\n@xtfile=\pret@gs \openin\t@gsin=#1.tgs \message{#1} \re@der 
   \closein\t@gsin\fi \n@xtfile}

\newcount\sectno\sectno=0\newcount\subsectno\subsectno=0
\newif\ifultr@local \def\ultralocal{\ultr@localtrue}
\def\firstpart{\number\sectno}
\def\lastpart#1{\ifcase#1 \or a\or b\or c\or d\or e\or f\or g\or h\or 
   i\or k\or l\or m\or n\or o\or p\or q\or r\or s\or t\or u\or v\or w\or 
   x\or y\or z \fi}

\def\resetall{\global\advance\sectno by 1\subsectno=0
   \gdef\firstpart{\number\sectno}\r@s@t}
\def\resetsub{\global\advance\subsectno by 1
   \gdef\firstpart{\number\sectno.\number\subsectno}\r@s@t}
\def\newsection#1\par{\resetall\vskip0pt plus.3\vsize\penalty-250
   \vskip0pt plus-.3\vsize\bigskip\bigskip
   \message{#1}\leftline{\bf#1}\nobreak\bigskip}
\def\subsection#1\par{\ifultr@local\resetsub\fi
   \vskip0pt plus.2\vsize\penalty-250\vskip0pt plus-.2\vsize
   \bigskip\smallskip\message{#1}\leftline{\bf#1}\nobreak\medskip}

\def\t@gsoff#1,{\def\@{#1}\ifx\@\empty\let\next=\relax\else\let\next=\t@gsoff
   \def\@@{p}\ifx\@\@@\else
   \expandafter\gdef\csname#1cite\endcsname##1{\zeigen{##1}}
   \expandafter\gdef\csname#1page\endcsname##1{?}
   \expandafter\gdef\csname#1tag\endcsname##1{\zeigen{##1}}\fi\fi\next}
\def\verbatimtags{\ifx\all\relax\else\expandafter\t@gsoff\all,\fi}
\def\zeigen#1{\hbox{$\langle$}#1\hbox{$\rangle$}}

\def\(#1){\edef\dot@g{\ifmmode\ifinner(\hbox{\noexpand\etag{#1}})
   \else\noexpand\eqno(\hbox{\noexpand\etag{#1}})\fi
   \else(\noexpand\ecite{#1})\fi}\dot@g}

\newif\ifbr@ck
\def\eat#1{}
\def\[#1]{\br@cktrue[\br@cket#1'X]}
\def\br@cket#1'#2X{\def\temp{#2}\ifx\temp\empty\let\next\eat
   \else\let\next\br@cket\fi
   \ifbr@ck\br@ckfalse\br@ck@t#1,X\else\br@cktrue#1\fi\next#2X}
\def\br@ck@t#1,#2X{\def\temp{#2}\ifx\temp\empty\let\neext\eat
   \else\let\neext\br@ck@t\def\temp{,}\fi
   \def\teemp{#1}\ifx\teemp\empty\else\rcite{#1}\fi\temp\neext#2X}
\def\resetbr@cket{\gdef\[##1]{[\rtag{##1}]}}
\def\references{\resetbr@cket\newsection References\par}

\newtoks\symb@ls\newtoks\s@mb@ls\newtoks\p@gelist\n@wcount\ftn@mber
    \ftn@mber=1\newif\ifftn@mbers\ftn@mbersfalse\newif\ifbyp@ge\byp@gefalse
\def\defm@rk{\ifftn@mbers\n@mberm@rk\else\symb@lm@rk\fi}
\def\n@mberm@rk{\xdef\m@rk{{\the\ftn@mber}}%
    \global\advance\ftn@mber by 1 }
\def\rot@te#1{\let\temp=#1\global#1=\expandafter\r@t@te\the\temp,X}
\def\r@t@te#1,#2X{{#2#1}\xdef\m@rk{{#1}}}
\def\b@@st#1{{$^{#1}$}}\def\str@p#1{#1}
\def\symb@lm@rk{\ifbyp@ge\rot@te\p@gelist\ifnum\expandafter\str@p\m@rk=1 
    \s@mb@ls=\symb@ls\fi\write\f@nsout{\number\count0}\fi \rot@te\s@mb@ls}
\def\byp@ge{\byp@getrue\n@wwrite\f@nsin\openin\f@nsin=\jobname.fns 
    \n@wcount\currentp@ge\currentp@ge=0\p@gelist={0}
    \re@dfns\closein\f@nsin\rot@te\p@gelist
    \n@wread\f@nsout\openout\f@nsout=\jobname.fns }
\def\m@kelist#1X#2{{#1,#2}}
\def\re@dfns{\ifeof\f@nsin\let\next=\relax\else\read\f@nsin to \f@nline
    \ifx\f@nline\v@idline\else\let\t@mplist=\p@gelist
    \ifnum\currentp@ge=\f@nline
    \global\p@gelist=\expandafter\m@kelist\the\t@mplistX0
    \else\currentp@ge=\f@nline
    \global\p@gelist=\expandafter\m@kelist\the\t@mplistX1\fi\fi
    \let\next=\re@dfns\fi\next}
\def\symbols#1{\symb@ls={#1}\s@mb@ls=\symb@ls} 
\def\bigsymbol{\textstyle}
\symbols{\bigsymbol\ast,\dagger,\ddagger,\sharp,\flat,\natural,\star}
\def\ftnumbers{\ftn@mberstrue} \def\ftsymbols{\ftn@mbersfalse}
\def\paginal{\byp@ge} \def\resetftnumbers{\ftn@mber=1}
\def\ftnote#1{\defm@rk\expandafter\expandafter\expandafter\footnote
    \expandafter\b@@st\m@rk{#1}}

\long\def\jump#1\endjump{}
\def\ssum{\mathop{\lower .1em\hbox{$\textstyle\Sigma$}}\nolimits}

\def\qed{\nobreak\kern 1em \vrule height .5em width .5em depth 0em}
\def\newneq{\hbox{\rlap{\hbox to 1\wd9{\hss$=$\hss}}\raise .1em 
   \hbox to 1\wd9{\hss$\scriptscriptstyle/$\hss}}}
\def\subsetne{\setbox9 = \hbox{$\subset$}\mathrel{\hbox{\rlap
   {\lower .4em \newneq}\raise .13em \hbox{$\subset$}}}}
\def\supsetne{\setbox9 = \hbox{$\subset$}\mathrel{\hbox{\rlap
   {\lower .4em \newneq}\raise .13em \hbox{$\supset$}}}}

\def\vbar{\mathchoice{\vrule height6.3ptdepth-.5ptwidth.8pt\kern-.8pt}
   {\vrule height6.3ptdepth-.5ptwidth.8pt\kern-.8pt}
   {\vrule height4.1ptdepth-.35ptwidth.6pt\kern-.6pt}
   {\vrule height3.1ptdepth-.25ptwidth.5pt\kern-.5pt}}
\def\f@dge{\mathchoice{}{}{\mkern.5mu}{\mkern.8mu}}
\def\b@c#1#2{{\rm \mkern#2mu\vbar\mkern-#2mu#1}}
\def\b@b#1{{\rm I\mkern-3.5mu #1}}
\def\b@a#1#2{{\rm #1\mkern-#2mu\f@dge #1}}
\def\bb#1{{\count4=`#1 \advance\count4by-64 \ifcase\count4\or\b@a A{11.5}\or
   \b@b B\or\b@c C{5}\or\b@b D\or\b@b E\or\b@b F \or\b@c G{5}\or\b@b H\or
   \b@b I\or\b@c J{3}\or\b@b K\or\b@b L \or\b@b M\or\b@b N\or\b@c O{5} \or
   \b@b P\or\b@c Q{5}\or\b@b R\or\b@a S{8}\or\b@a T{10.5}\or\b@c U{5}\or
   \b@a V{12}\or\b@a W{16.5}\or\b@a X{11}\or\b@a Y{11.7}\or\b@a Z{7.5}\fi}}

\catcode`\X=11 \catcode`\@=12

\expandafter\ifx\csname citeadd.tex\endcsname\relax
\expandafter\gdef\csname citeadd.tex\endcsname{}
\else \message{Hey!  Apparently you were trying to
\string twice.   This does not make sense.} 
\errmessage{Please edit your file (probably \jobname.tex) and remove
any duplicate ``\string\input'' lines} \fi

\def\sciteu{\sciteerror{undefined}}
\def\sciteuphantom{\complainaboutcitation{undefined}}

\def\sciteerror#1#2{{\mathortextbf{\scite{#2}}}\complainaboutcitation{#1}{#2}}
\def\mathortextbf#1{\hbox{\bf #1}}
\def\complainaboutcitation#1#2{%
\vadjust{\line{\llap{---$\!\!>$ }\qquad scite$\{$#2$\}$ #1\hfil}}}

\sectno=-1   
\localtags
\ifx\shlhetal\undefinedcontrolsequence\let\shlhetal\relax\fi
\NoBlackBoxes
\define\mr{\medskip\roster}
\define\sn{\smallskip\noindent}
\define\mn{\medskip\noindent}
\define\bn{\bigskip\noindent}
\define\ub{\underbar}
\define\wilog{\text{without loss of generality}}
\define\ermn{\endroster\medskip\noindent}

\define\dbcu{\dsize\bigcup}
\define\nl{\newline}
\documentstyle {amsppt}
\topmatter
\title {pcf Theory: Applications} \endtitle
\rightheadtext{$pcf$ Theory Applications}
\author {Saharon Shelah \thanks {\null\newline
Partially supported by the basic research fund, Israeli Academy \null\newline
This version came from sections of Sh580 \newline
I thank Alice Leonhardt for the excellent typing \newline
Latest Revision -  98/Apr/10 \null\newline
Pub. No. 589 \null\newline
Saharon references to 513:p.32,34, see \S2 4.8} \endthanks} \endauthor
\affil {Institute of Mathematics \\
The Hebrew University \\
Jerusalem, Israel
\medskip
Rutgers University \\
Department of Mathematics \\
New Brunswick, NJ  USA} \endaffil
\abstract{We deal with several pcf problems; we characterize another version
of exponentiation: number of $\kappa$-branches in a tree with $\lambda$
nodes, deal with existence of independent sets in stable theories, possible
cardinalities of ultraproducts and the depth of ultraproducts of Boolean 
Algebras.
Also we give cardinal invariants for each $\lambda$ with a pcf restriction and
investigate further $T_D(f)$.  The sections can be read independently.} 
\endabstract
\endtopmatter
\document  

\expandafter\ifx\csname alice2jlem.tex\endcsname\relax
  \expandafter\gdef\csname alice2jlem.tex\endcsname{}
\else \message{Hey!  Apparently you were trying to
\string  twice.   This does not make sense.}
\errmessage{Please edit your file (probably \jobname.tex) and remove
any duplicate ``\string\input'' lines} \fi

\expandafter\ifx\csname bib4plain.tex\endcsname\relax
  \expandafter\gdef\csname bib4plain.tex\endcsname{}
\else \message{Hey!  Apparently you were trying to \string twice.   This does not make sense.}
\errmessage{Please edit your file (probably \jobname.tex) and remove
any duplicate ``\string\input'' lines} \fi

\def\renewcommand{\newcommand}	       
\edef\cite{\the\catcode`@}%
\catcode`@ = 11
\let\@oldatcatcode = \cite
\chardef\@letter = 11
\chardef\@other = 12
%
%
%
%
\def\@innerdef#1#2{\edef#1{\expandafter\noexpand\csname #2\endcsname}}%
%
%
\@innerdef\@innernewcount{newcount}%
\@innerdef\@innernewdimen{newdimen}%
\@innerdef\@innernewif{newif}%
\@innerdef\@innernewwrite{newwrite}%
%
%
%
\def\@gobble#1{}%
%
%
%
\ifx\inputlineno\@undefined
   \let\@linenumber = \empty 
\else
   \def\@linenumber{\the\inputlineno:\space}%
\fi
%
%
%
\def\@futurenonspacelet#1{\def\cs{#1}%
   \afterassignment\@stepone\let\@nexttoken=
}%
\begingroup 
\def\\{\global\let\@stoken= }%
\\ 
\endgroup
\def\@stepone{\expandafter\futurelet\cs\@steptwo}%
\def\@steptwo{\expandafter\ifx\cs\@stoken\let\@@next=\@stepthree
   \else\let\@@next=\@nexttoken\fi \@@next}%
\def\@stepthree{\afterassignment\@stepone\let\@@next= }%
%
%
%
\def\@getoptionalarg#1{%
   \let\@optionaltemp = #1%
   \let\@optionalnext = \relax
   \@futurenonspacelet\@optionalnext\@bracketcheck
}%
%
%
\def\@bracketcheck{%
   \ifx [\@optionalnext
      \expandafter\@@getoptionalarg
   \else
      \let\@optionalarg = \empty
      \expandafter\@optionaltemp
   \fi
}%
\def\@@getoptionalarg[#1]{%
   \def\@optionalarg{#1}%
   \@optionaltemp
}%
%
%
%
\def\@nnil{\@nil}%
\def\@fornoop#1\@@#2#3{}%
\def\@for#1:=#2\do#3{%
   \edef\@fortmp{#2}%
   \ifx\@fortmp\empty \else
      \expandafter\@forloop#2,\@nil,\@nil\@@#1{#3}%
   \fi
}%
\def\@forloop#1,#2,#3\@@#4#5{\def#4{#1}\ifx #4\@nnil \else
       #5\def#4{#2}\ifx #4\@nnil \else#5\@iforloop #3\@@#4{#5}\fi\fi
}%
\def\@iforloop#1,#2\@@#3#4{\def#3{#1}\ifx #3\@nnil
       \let\@nextwhile=\@fornoop \else
      #4\relax\let\@nextwhile=\@iforloop\fi\@nextwhile#2\@@#3{#4}%
}%
%
%
%
\@innernewif\if@fileexists
\def\@testfileexistence{\@getoptionalarg\@finishtestfileexistence}%
\def\@finishtestfileexistence#1{%
   \begingroup
      \def\extension{#1}%
      \immediate\openin0 =
         \ifx\@optionalarg\empty\jobname\else\@optionalarg\fi
         \ifx\extension\empty \else .#1\fi
         \space
      \ifeof 0
         \global\@fileexistsfalse
      \else
         \global\@fileexiststrue
      \fi
      \immediate\closein0
   \endgroup
}%
%
%
%
%
\def\bibliographystyle#1{%
   \@readauxfile
   \@writeaux{\string\bibstyle{#1}}%
}%
\let\bibstyle = \@gobble
%
%
\let\bblfilebasename = \jobname
\def\bibliography#1{%
   \@readauxfile
   \@writeaux{\string\bibdata{#1}}%
   \@testfileexistence[\bblfilebasename]{bbl}%
   \if@fileexists
      \nobreak
      \@readbblfile
   \fi
}%
\let\bibdata = \@gobble
%
%
\def\nocite#1{%
   \@readauxfile
   \@writeaux{\string\citation{#1}}%
}%
\@innernewif\if@notfirstcitation
%
%
\def\cite{\@getoptionalarg\@cite}%
%
%
\def\@cite#1{%
   \let\@citenotetext = \@optionalarg
   \printcitestart
   \nocite{#1}%
   \@notfirstcitationfalse
   \@for \@citation :=#1\do
   {%
      \expandafter\@onecitation\@citation\@@
   }%
   \ifx\empty\@citenotetext\else
      \printcitenote{\@citenotetext}%
   \fi
   \printcitefinish
}%
\def\@onecitation#1\@@{%
   \if@notfirstcitation
      \printbetweencitations
   \fi
   \expandafter \ifx \csname\@citelabel{#1}\endcsname \relax
      \if@citewarning
         \message{\@linenumber Undefined citation `#1'.}%
      \fi
      \expandafter\gdef\csname\@citelabel{#1}\endcsname{%
\strut
\vadjust{\vskip-\dp\strutbox
\vbox to 0pt{\vss\parindent0cm \leftskip=\hsize 
\advance\leftskip3mm
\advance\hsize 4cm\strut\openup-4pt 
\rightskip 0cm plus 1cm minus 0.5cm ?  #1 ?\strut}}
         {\tt
            \escapechar = -1
            \nobreak\hskip0pt
            \expandafter\string\csname#1\endcsname
            \nobreak\hskip0pt
         }%
      }%
   \fi
   \csname\@citelabel{#1}\endcsname
   \@notfirstcitationtrue
}%
%
%
\def\@citelabel#1{b@#1}%
%
%
\def\@citedef#1#2{\expandafter\gdef\csname\@citelabel{#1}\endcsname{#2}}%
%
%
%
\def\@readbblfile{%
   \ifx\@itemnum\@undefined
      \@innernewcount\@itemnum
   \fi
   \begingroup
      \def\begin##1##2{%
         \setbox0 = \hbox{\biblabelcontents{##2}}%
         \biblabelwidth = \wd0
      }%
      \def\end##1{}
      %
      %
      \@itemnum = 0
      \def\bibitem{\@getoptionalarg\@bibitem}%
      \def\@bibitem{%
         \ifx\@optionalarg\empty
            \expandafter\@numberedbibitem
         \else
            \expandafter\@alphabibitem
         \fi
      }%
      \def\@alphabibitem##1{%
         \expandafter \xdef\csname\@citelabel{##1}\endcsname {\@optionalarg}%
         \ifx\biblabelprecontents\@undefined
            \let\biblabelprecontents = \relax
         \fi
         \ifx\biblabelpostcontents\@undefined
            \let\biblabelpostcontents = \hss
         \fi
         \@finishbibitem{##1}%
      }%
      \def\@numberedbibitem##1{%
         \advance\@itemnum by 1
         \expandafter \xdef\csname\@citelabel{##1}\endcsname{\number\@itemnum}%
         \ifx\biblabelprecontents\@undefined
            \let\biblabelprecontents = \hss
         \fi
         \ifx\biblabelpostcontents\@undefined
            \let\biblabelpostcontents = \relax
         \fi
         \@finishbibitem{##1}%
      }%
      \def\@finishbibitem##1{%
         \biblabelprint{\csname\@citelabel{##1}\endcsname}%
         \@writeaux{\string\@citedef{##1}{\csname\@citelabel{##1}\endcsname}}%
         \ignorespaces
      }%
      %
      %
      \let\em = \bblem
      \let\newblock = \bblnewblock
      \let\sc = \bblsc
      \frenchspacing
      \clubpenalty = 4000 \widowpenalty = 4000
      \tolerance = 10000 \hfuzz = .5pt
      \everypar = {\hangindent = \biblabelwidth
                      \advance\hangindent by \biblabelextraspace}%
      \bblrm
      \parskip = 1.5ex plus .5ex minus .5ex
      \biblabelextraspace = .5em
      \bblhook
      \input \bblfilebasename.bbl
   \endgroup
}%
%
%
\@innernewdimen\biblabelwidth
\@innernewdimen\biblabelextraspace
%
%
%
\def\biblabelprint#1{%
   \noindent
   \hbox to \biblabelwidth{%
      \biblabelprecontents
      \biblabelcontents{#1}%
      \biblabelpostcontents
   }%
   \kern\biblabelextraspace
}%
%
%
%
\def\biblabelcontents#1{{\bblrm [#1]}}%
%
%
\def\bblrm{\rm}%
%
%
\def\bblem{\it}%
%
%
\def\bblsc{\ifx\@scfont\@undefined
              \font\@scfont = cmcsc10
           \fi
           \@scfont
}%
%
%
\def\bblnewblock{\hskip .11em plus .33em minus .07em }%
%
%
\let\bblhook = \empty
%
%
%
\def\printcitestart{[}
\def\printcitefinish{]}
\def\printbetweencitations{, }
\def\printcitenote#1{, #1}
%
%
%
\let\citation = \@gobble
%
%
%
\@innernewcount\@numparams
%
%
\def\newcommand#1{%
   \def\@commandname{#1}%
   \@getoptionalarg\@continuenewcommand
}%
%
%
\def\@continuenewcommand{%
   \@numparams = \ifx\@optionalarg\empty 0\else\@optionalarg \fi \relax
   \@newcommand
}%
%
%
\def\@newcommand#1{%
   \def\@startdef{\expandafter\edef\@commandname}%
   \ifnum\@numparams=0
      \let\@paramdef = \empty
   \else
      \ifnum\@numparams>9
         \errmessage{\the\@numparams\space is too many parameters}%
      \else
         \ifnum\@numparams<0
            \errmessage{\the\@numparams\space is too few parameters}%
         \else
            \edef\@paramdef{%
               \ifcase\@numparams
                  \empty  No arguments.
               \or ####1%
               \or ####1####2%
               \or ####1####2####3%
               \or ####1####2####3####4%
               \or ####1####2####3####4####5%
               \or ####1####2####3####4####5####6%
               \or ####1####2####3####4####5####6####7%
               \or ####1####2####3####4####5####6####7####8%
               \or ####1####2####3####4####5####6####7####8####9%
               \fi
            }%
         \fi
      \fi
   \fi
   \expandafter\@startdef\@paramdef{#1}%
}%
%
%
%
%
\def\@readauxfile{%
   \if@auxfiledone \else 
      \global\@auxfiledonetrue
      \@testfileexistence{aux}%
      \if@fileexists
         \begingroup
            \endlinechar = -1
            \catcode`@ = 11
            \input \jobname.aux
         \endgroup
      \else
         \message{\@undefinedmessage}%
         \global\@citewarningfalse
      \fi
      \immediate\openout\@auxfile = \jobname.aux
   \fi
}%
%
%
\newif\if@auxfiledone
\ifx\noauxfile\@undefined \else \@auxfiledonetrue\fi
%
%
%
%
\@innernewwrite\@auxfile
\def\@writeaux#1{\ifx\noauxfile\@undefined \write\@auxfile{#1}\fi}%
%
%
%
\ifx\@undefinedmessage\@undefined
   \def\@undefinedmessage{No .aux file; I won't give you warnings about
                          undefined citations.}%
\fi
%
%
\@innernewif\if@citewarning
\ifx\noauxfile\@undefined \@citewarningtrue\fi
%
%
%
\catcode`@ = \@oldatcatcode


\def\widestnumber#1#2{}

\def\rm{\fam0 \tenrm}

\def\fakesubhead#1\endsubhead{\bigskip\noindent{\bf#1}\par}


%
%
%

%

\font\textrsfs=rsfs10
\font\scriptrsfs=rsfs7
\font\scriptscriptrsfs=rsfs5

\newfam\rsfsfam
\textfont\rsfsfam=\textrsfs
\scriptfont\rsfsfam=\scriptrsfs
\scriptscriptfont\rsfsfam=\scriptscriptrsfs

\edef\oldcatcodeofat{\the\catcode`\@}
\catcode`\@11

\def\Cal@@#1{\noaccents@ \fam \rsfsfam #1}

\catcode`\@\oldcatcodeofat

\newpage

\head {Annotated Content} \endhead  \resetall \bigskip

\noindent
\S1 $\quad$ $T_D$ via true cofinalities
\roster
\item "{{}}"  [Assume $D$ is a filter on $\kappa,\mu = \text{ cf}(\mu) >
2^\kappa,f \in {}^\kappa\text{Ord}$, and: $D$ is  \newline
$\aleph_1$-complete or
$(\forall \sigma < \mu)(\sigma^{\aleph_0} < \mu)$.  We prove that if
$T_D(f) \ge \mu$ (i.e. there are $f_\alpha <_D f$ for $\alpha < \mu$ such that
$f_\alpha \ne_D f_\beta$ for $\alpha < \beta < \mu)$ \underbar{then} for some
$A \in D^+$ and regular $\lambda_i \in (2^\kappa,f(i)]$ we have: $\mu$ is the
true cofinality of $\dsize \prod_{i < \kappa} \lambda_i/(D+A)$].
\endroster
\medskip

\noindent
\S2 $\quad$  The tree revised power
\roster
\item "{{}}"  [We characterize by pcf more natural cardinal functions.
The main one is $\lambda^{\kappa,\text{tr}}$, the supremum on the number of
$\kappa$-branches of trees with $\lambda$ nodes, where $\kappa$ is regular
uncountable.  If $\lambda > \kappa^{\kappa,\text{tr}}$ it is the supremum on
max pcf$\{\theta_\zeta:\zeta < \kappa\}$ for an increasing sequence 
$\langle \theta_\zeta:
\zeta < \kappa \rangle$ of regular cardinals with $\zeta < \kappa \Rightarrow
\lambda \ge \text{ max pcf}\{\theta_\varepsilon:\varepsilon < \zeta\}$].
\endroster 
\medskip

\noindent
\S3 $\quad$ On the depth behaviour for ultraproducts
\roster
\item "{{}}"  [We deal with a problem of Monk on the depth of ultraproducts
of Boolean algebras; this continues \cite[\S3]{Sh:506}.  We try to
characterize for a filter $D$ on $\kappa$ and $\lambda_i = \text{ cf}
(\lambda_i) > 2^\kappa$, and $\mu = \text{ cf}(\mu)$, when does
$(\forall i < \kappa)[\lambda_i \le \text{Depth}^+(B_i)] \Rightarrow
\mu < \text{ Depth}^+(\dsize \prod_{i < \kappa} B_i/D)$ 
(where Depth$^+(B) = \cup\{\mu^+: \text{in }
B$ there is an increasing sequence of length $\mu\}$).  When $D$ is
$\aleph_1$-complete or $(\forall \sigma < \mu)[\sigma^{\aleph_0} < \mu]$ the 
characterization is reasonable: for some $A \in D^+$ and $\lambda'_i =
\text{ cf}(\lambda'_i) < \lambda_i$ we have $\mu = \text{ tcf } \dsize
\prod_{i < \kappa} \lambda'_i/(D+A)$.  We then proceed to look at
Depth$^{(+)}_h$ (closing under homomorphic images), and with more work
succeed.  We use results from \S1].
\endroster
\medskip

\noindent
\S4 $\quad$ On existence of independent sets for stable theories
\roster
\item "{{}}"  [Bays [Bays Ph.D.] has continued work in 
\cite{Sh:c} on existence of
independent sets (in the sense of non-forking) for stable theories. \newline
We connect those problems to pcf and shed some light.  Note that the
combinatorial Claim \scite{4.1} continues \cite[\S3]{Sh:430}].
\endroster
\medskip

\noindent
\S5 $\quad$ Cardinal invariants for general cardinals: 
restriction on the depth
\roster
\item "{{}}"  [We show that some (natural) cardinal invariants defined for any
regular $\lambda(> \aleph_0)$, as functions of $\lambda$ satisfies 
inequalities coming from pcf (more accurately norms for $\aleph_1$-complete
filters).  They are variants of depth, supremum of length of sequences from
${}^\lambda \lambda$, increasing in suitable sense and also the supremum on
$\lambda$-MAD families.  Constrast this with Cummings Shelah \nl
\cite{CuSh:541}.  Also we connect pcf and the ideal $I[\lambda]$; see
\scite{5.16}].
\endroster
\medskip

\noindent
\S6 $\quad$ The class of cardinal ultraproducts mod $D$
\roster
\item "{{}}"  [Let $D$ be an ultrafilter on $\kappa$ and let \newline
reg$(D) = \text{ Min}\{\theta:D \text{ is not } \theta \text{ regular}\}$,
so reg$(D)$ is regular itself.  We prove that if $\mu = \mu^{\text{reg}
(\theta)} + 2^\kappa$ then $\mu$ can be represented as $| \dsize 
\prod_{i < \kappa} \lambda_i/D|$, and for suitable $\mu$'s get $\mu$-like
such ultraproducts].
\endroster
\newpage

\head {\S1 $T_D$ via true cofinality} \endhead  \resetall 
\bigskip

\noindent
We improve here results of \cite[\S3]{Sh:506}.  See more related things in
\S6.  Our main result is \scite{1.4}.
\proclaim{\stag{1.1} Claim}  Assume
\medskip
\roster
\item "{$(a)$}"  $J$ is an $\aleph_1$-complete ideal on $\kappa$
\sn
\item "{$(b)$}"  $f \in {}^\kappa\text{Ord}$, each $f(i)$ an infinite ordinal
\sn
\item "{$(c)$}"  $T^2_J(f) \ge \lambda = \text{ cf}(\lambda) > \mu \ge \kappa$
(see \scite{1.1A}(1) below)
\sn
\item "{$(d)$}"  $\mu = 2^\kappa$, or at least
{\roster
\itemitem{ $(d)^-(i)$ }  if ${\frak a} \subseteq \text{ Reg}$, and \newline
$(\forall \theta \in {\frak a})(\mu \le \theta < \lambda \and \mu \le 
\theta < \underset{i < \kappa} {}\to \sup f(i))$  \nl
and $|{\frak a}| \le \kappa$, \ub{then} $|\text{pcf}({\frak a})| \le \mu$
\sn
\itemitem{ $(ii)$ }  $|\mu^\kappa/J| < \lambda$
\sn
\itemitem{ $(iii)$ }  $2^\kappa < \lambda$.
\endroster}
\ermn
\underbar{Then} for some $A \in J^+$ and $\bar \lambda = \langle \lambda_i:
i \in A \rangle$ such that $\mu \le \lambda_i = \text{ cf}(\lambda_i) \le 
f(i)$ we have
$\dsize \prod_{i \in A} \lambda_i/(J \restriction A)$ has true cofinality
$\lambda$.
\endproclaim
\bigskip

\remark{\stag{1.1A} Remark}  1) Remember 
$T^2_J(f) = \text{ Min}\{|F|:F \subseteq
\dsize \prod_{i < \kappa} f(i)$ and for every \newline
$g \in \dsize \prod_{i < \kappa} f(i)$ for some $g' \in F$ 
we have $\neg(g \ne_J g')\}$.  The inverse of the claim is immediate. \nl
2)  If $A_1 = \{i < \kappa:f(i) \ge \lambda\} \in J^+$ then the conclusion
is immediate, with $\lambda_i = \lambda$. \newline
3) Note if $A_2 = \{ i < \kappa:f(i) < (2^\kappa)^+\} \in J^+$ then
$T^2_J(f) \le 2^\kappa$.  If in addition $\kappa \backslash A_2 \in J$ then
any $\lambda$ satisfying the conclusion satisfies 
$\lambda \le 2^\kappa$. \newline
4) We can omit the assumption clause $(d)^-(iii)$ and weaken 
(here and in \scite{2.5}) the assumption ``$|\mu^\kappa/J| < \lambda$" 
(in clause (d)$^-$) we can just ask:
\medskip
\roster
\item "{$\bigoplus_{J,\mu,\lambda}$}"  there is 
$F \subseteq {}^\kappa \mu$ of cardinality
$< \lambda$ such that for every $g \in {}^\kappa \mu$ we can find 
$F' \subseteq F$ of cardinality $\le \mu$ such that for every $A \in J^+$ 
for some $f \in F'$ we have $\{i \in A:g(i) = f(i)\} \in J^+$, or even
\smallskip
\noindent
\item  "{$\bigoplus^-_{J,\mu,\lambda}$}"  we require the above only for all
$g \in G$, where $G \subseteq {}^\kappa \mu$ has cardinality $< \lambda$ and:
if $\langle \theta_i:i < \kappa \rangle$ is a sequence of regulars in
$[\aleph_0,\mu]$ and $g' \in \dsize \prod_{i < \kappa} \theta_i$ then for some
$g'' \in G$ we have $g' <_J g'' <_J \langle \theta_i:i < \kappa \rangle$.
\endroster
\medskip

\noindent
Considering $(d)^-(iii)$ in the proof we weaken $g_n \restriction A \in N$
for some $g',A' \subseteq \kappa$ from $g_n \restriction A =_J g' \restriction
A'$. \nl
5) Also in \scite{1.4} and \scite{1.5} we can replace the 
assumption $\lambda > 2^\kappa$
by the existence of $\mu$ satisfying $\lambda > \mu \ge \kappa$ such that
$(d)^-$ as weakened above holds. \newline
6) Note that we do not ask $(\forall \alpha < \lambda)
[|\alpha|^{< \text{ reg}(J)} < \lambda]$. \nl
7) Of course, we can apply the claim to $J \restriction A$ for every 
$A \in J^+$ hence $\{A/J:A \in J^+$, and for some $\bar \lambda = \langle
\lambda_i:i \in A \rangle$ such that $\mu \le \lambda_i = \text{ cf}
(\lambda_i) \le f(i)$ we have $\dsize \prod_{i \in A} \lambda_i/
(J \restriction A)$ has true cofinality $\lambda\}$ is dense in the Boolean
Algebra ${\Cal P}(\kappa)/J$.
\endremark
\bigskip

\remark{\stag{1.2} Remark}  The changes in the proof required for 
weakening in \scite{1.1} the clause $|\mu^\kappa/J| < \lambda$ to 
$\bigoplus^-_{J,\mu,\lambda}$ from \scite{1.1A}(4) are as follows.

As $J,\mu,\lambda \in N$ there are $F \subseteq {}^\kappa \mu,G \subseteq
{}^\kappa \mu$ as required in $\bigoplus^-_{J,\mu,\lambda}$ belonging to $N$
(hence $\subseteq N$).
After choosing $g^{n,1}$ and $B_n$ apply the assumption on $G$ to
$g^{n,3} \in {}^\kappa \mu$ when $g^{n,3} \restriction B_n = (g^{n,2} 
\restriction B_n)$ and $g^{n,3} \restriction (\kappa \backslash B_n)$ is
constantly zero and $\bar \theta = \langle \theta_i:i < \kappa \rangle$ where
$\theta_i = \text{ cf}(g_n(i))$ if $i \in B_n$ and $\theta_i = \aleph_0$ if
$i \in \kappa \backslash B_n$.

So we get some $g^{n,4} \in G$ such that $g^{n,3} <_J g^{n,4} <_J \langle
\theta_i:i < \kappa \rangle$.  As $G \in N$, \newline
$|G| < \lambda$ clearly $G \subseteq N$ hence $g^{n,4} \in G$.  
Let $F'_n$ be a subset of $F$ of cardinality $\le
\mu$ such that: for every $A \in J^+$ for some $f \in F'_n$ we have
$\{ i \in A:g^{n,4}(i) = f(i)\} \in J^+$.

Now continue as there but defining $g_{n+1}$ use $g^{n,4}$ instead
$g^{n,3}$ and choose ${\Cal P}^1_{n+1}$ as

$$
\biggl\{ \{ i < \kappa:g^{n,4}(i) = f(i)\}:f \in F'_n \biggr\}.
$$
\medskip

\noindent
The rest is straight.
\endremark
\bigskip

\noindent
Remember
\demo{\stag{1.3} Fact}  Assume 
\mr
\item "{$(a)$}"  $N \prec ({\Cal H}(\chi),\in,<^*_\chi)$ and $\mu < \lambda
< \chi$ and $\{\mu,\lambda\} \in N$,
\sn
\item "{$(b)$}"   $N \cap \lambda$ is an ordinal,
\sn
\item "{$(c)$}"  $i^* \le \mu$, and for $i < i^*$ we have 
${\frak a}_i \subseteq \text{ Reg} \backslash
\mu^+,|{\frak a}_i| \le \mu,\theta_i \in \text{pcf}({\frak a}_i) 
\cap \lambda$ and
$({\frak a}_i,\theta_i) \in N$, and let ${\frak a} = \dsize \bigcup_{i <
i^*} {\frak a}_i$.
\ermn
\ub{Then}
\mr
\item "{$(*)$}"  for every $g \in \Pi{\frak a}$ there is $f$ such that:
{\roster
\itemitem{ $(\alpha)$ }  $g < f \in \Pi {\frak a}$
\sn
\itemitem { $(\beta)$ }  $f \restriction {\frak b}_{\theta_i}[{\frak a}_i]
\in N$, and if $\theta_i = \text{ max pcf}({\frak a}_i)$ we have
$f \restriction {\frak a}_i \in N$.
\endroster}
\endroster
\enddemo
\bigskip

\demo{Proof}  By \cite[Ch.II,3.4]{Sh:g} or \cite[VIII,\S1]{Sh:g}.
\enddemo
\bigskip

\demo{Proof of \scite{1.1}}  Note that assuming $2^\kappa < \lambda$ slightly
simplify the proof, as then we can demand $g_{A,n} = g_n \restriction A$.
Assume toward contradiction that the conclusion
fails.  Without loss of generality $i < \kappa \Rightarrow f(i) > 0$.  
Let $\chi$ be large enough, and let $N$ be an
elementary submodel of $({\Cal H}(\chi),\in,<^*_\chi)$ of cardinality 
$< \lambda$ such that $\{f,\lambda,\mu\}$ belongs to $N$ and $N \cap 
\lambda$ is an ordinal and if we assume only clause $(d)^-$ then
\footnote{note we did not forget to ask $J \in N$, we just want to help
reading this as a proof of \scite{1.3A}, too; for the case $2^{|J|} 
\ge \lambda$ so there $J'$ does not necessarily belong to $N$.}
\mr
\item "{$\boxtimes$}"  for every $f \in {}^\kappa \mu$ for some $g \in
N \cap {}^\kappa \mu$ such that $f = g \text{ mod } J$ (if $J \in N$ this
is immediate).
\ermn
So we shall prove $F =: \left( \dsize \prod_{i < \kappa} f(i) \right) \cap N$
exemplifies that $T^2_J(f) \le |F|(< \lambda)$, thus giving a contradiction

So it suffices to prove
\medskip
\roster
\item "{$(*)$}"  for every $g \in \dsize \prod_{i < \kappa}f(i)$ for some
$g' \in F$ we have $\neg(g \ne_J g')$ i.e. \newline
$\{i < \kappa:g'(i) = g(i)\} \in J^+$.  
\endroster
\medskip

\noindent
Assume $g \in \dsize \prod_{i < \kappa} f(i)$ exemplifies the failure of
$(*)$.

We now define by induction on $n < \omega$ the function $g_n$ and the family
${\Cal P}_n$ such that:
\medskip
\roster
\widestnumber\item{$(viii)$}
\item "{$(i)$}"  $g_0 = f,g_n \in {}^\kappa\text{Ord}$, and
$g \le g_n$
\sn
\item "{$(ii)$}"  $g_{n+1} < g_n \text{ mod } J$
\sn
\item "{$(iii)$}"  ${\Cal P}_n$ is a family of $\le \mu$ members of $J^+$ 
\sn
\item "{$(iv)$}"  if $A \in {\Cal P}_n$ then $g_n \restriction A \in N$ hence
$A \in N$ but if $2^\kappa \ge \lambda$ we just assume that for some
$g_{A,n} \in \dsize \prod_{i \in A} f(i)$ we have $g_{A,n} = g_n \restriction
A \text{ mod } J$ and $g_{A,n} \in N$ hence $A \in N$
\sn
\item "{$(v)$}"  ${\Cal P}_0 = \{ \kappa\}$
\sn
\item "{$(vi)$}"  if $A \in {\Cal P}_n$ and $B \subseteq A$ and $B \in J^+$ 
\underbar{then} for some $A' \in {\Cal P}_{n+1}$ we have \newline
$A' \subseteq A \and A' \cap B \in J^+$
\sn
\item "{$(vii)$}"  $g < g_n \text{ mod } J$
\sn
\item "{$(viii)$}"  $g(i) \le g_n(i)$ and $g(i) < g_n(i) \Rightarrow
g_{n+1}(i) < g_n(i)$ \nl
and $g(i) = g_n(i) \Rightarrow g(i) = g_{n+1}(i)$ \nl
(not necessary for \scite{1.1}).
\endroster
\medskip

\noindent
If we succeed as ``$J \text{ is } \aleph_1$-complete (see assumption (a))" 
then by clause $(ii)$ we get a contradiction as $<_J$ is well founded.  
Also the case $n=0$ is easy by $(i) + (v)$.
\newline
(Note: Clause $(vii)$ holds as $g \in \dsize \prod_{i < \kappa}f(i)$).
So assume we have $g_n,{\Cal P}_n$ and we shall define $g_{n+1},{\Cal P}
_{n+1}$.  In $N$ there is a two-place function $\bold e$, written
$\bold e_\delta(i),\bold e_\delta(i)$ is defined iff
$\delta \in \{\alpha:\alpha \text{ a non-zero ordinal } \le 
\underset {i < \kappa}\to \sup f(i)\}$, and $i < \text{ cf}(\delta)$, and if
$\delta$ is a limit ordinal, then 
$\langle \bold e_\delta(i):i < \text{ cf}(\delta) \rangle$ is strictly
increasing with limit $\delta$ and $\bold e_{\alpha + 1}(0) = \alpha$, of
course, Dom$(\bold e_{\alpha +1}) = \{0\}$. \newline
We also know by assumption $(d)$ or $(d)^-(i)$ that
\medskip
\roster
\item "{$\bigotimes$}"  for every $A \in {\Cal P}_n$ we have, letting
${\frak a}^n_A =: \{\text{cf}(g_{A,n}(i)):i \in A\} \backslash \mu^+$, the set
pcf$({\frak a}^n_A)$ has at most $\mu$ members.
\endroster
\medskip

So ${\Cal Y} =: \{(A,{\frak a}^n_A,\theta):A \in {\Cal P}_n \text{ and } 
\theta \in \lambda \cap \text{ pcf}({\frak a}^n_A)\}$ has at most
$|{\Cal P}_n| \times \mu \le \mu \times \mu = \mu$ members (as $|{\Cal P}_n|
\le \mu$ and $|\text{pcf }{\frak a}^n_A| \le \mu$ by $\bigotimes$ above) 
so let $\{(A^n_\varepsilon,{\frak a}^n_\varepsilon,
\theta^n_\varepsilon):\varepsilon < \varepsilon^*_n\}$ list them with
$\varepsilon^*_n \le \mu$.  Clearly ${\frak a}^n_\varepsilon \in N$ 
(as $g_{A,n} \restriction A^n_\varepsilon \in N$ hence 
${\frak a}^n_\varepsilon \in N$ but $\mu + 1 \subseteq N,|\text{pcf}
({\frak a}^n_\varepsilon)| \le \mu$ so ${\Cal Y} \subseteq N$).  
For each $\varepsilon < \varepsilon^*_n$ we define 
$h^n_\varepsilon \in \Pi {\frak a}^n_\varepsilon$ by:

$$
\align
h^n_\varepsilon(\theta) = \text{ Min} \biggl\{ \zeta < \theta:&\text{ if }
i \in A^n_\varepsilon,g(i) < g_n(i), \text{ and} \\
  &\,\theta = \text{ cf}(g_n(i)) \text{ then } g(i) < \bold e_{g_n(i)}(\zeta)
\biggr\}.
\endalign
$$

\noindent
[Why is $h^n_\varepsilon$ well defined?  The number of possible $i$'s is
$\le |A^n_\varepsilon| \le \kappa \le \mu$, for each $i$ (satisfying $i \in
A^n_\varepsilon,g(i) < g_n(i)$ and cf$(g_n(i)) > \mu)$, every 
$\zeta < \theta$ large enough is OK as $\langle \bold e_{g_n(i)}(\zeta):
\zeta < \theta \rangle$ 
is increasing continuous with limit $g_n(i)$).  Lastly, $\theta =
\text{ cf}(\theta) > \mu$ (by the choice of ${\frak a}^n_\varepsilon$) so
all the demands together hold for every large enough $\zeta < \theta$].
\newline
Let ${\frak a}_n = \dsize \bigcup_{\varepsilon < \varepsilon^*_n}
{\frak a}^n_\varepsilon$ and let $h_n \in \Pi {\frak a}_n$ be defined by
$h_n(\theta) =
\sup\{h^n_\varepsilon(\theta):\varepsilon < \varepsilon^*_n \text{ and }
\theta \in {\frak a}^n_\varepsilon\}$, it is well defined by the argument
above.  So by \scite{1.3} there is a function 
$g^{n,1} \in \Pi {\frak a}_n$ such that:
\medskip
\roster
\item "{$(\alpha)$}"  $h_n < g^{n,1}$ 
\sn
\item "{$(\beta)$}"   $g^{n,1} \restriction {\frak b}_{\theta^n_\varepsilon}
[{\frak a}^n_\varepsilon] \in N$ (and $\theta^n_\varepsilon = 
\text{ max pcf}({\frak a}^n_\varepsilon) \Rightarrow 
{\frak b}_{\theta^n_\varepsilon}[{\frak a}^n_\varepsilon] = 
{\frak a}^n_\varepsilon$).
\endroster
\medskip

\noindent
Also we can define $g^{n,2} \in {}^\kappa \text{Ord}$ by:

$$
g^{n,2}(i) = \text{ Min}\{ \zeta < \text{ cf}(g_n(i)):
\bold e_{g_n(i)}(\zeta) \ge g(i)\}.
$$

\noindent
So letting $B_n = \{i:1 \le \text{ cf}(g_n(i)) \le \mu\}$ clearly
$g^{n,2} \restriction B_n \in {}^{(B_n)}\mu$.  Now if assumption (d) holds,
then $\mu^\kappa /J < \lambda$, hence $\mu^\kappa \subseteq N$ so we can find 
$g^{n,3} \in N$ such that $g^{n,2} = g^{n,3}$ mod $(J + (\kappa \backslash
B_n))$; if assumption (d) fails we still can get such $g^{n,3}$ by $\boxtimes$
above.  Lastly, we define $g_{n+1} \in {}^\kappa\text{Ord}$:

$$
g_{n+1}(i) = \cases \bold e_{g_n(i)} \left(g^{n,1}(\text{cf}(g_n(i)))
\right) \quad &\text{ \underbar{if} } \quad \text{ cf}(g_n(i)) > \mu
\text{ and } g_n(i) > g(i) \\
  \bold e_{g_n(i)} \left( g^{n,3}(\text{cf}(g_n(i))) \right) \quad
&\text{ \underbar{if} } \quad \text{ cf}(g_n(i)) \in [1,\mu] 
\text{ and } g_n(i) > g(i) \\
g_n(i) &\text{ \underbar{if} } \quad \,\,g(i) = g_n(i)
\endcases
$$

\noindent
and ${\Cal P}_{n+1} = ({\Cal P}^0_{n+1} \cup {\Cal P}^1_{n+1}) \backslash
J$ where

$$
{\Cal P}^0_{n+1} = \biggl\{ \{i \in A^n_\varepsilon:\text{cf}
(g_{A^n_\varepsilon,n}(i)) \in {\frak b}_{\theta^n_\varepsilon}
[{\frak a}^n_\varepsilon]\}:\varepsilon < \varepsilon^*_n \biggr\}
$$

\noindent
and

$$
{\Cal P}^1_{n+1} = 
\biggl\{ \{i < \kappa:i \in A^* \text{ and  cf}(g_{A^*,n}(i)) \le \mu\}:
A^* \in {\Cal P}_n \biggr\}.
$$

\noindent
(Note: possibly ${\Cal P}_{n+1} \cap J \ne \emptyset$, i.e.
$({\Cal P}^0_{n+1} \cup {\Cal P}^1_{n+1}) \cap J = \emptyset$ 
but this does not make problems).
\enddemo
\bigskip

\noindent
So let us check clauses $(i)-(viii)$. \newline
\underbar{Clause $(i)$}:  Trivial.
\bigskip

\noindent
\underbar{Clause $(ii)$}:  By the definition of $g_{n+1}(i)$ above it is
$< g_n(i)$ except when $g_n(i) = g(i)$, but by clause $(vii)$ we know that
$g < g_n \text{ mod } J$ hence necessarily \newline
$A = \{i < \kappa:g_n(i) = 0\}
\in J$, so really $g_{n+1} < g_n \text{ mod } J$.
\bigskip

\noindent
\underbar{Clause $(iii)$}:  $|{\Cal P}_{n+1}| \le |{\Cal P}_n| +
|\varepsilon^*_n| + \aleph_0$ and 
$|{\Cal P}_n| \le \mu$ by clause $(iii)$ for $n$ (i.e.
the induction hypothesis) and during the construction we show that
$|\varepsilon^*_n| = |{\Cal Y}| \le \mu$.
\bigskip

\noindent
\underbar{Clause $(iv)$}:  let $A \in {\Cal P}_{n+1}$ so we have two cases.
\bigskip

\noindent
\underbar{Case 1}:  $A \in {\Cal P}^0_{n+1}$.

So for some $\varepsilon < \varepsilon^*_n$ we have $(\theta^n_\varepsilon
\in \lambda \cap pcf({\frak a}^n_\varepsilon)$ and)
$A =: \{ i \in A^n_\varepsilon:\text{cf}(g_{A^n_\varepsilon,n}(i)) \in 
{\frak b}_{\theta^n_\varepsilon}[{\frak a}^n_\varepsilon]\}$.
Let $g_{A,n+1} \in \dsize \prod_{i \in A} f(i)$ be defined by $g_{A,n+1}(i) =
e = g_{A^n_\varepsilon,n}(\varepsilon)\bigl( g^{n,1}(\text{cf}
(g_{A^n_\varepsilon,n}(i)))\bigr)$.   \newline
By the choice of $g^{n,1} \in \Pi{\frak a}_n$ we have:

$$
g^{n,1} \restriction {\frak b}_{\theta^n_\varepsilon}
[{\frak a}^n_\varepsilon] \in N.
$$ 

\noindent
Now the set $A$ is definable from $A^n_\varepsilon,g_{A^n_\varepsilon,n}:
{\frak b}_{\theta^n_\varepsilon}[{\frak a}^n_\varepsilon]$ all of which
belongs to $N$ hence $A \in N$.
Also $A^n_\varepsilon \in N$ and clearly $g_{A,n+1}$ is definable
from the functions $g^{n,1} \restriction {\frak b}_{\theta_\varepsilon}
[{\frak a}^n_\varepsilon],g^{n,2},g_{A^n_\varepsilon,n},A^n_\varepsilon$ 
and the function $\bold e$ (see the definition of $g_{n+1}$ by cases), but 
all four are from $N$ so $g_{A,n+1} \in N$.  Lastly, $g_{n+1} \restriction
A \equiv_J g_{A,n}$ as $i \in A \and g_{A^n_\varepsilon,n}(i) = 
g_n(i) \and g_n(i) > g(i) \Rightarrow g_{n+1}(i) = g_{A,n+1}(i)$ and each of
the three assumptions fail only for a st of $i \in A$ which belongs to $J$.
\bigskip

\noindent
\underbar{Case 2}:  $A \in {\Cal P}^1_{n+1}$.

So for some $A^* \in {\Cal P}_n$ we have

$$
A = \{ i < \kappa:i \in A^* \text{ and  cf}(g_{A^*,n}(i)) \le \mu \}.
$$

\noindent
Let $g_{A,n+1}(i) \equiv \bold e_{g_{A,n}}(g^{n,3}(\text{cf}(g_{A^*,n}(i)))$.
Again, $g_{A,n+1} \in N,g_{A,n+1} \equiv_J g_{n+1} \restriction A$.
Looking at the definition of $g_{A,n+1}$, clearly $g_{A,n}$ is
definable from $g^{n,2} \in N,g_{A^*,n}$ and the function $\bold e$, all 
of which belong to $N$.
\bigskip

\noindent
\underbar{Clause $(v)$}:  Holds trivially.
\bigskip

\noindent
\underbar{Clause $(vi)$}:  Assume $A \in {\Cal P}_n$ and $B \subseteq A$
satisfies $B \in J^+$ (so also $A \in J^+$), we have to find $A' \in
{\Cal P}_{n+1}$, such that $A' \subseteq A \and A' \cap B \in J^+$.
\bigskip

\noindent
\underbar{Case 1}:  $B_1 = \{i \in B:\text{cf}(g_{A,n}(i)) \le \mu\} \in J^+$.
\newline
In this case $A' =: \{i \in A:\text{cf}(g_{A,n}(i)) \le \mu\} \in {\Cal P}^1
_{n+1} \subseteq {\Cal P}_{n+1}$ and $A' \cap B \in J^+$ by the assumption
of the case.
\bigskip

\noindent
\underbar{Case 2}:  For some $\varepsilon < \varepsilon^*_n$ we have
$A = A^n_\varepsilon$ and

$$
B_2 = \{i \in B:\text{cf}(g_{A,n}(i)) \in {\frak b}_{\theta^n_\varepsilon}
[{\frak a}^n_\varepsilon]\} \in J^+.
$$

\noindent
In this case $A' =: \{i \in A:\text{cf}(g_{A,n}(i)) \in 
{\frak b}_{\theta^n_\varepsilon}[{\frak a}^n_\varepsilon]\} \in J^+$ 
belongs to ${\Cal P}^1_{n+1} \subseteq {\Cal P}_{n+1}$, is $\subseteq A$ 
and $B \cap A' \in J^+$ by the assumption of the case (remember
$g <_J g_n$).
\bigskip

\noindent
\underbar{Case 3}:  Neither Case 1 nor Case 2. \newline
So $B_3 = B \backslash B_1 \in J^+$ and let $\lambda_i = \text{ cf}
(g_{A,n}(i))$.

We shall show that $\dsize \prod_{i \in B_3} \text{ cf}(g_{A,n}(i))/J$ is 
$\lambda$-directed.  This suffices as letting \newline 
$\lambda_i =: \text{ cf}(g_{A,n}(i)) \in (\mu,f(i))]$, by
\cite[II,\S1]{Sh:g} for some $\lambda'_i = \text{ cf}(\lambda'_i) \le 
\lambda_i$, we have \newline
lim inf$_{J \restriction B_3}\langle \lambda'_i:i \in B_3 \rangle = 
\text{ lim inf}_{J \restriction B_3}\langle \lambda_i:i \in B_2 \rangle$
and \newline
$\lambda = \text{ tcf } \dsize \prod_{i \subseteq B_3} \lambda'_i /(J
\restriction B_3)$ and this shows that the conclusion of \scite{1.1} holds, 
but we are under the assumption it fails so the $\lambda$-directedness really 
suffice.

Now $i \in B \backslash B_1 \Rightarrow \lambda_i = \text{ cf}(g_n(i)) > \mu$;
and if $\dsize \prod_{i \in B_3} \lambda_i/J$ is not $\lambda$-directed, by
\cite{Sh:g},I,\S1 for some $B_4 \subseteq B_3$ and $\theta = \text{ cf}
(\theta) < \lambda$ we have: $B_4 \in J^+$ and $\dsize \prod_{i \in B_4}
\lambda_i/J$ has true cofinality $\theta$.  Hence $\theta \in \text{ pcf}
\{\text{cf}(g_{A,n}(i)):i \in A \text{ and cf}(g_n(i)) > \mu\}$, and as 
$\theta > \mu$, for some $\varepsilon < \varepsilon^*_n$ we have $A = 
A^n_\varepsilon$ and $\theta = \theta^n_\varepsilon$ so 
$A' = \{i \in A:\text{cf}(g_{A,n}(i)) \in {\frak b}_{\theta_\varepsilon}
[{\frak a}^n_\varepsilon]\}$ is as required in case 2
on $B_2$ (note: we could have restricted ourselves to $\theta$'s like that).
\bigskip

\noindent
\underbar{Clause $(vii)$}:  By the choice of $g^{n,1},g^{n,2}$ and $g^n$
clearly $i < \kappa \and g(i) < g_n(i) \Rightarrow g(i) \le g_{n+1}(i)$.
As $g < g_n \text{ mod }D$ it suffices to prove $B =: \{i:g(i) = g_{n+1}
(i)\} \in J$.  If not, we choose by induction on $\ell \le n+1$ a member
$B_\ell$ of ${\Cal P}_\ell$ such that $B_\ell \cap B \in J^+$.  For $\ell = 0$
let $B_\ell = \kappa \in {\Cal P}_0$, for $\ell + 1$ apply clause $(vi)$ for
$\ell$ (even when $\ell = n$ we have just proved it).  So $B_{n+1} \cap B
\in J^+$ and $g_{n+1} \restriction (B_{n+1} \cap B) = g \restriction
(B_{n+1} \cap B)$ hence $\neg(g_{n+1} \restriction B_{n+1} \ne_J g_n
\restriction B_{n+1})$ but $g_{n+1} \restriction B_{n+1} \in N$ so we
have contradicted the choice of $g$ as contradicting $(*)$. 
\bn
\ub{Clause $(viii)$}:  Easy.  \hfill$\square_{\scite{1.1}}$
\bigskip

\proclaim{\stag{1.3A} Claim}  Assume
\medskip
\roster
\widestnumber\item{$(d)^-(ii)$}
\item "{$(a)$}"  $J$ is an ideal on $\kappa$
\sn
\item "{$(b)$}"  $f \in {}^\kappa\text{Ord}$, each $f(i)$ an infinite ordinal
\sn
\item "{$(c)$}"  $T^2_J(f) \ge \lambda = \text{ cf}(\lambda) > \mu > \kappa$
\sn
\item "{$(d)$}"  $\mu = (2^\kappa)^+$ or at least
\sn
\item "{$(d)^-(i)$}"  if ${\frak a} \subseteq \text{ Reg}$, and \newline
$(\forall \theta \in
{\frak a})(\mu \le \theta < \lambda \and \mu \le \theta < f(i))$ 
\newline
and $|{\frak a}| \le \kappa$ then $|\text{pcf}({\frak a})| \le \mu$
\sn
\item "{$(ii)$}"  $|\mu^\kappa/J| < \lambda \vee (\forall g \in {}^\kappa\mu)
[|\Pi g/J| < \lambda]$ and $\mu$ is regular
\sn
\item "{$(e)$}"  $\alpha < \lambda \Rightarrow |\alpha|^{\aleph_0} < \lambda$.
\endroster
\medskip

\underbar{Then} for some $A \in J^+$ and $\bar \lambda = \langle \lambda_i:
i \in A \rangle$ such that $\mu \le \text{ cf}(\lambda_i) = \lambda_i \le
f(i)$ we have $\dsize \prod_{i \in A} \lambda_i/J$ has true cofinality 
$\lambda$.
\endproclaim  
\bigskip

\demo{Proof}  We repeat the proof of \scite{1.1} but we choose $N$ such that
${}^\omega N \subseteq N$, (possible by assumption (e) as $\lambda$ is
regular), and let $F =: (\dsize \prod_{i < \kappa} f(i)) \cap N$.  
If $2^\kappa < \lambda$ then clearly

$$
\align
F = \biggl\{ g \in \dsize \prod_{i < \kappa} f(i):&\text{ for some 
partition } \langle A_n:n < \omega \rangle \text{ of } \kappa \text{ and }\\
  &\,g_n \in N \cap \dsize \prod_{i < \kappa} f(i) \text{ we have }\\
  &\,g = \dsize \bigcup_{n < \omega} (g_n \restriction A_n) \biggr\}.
\endalign
$$

\noindent
Then assume $(*)$ (from the proof of \scite{1.1}) fails and 
$g \in \dsize \prod_{i < \kappa}f(i)$ exemplifies
it and we let $J'$ be the ideal $J' = \{A \subseteq \kappa:g \restriction
A = g' \restriction A \text{ for some } g' \in F\}$.

Clearly $J'$ is $\aleph_1$-complete, $J' \subseteq J$ (as $g$ is a 
counterexample to $(*)$ and the representation of $F$ above) and we 
continue as there getting the conclusion for $J'$ hence for $J$.

If $2^\kappa \ge \lambda$, so
\mr
\item "{$\bigotimes$}"  for $g \in \dsize \prod_{i < \kappa} f(i)$ and
$A \in J^+$ we have $(i) \Leftrightarrow (ii)$ where:
{\roster
\itemitem{ (i) }  there are $g'_n \in F$ for $n < \omega$ such that 
$\{i < \kappa:\dsize \bigvee_{n < \omega} g(i) = g'_n(i)\} \supseteq A
\text{ mod } J$
\sn
\itemitem{ (ii) }  for some $g' \in F$ we have
$\{i < \kappa:g(i) = g'(i)\} \supseteq A \text{ mod } J$
\endroster}
\endroster
\medskip

\noindent
[why?  $\Leftarrow$ is trivial; now $\Rightarrow$ holds as 
$g_n \in N$ also $\langle g_n:n < \omega \rangle \in N$ hence
$\langle \{g_n(i):n < \omega\}:i < \kappa \rangle \in N$ and use
$\omega^\kappa/J \le \mu^\kappa/J < \lambda$ (or just 
$\bigoplus_{J,\mu,\lambda}$ from \scite{1.1A}(4).] \nl
\sn
Let $g \in \dsize \prod_{i < \kappa} f(i)$ be such that 
$f \in N \cap \dsize \prod_{i < \kappa} f(i) \Rightarrow g \ne_J f$.
Now we repeat the proof of \scite{1.1} with our $\kappa,f,\lambda,N,F,g$ this
time using the demands in clause (viii) (i.e. $g(i) \le g_n(i)$).  The proof
does not change except that we do not get a contradiction from $n < \omega
\Rightarrow g_{n+1} <_J g_n$.  However, for each $i < \kappa,\langle g_n(i):
n < \omega \rangle$ is non-increasing (by clause (viii)) hence eventually
constant and by that clause eventually equal to $g(i)$.  So clause (i) of
$\bigotimes$ above holds hence clause (ii) so we are \nl
done.  \hfill$\square_{\scite{1.3A}}$
\enddemo  
\bigskip

\demo{\stag{1.4} Conclusion}  Assume $J$ is an ideal on 
$\kappa,f \in {}^\kappa \text{Ord},i < \kappa \Rightarrow f(i) > 
2^\kappa$, \nl
$\lambda = \text{ cf}(\lambda) > 2^\kappa$, and
\medskip
\roster
\item "{$(*)$}"  $J$ is $\aleph_1$-complete or $(\forall \alpha < \lambda)
(|\alpha|^{\aleph_0} < \lambda)$.
\endroster
\medskip

\noindent
Then $(a) \Leftrightarrow (b) \Leftrightarrow (b)^+ \Leftrightarrow (c)
\Leftrightarrow (c)^+$ where
\medskip
\roster
\widestnumber\item{$(b)^+$}
\item "{$(a)$}"  for some $A \in J^+$ we have $T^2_{J \restriction A}(f
\restriction A) \ge \lambda$
\sn
\item "{$(b)$}"  for some $A \in J^+$ and $\lambda_i = \text{ cf}(\lambda_i)
\in (2^\kappa,f(i)]$ (for $i \in A)$ we have \nl
$\dsize \prod_{i \in A} \lambda_i/(J \restriction A)$ is $\lambda$-directed
\sn
\item "{$(b)^+$}"  like $(b)$ but $\dsize \prod_{i \in A} \lambda_i/(J
\restriction A)$ has true cofinality $\lambda$
\sn
\item "{$(c)$}"  for some $A \in J^+$, and $\bar n = \langle n_i:i < \kappa
\rangle \in {}^\kappa \omega$ and ideal $J^*$ on \newline
$A^* = \dsize \bigcup_{i \in A}(\{i\} \times n_i)$ satisfying

$$
(\forall B \subseteq A)[B \in J \Leftrightarrow \dsize \bigcup_{i \in B}
(\{i\} \times n_i) \in J^*]
$$

\noindent
and regular cardinals $\lambda_{(i,n)} \in (2^\kappa,f(i)]$ we have
$\dsize \prod_{(i,n) \in A^*} \lambda_{(i,n)}/J^*$ is \newline
$\lambda$-directed
\sn
\item "{$(c)^+$}"  as in $(c)$ but $\dsize \prod_{(i,n) \in A^*}
\lambda_{(i,n)}/J^*$ has true cofinality $\lambda$.
\endroster
\enddemo
\bigskip

\demo{Proof}  Clearly $(b)^+ \Rightarrow (b),(b) \Rightarrow (c),(b)^+
\Rightarrow (c)^+$ and $(c)^+ \Rightarrow (c)$.  Also 
$(b) \Rightarrow (b)^+$ by 
\cite[II,1.5B]{Sh:g}, and similarly $(c) \Rightarrow (c)^+$.  Now we prove
$(c) \Rightarrow (a)$; let \newline
$\lambda_i = \text{ max}\{\lambda_{(i,n)}:n < n_i\}$ and let $g_i$ be a
one-to-one function from \newline
$\dsize \prod_{n < n_i} \lambda_{(i,n)}$ into $\lambda_i$ and let 
$\langle f_\alpha:\alpha < \lambda \rangle$ be a $<_{J^*}$-increasing sequence
in \nl
$\dsize \prod_{(i,n) \in A^*} \lambda_{(i,n)}$.  Define $f^*_\alpha \in
\dsize \prod_{i \in A} \lambda_i$ by $f^*_\alpha(i) = g_\alpha \left(
f_\alpha \restriction (\{i\} \times n_i) \right)$.  So if $\alpha < \beta$,
then

$$
\biggl\{ i \in A:f^*_\alpha(i) = f^*_\beta(i) \biggr\} =
\biggl\{ i:\dsize \bigwedge_{n < n_i} f_\alpha((i,n)) = f_\beta(i,n) \biggr\}
$$

\noindent
so by the assumption on $J^*$ we get $f_\alpha \ne_{J^*} f_\beta$ hence
$\{ f^*_\alpha:\alpha < \lambda\}$ is as required in clause $(a)$.

Lastly $(a) \Rightarrow (b)$ by 1.1 (in the case $J$ is $\aleph_1$-complete)
or \scite{1.3A} (in the case $(\forall \alpha < \lambda)(|\alpha|^{\aleph_0} 
< \lambda))$.  We have gotten enough implications to prove the conclusions.
\nl
${{}}$ \hfill$\square_{\scite{1.4}}$
\enddemo
\bigskip

\demo{\stag{1.5} Conclusion}  Let $D$ be an ultrafilter on $\kappa$.  If
$\biggl| \dsize \prod_{i < \kappa} f(i)/D \biggr| \ge \lambda 
= \text{ cf}(\lambda) > 2^\kappa$ and $(\forall \alpha < \lambda)
[|\alpha|^{\aleph_0} < \lambda]$, \underbar{then} for some regular 
$\lambda_i \le f(i)$ (for $i < \kappa$) we have 
$\lambda = \text{ tcf}(\dsize \prod_{i < \kappa} \lambda_i/D)$.
\enddemo
\bigskip

\remark{Remark}  On $|\dsize \prod_{i < \kappa} \lambda_i/D|$, see
\cite[3.9B]{Sh:506}.
\endremark
\newpage

\head {\S2 The tree revised power} \endhead  \resetall
\bigskip

\definition{\stag{2.1} Definition}  For 
$\kappa$ regular and $\lambda \ge \kappa$ let

$$
\lambda^{\kappa,\text{tr}} = \sup\{|\text{lim}_\kappa(T)|:
T \text{ a tree with } \le \lambda \text{ nodes and } \kappa \text{ levels}\}
$$

\noindent
where $\text{lim}_\kappa(T)$ is the set of $\kappa$-branches of $T$; and let
when $\lambda \ge \mu \ge \kappa$ and $\theta \ge 1$

$$
\align
\lambda^{\langle \kappa,\theta \rangle} = \text{ Min} \biggl\{
\mu:&\text{ if } T \text{ is a tree with } \lambda \text{ nodes and } \kappa
\text{ levels}, \\
  &\,\text{\underbar{then} there is } {\Cal P} \in \bigl[ [T]^\theta \bigr]
^\mu \text{ such that} \\
  &\,\eta \in \text{ lim}_\kappa(T) \Rightarrow (\exists A \in {\Cal P})
(\eta \subseteq A) \biggr\}.
\endalign
$$
\medskip
\roster
\item "{${{}}$}"  $\quad \lambda^{\langle \kappa \rangle} = 
\lambda^{\langle \kappa,\kappa \rangle}$.
\endroster
\medskip

\noindent
Recall $[A]^\kappa =: \{B:B \subseteq A \text{ and } |B| = \kappa\}$.
\enddefinition
\bigskip

\remark{\stag{2.1A} Remark}  1) Clearly $\lambda^{\langle \kappa,\theta
\rangle} \le \lambda^{\kappa,\text{tr}} 
\le \lambda^{\langle \kappa,\theta \rangle} + \theta^\kappa$. \nl
2)  If $\kappa = \aleph_0$ then obviously $\lambda^{\kappa,\text{tr}} =
\lambda^\kappa$. \nl
3)  Of course, $\lambda^{\langle \kappa,\theta \rangle} \le \text{ cov}
(\lambda,\theta^+,\kappa^+,\kappa)$ and $\kappa \le \theta \le \sigma \le
\lambda \Rightarrow \lambda^{\langle \kappa,\theta \rangle} \le
\lambda^{\langle \kappa,\sigma \rangle} + \text{cov}(\lambda,\theta^+,
\kappa^+,\kappa)$.  (See \cite{Sh:g} if these concepts are unfamiliar.)
\endremark
\bigskip

\proclaim{\stag{2.2} Theorem}  Let $\kappa$ be regular uncountable
$\, \le \lambda$.  \underbar{Then} the following cardinals are equal:
\medskip
\roster
\item "{$(i)$}"  $\lambda^{\langle \kappa \rangle}$
\sn
\item "{$(ii)$}"  $\lambda + \sup\{\text{\rm max pcf}({\frak a}):
{\frak a} \subseteq \text{ Reg } \cap \lambda \backslash \kappa,
{\frak a} = \{\theta_\zeta:\zeta < \kappa\}$ strictly increasing, \nl

$\qquad \qquad \qquad \qquad$ and if $\xi < \kappa$ then
{\rm max pcf}$(\{\theta_\zeta:\zeta < \xi\}) \le \theta_\xi \le \lambda\}$.
\endroster
\endproclaim
\bigskip

\demo{Proof}  \underbar{First inequality}.  
Cardinal of (i) (i.e. $\lambda^{\langle \kappa \rangle}$)
is $\le$ cardinal of (ii). \newline
Assume not and let $\mu$ be the cardinal from clause (ii) so
$\mu \ge \lambda$. \newline
Let $T$, a tree with $\kappa$ levels and $\lambda$ nodes, exemplify 
$\lambda^{\langle \kappa \rangle} > \mu$.  Without loss of 
generality $T \subseteq {}^{\kappa >}\lambda$ and $<_T = \triangleleft 
\restriction T$. \newline
Let $\{T,\kappa,\lambda,\mu\} \in {\frak B}_n \prec ({\Cal H}(\chi),
\in <^*_\chi),
\mu +1 \subseteq {\frak B}_n,\|{\frak B}_n\| = \mu$ for $n < \omega$,
\newline
${\frak B}_n \in {\frak B}_{n+1},{\frak B}_n \prec {\frak B}_{n+1}$ and let
${\frak B} =: \dsize \bigcup_{n < \omega} {\frak B}_n$.  So 
${\Cal P} =: {\frak B} \cap [T]^{\le \kappa}$ cannot exemplify (i).  
So there is $\eta \in \text{ lim}_\kappa(T)$ such
that $(\forall A \in {\Cal P})[\{\eta \restriction \zeta:\zeta < \kappa\}] 
\nsubseteq A]$.
\medskip

\noindent
We choose by induction on $n,N^0_n,N^1_n$ such that:
\medskip
\roster
\item "{$(a)$}"  $N^0_n \prec N^1_n \prec {\frak B}_{3n}$
\sn
\item "{$(b)$}"  $N^1_0 = \text{ Sk}_{{\frak B}_0}
(\{\zeta:\zeta < \kappa\} \cup \{ \eta \restriction \zeta:\zeta < \kappa\}
\cup \{ \kappa,\mu,\lambda,T\})$ and \newline
$N^0_0 = \text{ Sk}_{{\frak B}_0}(\{\zeta:
\zeta < \kappa\} \cup \{ \kappa,\mu,\lambda,T\})$
\sn
\item "{$(c)$}"  $\|N^\ell_n\| = \kappa$
\sn
\item "{$(d)$}"  $N^0_n \in {\frak B}_{n+1}$
\sn
\item "{$(e)$}"  $N^1_n = \text{ Sk}_{{\frak B}_{3n}}(N^0_n \cup 
\{ \eta \restriction \zeta:\zeta < \kappa\})$
\sn
\item "{$(f)$}"  $\theta \in \lambda^+ \cap \text{ Reg} \cap N^0_n
\backslash \kappa^+ \Rightarrow \sup(N^0_{n+1} \cap \theta) >
\sup(N^1_n \cap \theta)$.
\endroster
\medskip

\noindent
Let us carry the induction.
\sn
\underbar{For $n=0$}:  No problem.
\bigskip

\noindent
\underbar{For $n+1$}:  Let ${\frak a} =: N^0_n \cap \text{ Reg } \cap 
\lambda^+ \backslash \kappa^+$, so ${\frak a}^n \in {\frak B}_{n+1}$ and 
${\frak a}^n$ is a set of cardinality $\le \kappa$ of regular cardinals 
$\in (\kappa,\lambda^+)$.  

Let $g^n \in \Pi {\frak a}^n$ be defined by $g^n(\theta) =: \sup(N^1_n \cap 
\theta)$.
\mr
\item "{$(*)_1$}"   It is enough to prove that ${\frak a}^n$ is a member of
\endroster

$$
I^n = \{ {\frak b} \subseteq {\frak a}^n:\text{ for some }f \in 
(\Pi {\frak a}^n) \cap {\frak B}_{n+1} \text{ we have } 
g^n \restriction {\frak b} < f\}
$$
\mn
so we need to show ${\frak a}^n \in I^n$. \newline
Clearly
\medskip
\roster
\item "{$(*)_2$}"   $J_{\le \mu}[{\frak a}^n] \subseteq I^n$ 
(in particular all singletons are in $I^n$).
\endroster
\enddemo
\bigskip

\noindent
\underbar{Fact}:  There is $f^* \in {\frak B}_{n+1} \cap \Pi 
{\frak a}^n$ such that:

$$
{\frak b}^n =: \{ \theta \in {\frak a}^n:f^*(\theta) < g^n(\theta)\}
$$
\medskip

\noindent
satisfies

$$
[{\frak b}^n]^{< \kappa} \subseteq J_{\le \lambda}[{\frak a}^n]
$$
\mn
(yes! not $J_{\le \mu}[{\frak a}^n]$).
\bigskip

\demo{Proof}  In ${\frak B}_{n+1}$ there is a list $\{a_{n,\varepsilon}:
\varepsilon < \kappa\}$ of $N^0_n$.  For each $\nu \in T$ let $\nu$ be of
level $\zeta$ and let $N^1_{n,\nu} = \text{ Sk}_{{\frak B}_n}
(\{(a_{n,\varepsilon},\nu \restriction \varepsilon):\varepsilon < \zeta\})$.
So the function $\nu \mapsto N^1_{n,\nu}$ (i.e. the set of pairs
$\langle (\nu,N^1_{n,\nu}):\nu \in T \rangle$ belongs to ${\frak B}_{n+1}$.
Clearly $\langle N^1_{n,\eta \restriction \zeta}:\zeta < \kappa \rangle$ is
increasing continuous with union $N^1_n$.  Let
$g^1_{n,\nu} \in \Pi({\frak a}^n \cap N^1_{n,\nu})$ be defined by
$g^1_{n,\nu}(\theta) = \sup(\theta \cap N^1_{n,\nu})$, so 
$\{({\frak a}^n \cap N^1_{n,\nu},g^1_{n,\nu}):\nu \in T\} \in 
{\frak B}_{n+1}$.  Now $\Pi {\frak a}^n / J_{\le \lambda}[{\frak a}^n]$ is
$\lambda^+$-directed, hence as $|T| \le \lambda$ there is 
$f^* \in \Pi {\frak a}^n$ such that:
\medskip
\roster
\item "{$(*)_3$}"  $\nu \in T \Rightarrow g^1_{n,\nu} <_{J_{\le \lambda}
[{\frak a}^n]} f^*$,
\endroster
\medskip

\noindent
and by the previous sentence without loss of generality $f^* \in 
{\frak B}_{n+1}$.  Note that for $\theta \in {\frak a}^n$ the sequence
$\langle g^1_{n,\eta \restriction \zeta}(\theta):\zeta < \kappa \rangle$
is non-decreasing with limit $g^n(\theta)$.
\newline
Let ${\frak c} = \{\theta \in {\frak a}^n:f^*(\theta) < g^n(\theta)\}$,
now note
\medskip
\roster
\item "{$(*)_4$}"  if $\theta \in {\frak c}$ then for every $\zeta < \kappa$
large enough, $f^*(\theta) < g^1_{n,\eta \restriction \zeta}(\theta)$.
\endroster
\medskip

\noindent
Hence ${\frak c}' \in [{\frak c}]^{< \kappa} \Rightarrow
{\frak c}' \in J_{\le \lambda}[{\frak a}^n]$ as required in the fact
\newline
(why the implication?  because if ${\frak c}' \subseteq {\frak c},
|{\frak c}| < \kappa$ then by $(*)_4$ for some $\zeta < \kappa$ we have
$f^* \restriction {\frak c}' < g'_{n,\eta \restriction \zeta} \restriction
{\frak c}'$ which by $(*)_3$ gives ${\frak c}' \in J_{\le \lambda}
[{\frak a}^n]$).  So let ${\frak b}^n = {\frak c}$. 
\hfill$\square_{\text{fact}}$
\medskip

Now if ${\frak b}^n$ is in $J_{\le \mu}[{\frak a}^n]$ by $(*)_1 + (*)_2$ 
above we can finish the induction step. \newline
If not, some $\tau^* \in \text{ Reg } \backslash \mu^+$ satisfies $\tau^* \in
\text{ pcf}({\frak b}^n)$; let $\langle {\frak c}_\zeta:\zeta < 
\kappa \rangle$ be an increasing continuous sequence of subsets of 
${\frak a}^n$ each of cardinality $< \kappa$ such that ${\frak b}^n = 
\dsize \bigcup_{\zeta < \kappa} {\frak c}_\zeta$ and so (by the face above)
$\zeta < \kappa \Rightarrow \tau^* > \lambda \ge 
\text{ max pcf}({\frak c}_\zeta)$.  
We know that this implies that for some club $E$ of $\kappa$ and
$\theta_\zeta \in \text{ pcf}({\frak c}_\zeta)$, for $\zeta \in E$,
$\tau^* \in \text{ pcf}_{\kappa\text{-complete}}(\{\theta_\zeta:\zeta 
\in E\})$ and $\langle \theta_\zeta:\zeta \in E \rangle$ is strictly
increasing and max pcf$\{\theta_\zeta:\zeta \in E \cap \xi\} \le
\theta_\xi$ for $\xi \in E$, by \cite[Ch.VIII,1.5(2),(3),p.317]{Sh:g}.

Now max pcf$\{\theta_\varepsilon:\varepsilon \in \zeta \cap E\} \le
\text{ max pcf}({\frak c}_\zeta) \le \lambda$ so $\mu < \tau^* \le$ 
the cardinal from
clause (i) of \scite{2.2}, against an assumption.  So we have carried the 
inductive step in defining $N^0_n,N^1_n$.
\medskip

So $N^0_n,N^1_n$ are well defined for every $n$, 
clearly $\dsize \bigcup_{n < \omega} N^0_n \cap
\lambda = \dsize \bigcup_{n < \omega} N^1_n \cap \lambda$ \nl
(see \cite[Ch.IX,3.3A,p.379]{Sh:g}) hence 
$\dsize \bigcup_{n < \omega} N^0_n \cap T = \dsize
\bigcup_{n < \omega} N^1_n \cap T$, hence for some $n,N^0_n \cap
\{\eta \restriction \zeta:\zeta < \kappa\}$ has cardinality $\kappa$.
Now

$$
A = \{ \nu \in T:\text{ for some } \rho \text{ we have } 
\nu \triangleleft \rho \in N^0_n\}
$$
\medskip

\noindent
belongs to ${\frak B}_{n+1} \cap [T]^\kappa$ and 
$\{ \eta \restriction \zeta:\zeta < \kappa\} \subseteq A$. \newline
Contradiction to the choice of $\eta$.
\enddemo
\bigskip

\noindent
\underbar{Second inequality}  Cardinal of $(ii) \le$ cardinal of $(i)$.
\newline
By the proof of \cite[II,3.5]{Sh:g}. \hfill$\square_{\scite{2.2}}$
\bigskip

\definition{\stag{2.3} Definition}   1)  Assume 
$I \subseteq J \subseteq {\Cal P}
(\kappa),I$ an ideal on $\kappa,J$ an ideal or the complement of a filter on
$\kappa$, e.g. $J = {\Cal P}^-(\kappa) = {\Cal P}(\kappa)
\backslash \{\kappa\}$ stipulating \nl
$f \ne_J g \Leftrightarrow \{i < \kappa:f(i) = g(i)\} \in J$.  We let

$$
T^+_{I,J}(f,\lambda) = \sup\{|F|^+:F \in {\Cal F}_{I,J}(f,\lambda)\}
$$
\mn
and

$$
T_{I,J}(f,\lambda) = \sup\{|F|:F \in {\Cal F}_{I,J}(f,\lambda)\}
$$
\mn

$$
\align
{\Cal F}_{I,J}(f,\lambda) = \{ F:&F \subseteq \dsize \prod_{i < \kappa}
f(i) \text{ and } f \ne g \in F \Rightarrow f \ne_J g \\
  &\text{and } A \in I \Rightarrow \lambda \ge |\{ f \restriction A:f \in
F\}| \}.
\endalign
$$
\mn
2) For $J$ an ideal on $\kappa,\theta \ge \kappa$ and $f \in {}^\kappa
(\text{Ord} \backslash \{0\})$, we let

$$
\align
\bold U_J(f,\theta) = \text{ Min} \bigl\{ |{\Cal P}|:&{\Cal P} \subseteq
[\text{sup Rang}(f)]^\theta \text{ and for every } g \in \dsize
\prod_{i < \kappa} f(i) \\
  &\text{for some } a \in {\Cal P} \text{ we have }
\{i < \kappa:g(i) \in a\} \in J^+ \bigr\}.
\endalign
$$
\mn
If $\theta = \kappa$ (= Dom$(J)$), then we may omit $\theta$.  If $f$ is
constantly $\lambda$ we may write $\lambda$ instead of $f$. \nl
3) For $I \subseteq J,I$ ideal on $\kappa,J$ an ideal or compliment of a
filter on $\kappa,\mu \ge \theta \ge \kappa$ and
$f \in {}^\kappa(\text{Ord} \backslash \{0\})$ let

$$
\bold U_{I,J}(f,\theta,\mu) = \sup\{\bold U_J(F,\theta):F \in
{\Cal F}^-_I(f,\mu)\}
$$
\mn
where

$$
\align
{\Cal F}^-_I(f,\mu) = \bigl\{F:&F \subseteq \dsize \prod_{i < \kappa}
f(i) \text{ and} \\
  &A \in I \Rightarrow \mu \ge |\{f \restriction A:f \in F\}| \bigr\}
\endalign
$$

$$
\align
\bold U_J(F,\theta) = \text{ Min} \bigl\{ |{\Cal P}|:&{\Cal P} \subseteq
[\text{sup Rang}(f)]^\theta \text{ and for every } f \in F \\
  &\text{for some } a \in {\Cal P} \text{ we have }
\{i < \kappa:f(i) \in a\} \in J^+ \bigr\}.
\endalign
$$
\enddefinition
\bigskip

\demo{\stag{2.4} Fact}  Let $\lambda \ge \kappa = \text{ cf}(\kappa) >
\aleph_0$. \newline
1)  $\lambda^{\kappa,\text{tr}} = 
T_{J^{bd}_\kappa,{\Cal P}^-(\kappa)}(\lambda,\lambda)$ and
$\lambda^{\langle \kappa,\theta \rangle} \le \bold U_{J^{bd}_\kappa}
(\lambda,\theta)$. \newline
2)  If $\lambda \ge \mu$, \underbar{then} 
$\lambda^{\kappa,\text{tr}} \ge \mu^{\kappa,\text{tr}}$ and
$\lambda^{<\kappa>} \ge \mu^{<\kappa>}$. \newline
3)  $\lambda^{\kappa,\text{tr}} = \lambda^{\langle \kappa \rangle} +
\kappa^{\kappa,\text{tr}}$. \newline
4)  Assume $I \subseteq J$ are ideals on $\kappa$.  Then
$T^+_I(f,\lambda) > \mu$ if:
\medskip
\roster
\widestnumber\item{ $(iii)$ }
\item "{$(i)$}"  each $f(i)$ is a regular cardinal $\lambda_i \in (\kappa,
\lambda)$
\sn
\item "{$(ii)$}"  $\dsize \prod_{i < \kappa}f(i)/J$ is $\mu$-directed
\sn
\item "{$(iii)$}"  for some $A_\zeta \subseteq \kappa$ for $\zeta < \zeta^*
< \underset {j < \kappa}\to {\text{Min}}f(j)$ we have: \newline
max pcf$\{f(i):i \in A_\zeta\} \le \lambda$ (hence cf$\biggl( \dsize \prod
_{i \in A_\zeta} f(i)\biggr) \le \lambda$) and $\{ A_\zeta:\zeta < \zeta^*\}$
generates an ideal on $\kappa$ extending $I$ but included in $J$. 
\ermn
5) $\bold U_J(\lambda) \le \bold U_J(\lambda,\theta) \le 
\bold U_J(\lambda) + \text{ cf}([\theta]^\kappa,\subseteq) \le 
\bold U_J(\lambda) + \theta^\kappa$ and $T_I(f) \le \bold U_I(f) + 2^\kappa$ 
and $\bold U_{I,J}(f,\lambda) \le T_{I,J}(f,\lambda) \le \bold U_{I,J}
(f,\lambda) + 2^\kappa$ where $I \subseteq J$ are ideals on $\kappa$. \nl
Also obvious monotonicity properties (in $I,J,\lambda,\theta,f$) holds.
\enddemo
\bigskip

\demo{Proof}  1) Easy.  Let us prove the first equation.
First assume $F \in {\Cal F}_{J^{bd}_\kappa,{\Cal P}^-(\kappa)}
(\lambda,\lambda)$, and we define a tree as follows: for $i < \kappa$ the
$i$th level is

$$
T_i = \{f \restriction i:f \in F\}
$$
\mn
and

$$
T = \dsize \bigcup_{i < \kappa} T_i, \text{ with the natural order }
\subseteq.
$$

\noindent 
Clearly $T$ is a tree with $\kappa$ levels, the $i$-th level being $T_i$.
\newline
By the definition of ${\Cal F}_{J^{bd}_\kappa,{\Cal P}^-(\kappa)}(\lambda,
\lambda)$ as $i < \kappa \Rightarrow \{j:j < i\} \in J^{bd}_\kappa$,
clearly $|T_i| \le \lambda$.  Now for each $f \in F$, clearly
$t_f =: \langle (f \restriction i):i < \kappa \rangle$ is a $\kappa$-branch
of $T$, and $f_1 \ne f_2 \in F \Rightarrow t_{f_1} \ne t_{f_2}$ so $T$ has
at least $|F|\,\kappa$-branches.

The other direction is easy, too.  Note that the proof gives $=^+$; i.e.,
the supremum is obtained in one side iff it is obtained in the other side.
\newline
2)  If $T$ is a tree with $\mu$ nodes and $\kappa$ levels then we can add
$\lambda$ nodes adding $\lambda$ branches. Also the other inequality is
trivial.  \newline
3)  First $\lambda^{\kappa,\text{tr}} \ge \lambda^{\langle \kappa \rangle}$ 
because if $T$ is a tree with $\lambda$ nodes and $\kappa$ levels, then we 
know $|\text{lim}_\kappa(T)| \le \lambda^{\kappa,\text{tr}}$, hence 
${\Cal P} = \{t:t \text{ is a } \kappa \text{-branch of }T\}$ has 
cardinality $\le \lambda^{\kappa,\text{tr}}$ and satisfies the 
requirement in the definition of $\lambda^{<\kappa>}$.

Second $\lambda^{\kappa,\text{tr}} \ge \kappa^{\kappa,\text{tr}}$ 
by part (2) of \scite{2.4}.

Lastly, $\lambda^{\kappa,\text{tr}} \le \lambda^{<\kappa>} + 
\kappa^{\kappa,\text{tr}}$ 
because if $T$ is a tree with $\lambda$ nodes and $\kappa$ levels, we know
by Definition \scite{2.1} that there is ${\Cal P} \subseteq [T]^\kappa$ of 
cardinality
$\le \lambda^{<\kappa>}$ such that every $\kappa$-branch of $T$ is included
in some $A \in {\Cal P}$, without loss of generality $x <_T y \in A \in
{\Cal P} \Rightarrow x \in A$; so

$$
\align
|\text{lim}_\kappa(T)| &= |\{t:t \text{ a } \kappa \text{-branch of } T\}| \\
  &= | \dsize \bigcup_{A \in {\Cal P}} \{t \subseteq A:t \text{ a }
\kappa \text{-branch of } T\}| \\
  &\le \dsize \sum_{A \in {\Cal P}} |\text{lim}_\kappa(T \restriction A)| \\
  &\le |{\Cal P}| + \kappa^{\kappa,\text{tr}} \le \lambda^{<\kappa>} +
\kappa^{\kappa,tr}.
\endalign
$$

\noindent
4)  Like the proof of \cite[Ch.II,3.5]{Sh:g}. \nl
5)  Left to the reader.   \hfill$\square_{\scite{2.4}}$
\enddemo
\bigskip

\proclaim{\stag{2.5} Lemma}  Assume
\medskip
\roster
\item "{$(a)$}"  $I \subseteq J$ are ideals on $\kappa$
\sn
\item "{$(b)$}"  $I$ is generated by $\le \mu^*$ sets, 
$\mu^* \ge \kappa$
\sn
\item "{$(c)$}"  $T^+_{I,J}(f,\lambda) > \mu = \text{ cf}(\mu) > 
\mu^* \ge T_{I,J}(\mu^*,\kappa)$
\sn
\item "{$(d)$}"  $\kappa$ is not the union of countably many members of $I$.
\endroster
\medskip

\noindent
\underbar{Then} We can find $A_0 \subseteq A_1 \subseteq \cdots \subseteq
A_n \subseteq \ldots$ from $I^+$ with union $\kappa$, such that for each
$n$ there is $\langle \lambda^n_i:i \in A_n \rangle,\mu^* < \lambda^n_i
= \text{ cf}(\lambda^n_i) \le f(i)$ such that:

$$
\dsize \prod_{i \in A_n} \lambda^n_i/J \text{ is } \mu \text{-directed}
$$

$$
A \subseteq A_n,A \in I \Rightarrow \text{ cf}(\dsize \prod_{i \in A}
\lambda^n_i) \le \lambda.
$$
\endproclaim
\bigskip

\remark{\stag{2.5A} Remark}  The point in the proof is that if $I$ is generated
by $\{B_\gamma:\gamma < \gamma^* \le \mu^*\}$, and 
$\{ \eta_\alpha:\alpha < \mu^+\}$ are distinct branches and 
$f \in {}^A(\lambda + 1 \backslash \{ 0 \}),A \subseteq \kappa$ and 
$i \in A \Rightarrow \text{ cf}(f(i)) > \mu^*$, then for some $g < f$ for 
every $\gamma < \gamma^*$ and $\alpha < \mu^+$, \nl
$\{i < \gamma:\text{if }\eta_\alpha(i) < f(i) \text{ then } 
\eta_\alpha(i) < g(i)\} = \gamma \text{ mod }
J_{< \lambda^+}(f \restriction \gamma)$.
\endremark
\bigskip

\demo{Proof}  Similar to the proof of \scite{1.1} adding the main point of
the proof of \scite{2.2}, the ``fact" there. 
\enddemo
\bigskip

\noindent
\ub{We can further generalize}
\definition{\stag{2.6} Definition}  
For $I \subseteq J \subseteq {\Cal P}(\kappa)$, function $f^* \in {}^\kappa
\text{Reg}$ and $\lambda$, we let 

$$
\align
{\Cal F}^1_{(I,J,\lambda)}(f^*) = \biggl\{
F \subseteq \dsize \prod_{i < \kappa} f^*(i):&\text{if }A \in J 
\text{ then} \\
  &\lambda \ge | \{(f \restriction A)/I:f \in {\Cal F}\}| \biggr\}
\endalign
$$
\mn
(so $I$ is \wilog \, an ideal on $\kappa$ and this is just
${\Cal F}^-_I(f^*,\lambda))$

$$
\align
{\Cal F}^2_{(I,J,\lambda)}(f^*) = \biggl\{
F \subseteq \dsize \prod_{i < \kappa}f^*(i):&\text{if }A \in J, \text{ and }
f,g \in F \text{ are distinct} \\
  &\text{then } \{i \in A:f(i) = g(i)\} \in I \biggr\}
\endalign
$$

$$
\align
{\Cal F}^3_{(I,J,\lambda,\bar \theta)}(f^*) = \biggl\{
F \subseteq \dsize \prod_{i < \kappa}f^*(i):&\text{if }A \in J, \text{ then
for some} \\
  &G \subseteq \dsize \prod_{i \in A} [f^*(i)]^{\theta_i} 
\text{ of cardinality } \le
\lambda \text{ we have} \\
  &(\forall f \in F)(\exists g \in G)\{i \in A:f(i) \notin g(i)\} \in I
\biggr\}.
\endalign
$$
\medskip

\noindent
If $\Xi$ is a set of such tuples ${\Cal F}^\ell_\Xi(f^*) = 
\dsize \bigcap_{\Upsilon \in \Xi} {\Cal F}^\ell_\Upsilon(f^*)$  \newline
If in all the tuples $\lambda$ is the third element, we write triples and
$f,\lambda$ instead of $f$.

For any ${\Cal F}^\ell_\Upsilon$ we let $T^\ell_\Upsilon(f^*) = \sup
\{|F|:F \in {\Cal F}^\ell_x(f^*)\}$ but: instead of $T$ we have $F \in
{\Cal F}_I(f)$ exemplifying $\bold U_{I,J}(f,\lambda) > \mu$; i.e.
$\bold U_{I,J}(F,\lambda) > \mu$.  Then $\eta \in F$ satisfies $(\forall A
\in {\Cal P})[\{i:\eta(i) \in A\} \in J]$.  We choose $N^0_n,N^1_n$
satisfying (a)-(f) with $\gamma_n =1$.
\enddefinition
\newpage

\head {\S3 On the depth behaviour for ultraproducts} \endhead  \resetall
\bigskip

\noindent
The problem \footnote""{Done 24/Feb/95-Proof read 4/4/95} 
originates from Monk \cite{M} and see on it Roslanowski
Shelah \nl
\cite{RoSh:534} and then \cite[\S3]{Sh:506} but the presentation is
self-contained. \newline
We would like to have (letting $B_i$ denote Boolean algebra), for $D$ an
ultrafilter \nl
on $\kappa$:

$$
\text{Depth}(\dsize \prod_{i < \kappa} B_i/D) \ge \biggl| \dsize \prod_{i <
\kappa} \text{ Depth}(B_i)/D \biggr|.
$$
\medskip

\noindent
(If $D$ is just a filter, we should use $T_D$ instead of product in the right
side).  Because of the problem of attainment (serious by 
Magidor Shelah \cite{MgSh:433}), we rephrase the
question:
\medskip
\roster
\item "{$\bigotimes$}"  for $D$ an ultrafilter on $\kappa$, does $\lambda_i 
< \text{ Depth}^+(B_i)$ for $i < \kappa$ imply

$$
\biggl| \dsize \prod_{i < \kappa} \lambda_i/D \biggr| < 
\text{ Depth}^+(\dsize \prod_{i < \kappa} B_i/D)
$$

\noindent
at least when $\lambda_i > 2^\kappa$;
\medskip
\noindent
\item "{$\bigotimes'$}"  for $D$ a filter on $\kappa$ does $\lambda_i <
\text{ Depth}^+(B_i)$ for $i < \kappa$ imply (assuming \newline
$\lambda_i > 2^\kappa$ for simplicity):

$$
\align
\mu = \text{ cf}(\mu) &< T^+_{D+A}(\langle \lambda_i:i < \kappa \rangle)
\text{ for some } A \in D^+ \Rightarrow \\
  &\mu < \text{ Depth}^+
(\dsize \prod_{i < \kappa} B_i/(D+A)) \text{ for some } A \in D^+.
\endalign
$$
\endroster
\medskip

\noindent
As found in \cite{Sh:506}, this actually is connected to a pcf problem,
whose answer under reasonable restrictions is \scite{1.4}.  So now we can 
clarify the connections.

Also, by changing the invariant (closing under homomorphisms, see \cite{M})
we get a nicer result; this shall be dealt with here.  

The results here (mainly \scite{3.3}) supercede \cite[3.26]{Sh:506}.
\bigskip

\definition{\stag{3.1} Definition}  1) For a partial order $P$ (e.g. a Boolean
algebra) let \newline
$\text{Depth}^+(P) = \text{ min}\{ \lambda:\text{we cannot find } a_\alpha
\in P$ for $\alpha < \lambda$ such that  \newline

$\qquad \qquad \qquad \qquad \quad \alpha < \beta \Rightarrow a_\alpha
<_P a_\beta\}$. \newline
2) For a Boolean algebra $B$ let \newline
$D^+_h(B) = \text{ Depth}^+_h(B) = \sup
\{\text{Depth}^+(B'):B'$ is a homomorphic image of $B\}$. \newline
3)  $\text{Depth}(P) = \sup\{\mu$: there are $a_\alpha \in P$ for $\alpha <
\mu$ such that \newline

$\qquad \qquad \qquad \qquad \quad \alpha < \beta < \mu 
\Rightarrow a_\alpha <_P a_\beta\}$.
\newline
4)  $\text{Depth}_h(P) = D_h(P) = \sup\{\text{Depth}(B'):B'$ is a 
homomorphic image of $B\}$. \newline
5) We write $D_r$ or $Dh_r$ or $\text{Depth}_r$ if we restrict ourselves 
to regular cardinals.  Of course we could have looked at the ordinals.
\enddefinition
\bigskip

\definition{\stag{3.1A} Definition}  1) For a linear order ${\Cal I}$, let the
interval Boolean algebra, $BA[{\Cal I}]$ be the Boolean algebra of subsets
of ${\Cal I}$ generated by $\{[s,t)_{\Cal I}:s < t \text{ are from}$
\newline
$\{- \infty\} \cup {\Cal I} \cup \{+ \infty\}\}$. \newline
2) For a Boolean algebra $B$ and regular $\theta$, 
let com$_{< \theta}(B)$ be the $(< \theta)$-completion of $B$, that is the
closure of $B$ under the operations $-x$ and 
$\dsize \bigvee_{i < \alpha} x_i$ for $\alpha < \theta$ inside the 
completion of $B$.
\enddefinition
\bigskip

\demo{\stag{3.2} Fact}  1) If $B$ is the interval 
Boolean algebra of the ordinal $\gamma \ge \omega$ \underbar{then}
\medskip
\roster
\item "{$(a)$}"  $D^+_h(B) = |\gamma|^+$ 
\sn
\item "{$(b)$}"  Depth$^+(B) = |\gamma|^+$.
\endroster
\medskip

\noindent
2) If $B'$ is a subalgebra of a homomorphic image of $B$, then
$D_h^+(B) \ge D_h^+(B')$. \newline
3) If $D' \supseteq D$ are filters on $\kappa$ and for $i < \kappa,B'_i$ is a subalgebra of
a homomorphic image of $B_i$ then:
\medskip
\roster
\item "{$(\alpha)$}"  $\dsize \prod_{i < \kappa} B'_i/D'$ is a subalgebra of
a homomorphic image of $\dsize \prod_{i < \kappa} B_i/D$, hence
\sn
\item "{$(\beta)$}"  $D_h^+(\dsize \prod_{i < \kappa} B_i/D) \ge D_h^+
(\dsize \prod_{i < \kappa} B'_i/D')$.
\ermn
4) In parts (2), (3) we can replace $D_h$ by $D$ if we omit ``homomorphic
image".
\enddemo
\bigskip

\demo{Proof}  Straightforward.
\enddemo
\bigskip

\proclaim{\stag{3.2A} Claim}  1) If $D$ is a filter on $\kappa$ and for
$i < \kappa$, $B_i$ a Boolean algebra, \newline
$\lambda_i < \text{ Depth}^+_h(B_i)$ \underbar{then}
\medskip
\roster
\item "{$(a)$}"  $\text{Depth}^+_h(\dsize \prod_{i < \kappa} B_i/D) \ge \sup
_{D_1 \supseteq D} \left(\text{tcf }(\dsize 
\prod_{i < \kappa} \lambda_i/D_1) \right)^+$ \newline
(i.e. $\sup$ on the cases tcf is well defined)
\smallskip
\noindent
\item "{$(b)$}"  $\text{Depth}^+_h(\dsize \prod_{i < \kappa} B_i/D) 
\text{ is }
\ge \text{ Depth}^+_h({\Cal P}(\kappa)/D)$ and is at least

$$
\sup\{[\text{tcf }(\dsize \prod_{i < \kappa} \lambda'_i/D_1)]^+:\lambda'_i
< \text{Depth}^+(B_i),D_1 \supseteq D \}.
$$
\endroster
\medskip

\noindent
2) $\mu < \text{ Depth}^+_h(B)$ \underbar{iff} for some $a_i \in B$ for
$i < \mu$ we have that: $\alpha < \beta < \mu,n < \omega$, and 
$\alpha_\ell < \beta_\ell < \mu$ for $\ell < n$ together imply that 
$B \models ``(a_\beta - a_\alpha) - 
\dsize \bigcup_{\ell < n} (a_{\alpha_\ell} - a_{\beta_\ell}) > 0"$. \nl
3) Let $A \in D^+$ ($D$ a filter on $\kappa$).  In 
$\dsize \prod_{i < \kappa} B_i/D$ there is a chain of order type
$\Upsilon$ \ub{if} in $\dsize \prod_{i < \kappa} B_i/(D+A)$ 
there is such a chain.  If $\Upsilon = \lambda$; cf$(\lambda) > 2^kappa$
also the inverse is true.  \newline
4) If $\mu < \text{ Depth}^+(\dsize \prod_{i < \kappa} B_i/D)$ and
cf$(\mu) > 2^\kappa$, then we can find $A \in D^+$ and \nl
$f_\alpha \in \dsize \prod_{i < \kappa} B_i$ for $\alpha < \mu$ such 
that letting $D^* = D+A$: \newline
$\alpha < \beta < \mu \Rightarrow (\dsize \prod_{i < \kappa} B_i/D^*)
\models f_\alpha/D^* < f_\beta/D^*$ moreover $f_\alpha <_{D^*} f_\beta$. \nl
5) Like (1) replacing Depth$^+_h$ by Depth$^+,D_1 \supseteq D$ by
$\{D+A:A \in D^+\}$.
\endproclaim
\bigskip

\demo{Proof}  Check, e.g.: \newline
2) \ub{The ``if" direction}:
\sn
Let $I$ be the ideal of $B$ generated by $\{a_\alpha - a_\beta:\alpha <
\beta < \mu\},h:B \rightarrow B/I$ the canonical homomorphism, so $\langle
a_\alpha/I:\alpha < \mu \rangle$ is strictly increasing in $B/I$.
\mn
\ub{The ``only if" direction}:
\sn
Let $h$ be a homomorphism from $B$ onto $B_1$ and $\langle b_\alpha:\alpha <
\mu \rangle$ be a (strictly) increasing sequence of elements of $B_1$.  Choose
$a_\alpha \in B$ such that $h(a_\alpha) = b_\alpha$, so $\alpha < \beta
\Rightarrow a_\alpha \backslash a_\beta \in \text{ Ker}(h)$ but $a_\alpha
\notin \text{ Ker}$. \nl
3) The first implication is trivial, the second follows from part (4). \nl
4) First, assume $\mu$ is regular.  
Let $\langle f_\alpha/D:\alpha < \mu \rangle$ exemplify 
$\mu < \text{ Depth}^+(\dsize \prod_{i < \kappa} B_i/D)$. \ub{Then}
$\alpha < \beta < \mu \Rightarrow f_\alpha \le_D f_\beta \and \neg(f_\alpha
=_D f_\beta)$, so for each $\alpha \langle \{i < \kappa:f_\alpha(i) =
f_\beta(i)\}/D:\beta < \mu \rangle$ is decreasing and $|2^\kappa/D| <
\mu = \text{ cf}(\mu)$ hence for some $\beta_\alpha \in (\alpha,\mu)$
we have $(\forall \beta)(\beta_\alpha \le \beta < \mu \Rightarrow \{i <
\kappa:f_\alpha(i) \ne f_{\beta_\alpha}(i)\} = \{i < \kappa:f_\alpha(i) \ne
f_\beta(i)\} \text{ mod } D$ (as $f_\gamma/D$ is increasing).  So
$\langle\{i:f_\alpha(i) = f_{\beta_\alpha}(i)\}/D:\alpha < \mu \rangle$ is
decreasing and $|2^\kappa/D| \le 2^\kappa < \mu$, hence for some $A^* 
\subseteq \kappa$ the set $E = \{ \alpha < \mu:\{i < \kappa:f_\alpha(i) <
f_{\beta_\alpha}(i)\} = A^* \text{ mod } D\}$ is unbounded and even
stationary in $\mu$.  Let $D^* = 
D+A^*$, so for $\alpha < \beta < \mu$ we have $f_\alpha \le_D f_\beta$ hence
$f_\alpha \le_{D^*} f_\beta$, but $\alpha \in E \and \beta \ge \beta_\alpha
\Rightarrow f_\alpha \ne_{D^*} f_\beta$.  
Hence there is $E' \subseteq \{ \delta \in E:(\forall \alpha < \delta \cap E)
(\beta_\alpha < \delta)\}$ is unbounded in $\mu$ and clearly 
$(\forall \alpha,\beta)(\alpha < \beta \and \alpha \in E' \and 
\beta \in E' \Rightarrow f_\alpha <_{D^*} f_\beta)$.

So $\{f_\alpha:\alpha \in E'\}$ exemplifies the conclusion. \newline
Second, if $\mu$ is singular, let $\mu = \dsize \sum_{\zeta < \text{ cf}
(\mu)} \mu_\zeta,\mu_\zeta > 2^\kappa;\mu_\zeta$ strictly increasing and
each $\mu_\zeta$ is regular.  So
given $\langle f_\alpha:\alpha < \mu \rangle$, for each $\zeta < 
\text{ cf}(\mu)$ we can find $E_\zeta \subseteq \mu^+_\zeta$ of cardinality
$\mu^+_\zeta$ and $A_\zeta \in D^+$ such that $\alpha \in E_\zeta \and
\beta \in E_\zeta \and \alpha < \beta \Rightarrow f_\alpha <_{D+A_\zeta}
f_\beta$.  For some $A$, cf$(\mu) = \sup\{ \zeta:A_\zeta = A\}$; so
$A$ and the $f_\alpha$'s for $\alpha \cup \{E_\zeta \backslash
\{\text{Min}(E_\zeta)\}:\zeta < \text{ cf}(\mu)$ is such that $A_\zeta = A\}$ 
are as required.  \hfill$\square_{\scite{3.2A}}$
\enddemo
\bigskip

\noindent
We now give lower bound of depth of reduced products of Boolean algebras 
$B_i$ from the depths of the $B_i$'s.
\proclaim{\stag{3.3} First Main Lemma}  Let $D$ be a filter on $\kappa$ and
$\langle \lambda_i:i < \kappa
\rangle$ a sequence of cardinals $(> 2^\kappa)$ and $2^\kappa < \mu 
= \text{ cf}(\mu)$.
\underbar{Then}: \newline
1)  $(\alpha) \Leftrightarrow (\alpha)^+ \Leftrightarrow (\beta) 
\Leftrightarrow (\beta)^- \Leftrightarrow (\gamma)^+ \Rightarrow (\gamma)
\Rightarrow (\delta)$. \newline
2) If in addition $(\forall \sigma < \mu)(\sigma^{\aleph_0} < \mu) \vee (D$ 
is $\aleph_1$-complete) we also 
have $(\gamma) \Leftrightarrow (\gamma)^+ \Leftrightarrow (\delta)$ so all
clauses are equivalent \nl
where:
\medskip
\roster
\item "{$(\alpha)$}"  if $B_i$ is a Boolean algebra, $\lambda_i \le
\text{ Depth}^+(B_i)$ \underbar{then} $\mu < \text{ Depth}^+
(\dsize \prod_{i < \kappa} B_i/D)$
\sn
\item "{$(\beta)$}"  there are cardinals $\gamma_i < \lambda_i$ for 
$i < \kappa$ such that, letting $B_i$ be \newline
$BA[\gamma_i] =$ the interval
Boolean algebra of (the linear order) $\gamma_i$, we have \newline
$\mu < \text{ Depth}^+(\dsize \prod_{i < \kappa} B_i/D)$
\sn
\item "{$(\gamma)$}"  there are $\langle \langle \lambda_{i,n}:n < n_i
\rangle:i < \kappa \rangle$ where $\lambda_{i,n} = \text{ cf}(\lambda_{i,n}) <
\lambda_i$ and a non-trivial filter $D^*$ on 
$\dsize \bigcup_{i < \kappa} (\{ i \} \times n_i)$ such that:
{\roster
\itemitem{ $(i)$ }  $\mu = 
\text{ tcf}(\dsize \prod_{(i,n)} \lambda_{i,n}/D^*)$ 
\sn
\itemitem{ $(ii)$ }  for some $A^* \in D^+ \text{ we have}$ \nl
$D + A^* =  \{ A \subseteq \kappa:\text{the set } \dsize \bigcup_{i \in A}
(\{ i\} \times n_i) \text{ belongs to } D^*\}$
\endroster}
\item "{$(\delta)$}"  for some filter $D' = D+A,A \in D^+$ and cardinals
$\lambda'_i < \lambda_i$ \newline
we have $\mu \le T_{D'}(\langle \lambda_i:i < \kappa \rangle)$
\smallskip
\noindent
\item "{$(\beta)'$}"   like $(\beta)$ we allow $\gamma_i$ to be an ordinal
\sn
\item "{$(\beta)^-$}"  letting $B_i$ be the disjoint sum of $\{BA[\gamma]:
\gamma < \lambda_i\}$ we have: \newline
$\mu < \text{ Depth}^+(\dsize \prod_{i < \kappa} B_i/D)$. 
\item "{$(\gamma)^+$}"  for some filter $D^*$ of the form $D+A$ and
$\lambda'_i = \text{ cf}(\lambda'_i) < \lambda_i$ we have
$\mu = \text{ tcf}\left( \dsize \prod_{i < \kappa} \lambda'_i/D^* \right)$
\item "{$(\alpha)^+$}"  if $B_i$ is a Boolean algebra, $\lambda_i \le
\text{ Depth}^+(B_i)$ then for some $A \in D^+$ we have, setting
$D^* = D+A$, that $\mu < \text{ Depth}^+
\left( \dsize \prod_{i < \kappa} B_i,<_{D^*} \right)$; moreover for some
$f_\alpha \in \dsize \prod_{i < \kappa} B_i$ for $\alpha < \mu$ we have 
$\alpha < \beta \Rightarrow \{i:B_i \models f_\alpha(i) < f_\beta(i)\} 
= \kappa \text{ mod } D^*$.
\endroster
\endproclaim
\bigskip

\demo{Proof}  1) We shall prove $(\alpha) \Leftrightarrow (\beta)(\beta)'
\Rightarrow (\beta)^- \Rightarrow (\gamma)^+ \Rightarrow (\beta)$ and
$(\alpha)^+ \Leftrightarrow (\alpha)$ and 
$(\gamma)^+ \Rightarrow (\gamma) \Rightarrow (\delta)$.

This suffices. \newline
Now for $(\alpha)^+ \Rightarrow (\alpha)$ note that if ($\lambda_i,B_i$ for
$i < \kappa$ are given and) $A \in D^+$, \nl
$\langle f_\alpha:\alpha < \lambda
\rangle$ exemplify $(\alpha)^+$ then letting $f'_\alpha = (f_\alpha
\restriction A) \cup 0_{(\kappa \backslash A)}$; i.e. $f'_\alpha(i)$ is
$f_\alpha(i)$ when $i \in A$ and $0_{B_i}$ if $i \in \kappa \backslash A$,
easily $\langle f'_\alpha:\alpha < \lambda \rangle$ exemplify $(\alpha)$.
Next $(\alpha) \Rightarrow (\alpha)^+$ by \scite{3.2A}(4).

Now $(\beta) \Rightarrow (\beta)' \Rightarrow (\beta)^-$ holds trivially 
and for $(\beta)' \Rightarrow (\gamma)^+$ repeat the proof of 
\cite[3.24,p.35]{Sh:506} or the relevant part of the proof of \scite{3.4} 
below (with appropriate changes).  Also $(\beta)^- \Rightarrow (\beta)'$ as
in the proof of \scite{3.4} below.  Easily $(\gamma)^+ \Rightarrow (\beta)$;
also $(\beta) \Rightarrow (\alpha)$ because
\mr
\widestnumber\item{$(iii)$}
\item "{$(i)$}"  if $\gamma_i$ a cardinal $<$ Depth$^+(B_i)$, the Boolean
Algebra $BA[\gamma_i]$ can be embedded into $B_i$, and
\sn
\item "{$(ii)$}"  if $B'_i$ is embeddable into $B_i$ for $i < \kappa$ then
$B' = \dsize \prod_{i < \kappa} B'_i/D$ can be embedded into
$\dsize \prod_{i < \kappa} B_i/D$
\sn
\item "{$(iii)$}"  if $B'$ is embeddable into $B$ then Depth$^+(B') \le
\text{ Depth}^+(B)$.
\ermn
Now $(\alpha) \Rightarrow (\beta)$ trivially.  Also $(\gamma)^+ \Rightarrow 
(\gamma) \Rightarrow (\delta)$
trivially.  Next we note $(\beta) \Rightarrow (\delta)$, as if
$B_i = BA[\gamma_i]$ and $\gamma_i < \lambda_i$ and $\mu < \text{ Depth}^+
(\Pi B_i/D)$, then by \scite{3.2A}(4) there is a sequence $\langle f_\alpha:
\alpha < \mu \rangle$ satisfying $f_\alpha \in \dsize \prod_{i < \kappa}
B_i$ and $A^* \in D^+$ such that $\alpha < \beta < \mu \Rightarrow f_\alpha
<_{D+A} f_\beta$.  So $\{f_\alpha:\alpha < \mu\}$ exemplifies that
$T_{D+A}(\langle |B_i|:i < \kappa \rangle) \ge \mu$, as required in clause
$(\delta)$. \newline
2) Assume $(\forall \sigma < \mu)(\sigma^{\aleph_0} < \mu) \vee (D$
is $\aleph_1$-complete). \newline
Now \scite{1.4} gives $(\gamma) \Leftrightarrow (\gamma)^+ \Leftrightarrow
(\delta)$. \hfill$\square_{\scite{3.3}}$
\enddemo
\bigskip

\noindent
Now we turn to the other variant, $D^+_h$.
\proclaim{\stag{3.4} Second Main Lemma}  Let $D$ be a filter on $\kappa$ and
$\langle \lambda_i:i < \kappa
\rangle$ be a sequence of cardinals $(> 2^\kappa)$ and $2^\kappa < \mu 
= \text{ cf}(\mu)$. \underbar{Then} (see below on $(\alpha),\ldots$): \newline
1)  $(\alpha) \Leftrightarrow (\alpha)^+ \Leftrightarrow (\beta) 
\Leftrightarrow (\beta)^- \Leftrightarrow (\gamma)$ and $(\gamma)^+
\Rightarrow (\gamma) \Leftrightarrow (\beta) \Rightarrow (\delta)$.
\newline
2) If $(\forall \sigma < \mu)(\sigma^{\aleph_0} < \mu) \vee (D$ is 
$\aleph_1$-complete) we also have $(\beta) \Leftrightarrow (\gamma) 
\Leftrightarrow (\gamma)^+
\Leftrightarrow (\delta)$ (so all clauses are equivalent) \nl
where:
\medskip
\roster
\item "{$(\alpha)$}"  if $B_i$ is a Boolean algebra, $\lambda_i \le
\text{ Depth}^+_h(B_i)$ \underbar{then} $\mu < \text{ Depth}^+_h
(\dsize \prod_{i < \kappa} B_i/D)$
\sn
\item "{$(\beta)$}"  there are cardinals $\gamma_i < \lambda_i$ for 
$i < \kappa$ such that, letting $B_i$ be \newline
$BA[\gamma_i] =$ the interval
Boolean algebra of (the linear order) $\gamma_i$, we have \newline
$\mu < \text{ Depth}^+_h(\dsize \prod_{i < \kappa} B_i/D)$
\sn
\item "{$(\gamma)$}"  there are $\langle \langle \lambda_{i,n}:n < n_i
\rangle:i < \kappa \rangle$ where $\lambda_{i,n} = \text{ cf}(\lambda_{i,n}) <
\lambda_i$ and a non-trivial filter $D^*$ on 
$\dsize \bigcup_{i < \kappa} \{ i \} \times n_i$ such that:

$$
\mu = \text{ tcf}(\dsize \prod_{(i,n)} \lambda_{i,n}/D^*) \text{ and }
D \subseteq \{ A \subseteq \kappa:\text{the set } \dsize \bigcup_{i \in A}
\{ i\} \times n_i \text{ belongs to } D^*\}
$$

\noindent
\item "{$(\delta)$}"  for some filter $D^* \supseteq D$ and cardinals
$\lambda'_i < \lambda_i$ we have \newline
$\mu \le T_{D^*}(\langle \lambda_i:i < \kappa \rangle)$
\sn
\item "{$(\beta)'$}"  like $(\beta)$ but allowing $\gamma_i$ to be any
ordinal $< \lambda_i$
\sn
\item "{$(\beta)^-$}"  letting $B_i$ be the disjoint sum of $\{BA[\gamma]:
\gamma < \lambda_i\}$ (so Depth$^+(B_i) = \lambda_i$) we have: \newline
$\mu < \text{ Depth}^+_h(\dsize \prod_{i < \kappa} B_i/D)$ 
\sn
\item "{$(\gamma)^+$}"  there are $\lambda'_i = \text{ cf}(\lambda'_i) \in
(2^\kappa,\lambda_i)$ for $i < \kappa$ and filter $D^*_1 \supseteq D$
such that $\dsize \prod_{i \in A} \lambda'_i/D^*$ has true cofinality $\mu$
\sn
\item "{$(\alpha)^+$}"  if $B_i$ is a Boolean algebra, $\lambda_i \le
\text{ Depth}^+_h(B_i)$ \underbar{then} for some filter $D^* \supseteq D$ 
we have $\mu < \text{ Depth}^+_h 
\left( \dsize \prod_{i < \kappa} B_i/D^* \right)$.
\endroster
\endproclaim
\bigskip

\demo{Proof}  Now $(\beta) \Rightarrow (\beta)'$ trivially and $(\beta)' 
\Rightarrow (\beta)^-$ by \scite{3.2}(3) as 
$BA[\gamma_i]$ can be embedded into $B_i$, and similarly 
$(\beta) \Rightarrow (\alpha)$ by \scite{3.2}(3), and
$(\alpha) \Rightarrow (\beta)$ trivially.  Also $(\alpha) \Rightarrow
(\alpha)^+$ trivially and $(\alpha)^+ \Rightarrow (\alpha)$ easily (e.g. by
\scite{3.2}(3)). \newline
Also $(\gamma)^+ \Rightarrow (\beta)$ trivially and $(\gamma)^+ \Rightarrow
(\gamma)$ trivially.

We shall prove below 
$(\gamma) \Rightarrow (\beta)',(\beta) \Rightarrow (\gamma)$ and
$(\beta)^- \Rightarrow (\beta)'$. \newline
Together we have $(\alpha) \Rightarrow (\alpha)^+ \Rightarrow (\alpha)
\Rightarrow (\beta) \Rightarrow (\beta)' \Rightarrow (\beta)^- \Rightarrow 
(\beta)' \Rightarrow (\gamma) \Rightarrow (\beta) \Rightarrow (\alpha)$ 
and $(\gamma)^+ \Rightarrow (\gamma) \Rightarrow
(\delta)$; this is enough for part (1).

Lastly, to prove part (2) of \scite{3.4}, by part (1) it is enough to prove
$(\delta) \Rightarrow (\gamma)^+$ as in the proof of \scite{3.3}.
\enddemo
\bigskip

\noindent
\underbar{$(\gamma) \Rightarrow (\beta)$} \newline
So we have $\lambda_{i,n}$ (for $n < n_i,i < \kappa),D^*$ 
as in clause $(\gamma)$ 
and let $\langle g_\varepsilon:\varepsilon < \mu \rangle$ be \newline 
$<_{D^*}$-increasing cofinal in $\dsize \prod_{(i,n)} \lambda_{i,n}$ but
abusing notation we may write $g_\varepsilon(i,n)$ for $g_\varepsilon((i,n))$.
Let $\gamma_i =: \text{ max}\{ \lambda_{i,n}:
n < n_i\}$ and $B_i =: BA[\gamma_i]$, clearly $\gamma_i < \lambda_i$, a
(regular) cardinal as by
assumption $\lambda_{i,n} < \lambda_i \le \text{ Depth}^+(B_i)$.  
In $B_i$ we have a strictly increasing sequence of length $\gamma_i$.  
Without loss of generality $\{\lambda_{i,n}:n < n_i\}$ is with no repetition
(see \cite[I,1.3(8)]{Sh:g}) and $\lambda_{i,0} > \lambda_{i,1} > \cdots 
\lambda_{i,n_i-1}$.

So for each $i$ we can find 
$a_{i,n} \in B_i$ (for $n < n_i$) pairwise disjoint and \newline
$\langle a_{i,n,\zeta}:\zeta < \lambda_{i,n} \rangle$ (again in $B_i$)
strictly increasing and $< a_{i,n}$.  
\medskip

\noindent
Let $b_{i,\varepsilon} \in B_i$ be $\dsize \bigcup_{n < n_i}
a_{i,n,g_\varepsilon(i,n)}$ (it is a finite union of members of $B_i$ hence
a member of $B_i$).  Let $b_\varepsilon \in \dsize 
\prod_{i < \kappa} B_i/D$ be $b_\varepsilon = \langle 
b_{i,\varepsilon}:i < \kappa \rangle/D$.  Let $J$
be the ideal of \newline
$B =: \dsize \prod_{i < \kappa} B_i/D$ generated by
$\{ b_\varepsilon - b_\zeta:\varepsilon < \zeta < \mu \}$.  Clearly
$\varepsilon < \zeta < \mu \Rightarrow b_\varepsilon \le b_\zeta
\text{ mod } J$, so by \scite{3.2A}(2) what we have to prove is: assuming
$\varepsilon < \zeta < \mu,k < \omega$ and $\varepsilon_m < \zeta_m < \mu$ for
$m < k$, then $B \models ``b_\zeta - b_\varepsilon -
\dsize \bigcup_{m < k} (b_{\varepsilon_m} - b_{\zeta_m}) \ne 0"$.

Now 
$$
\align
Y =: \biggl\{
(i,n):&g_\varepsilon(i,n) < g_\zeta(i,n) \text{ and}\\
  &g_{\varepsilon_m}(i,n) < g_{\zeta_m}(i,n) \text{ for } m = 0,1,
\dotsc,k-1 \biggr\}
\endalign
$$
\medskip

\noindent
is known to belong to $D^*$, hence it is not empty so let $(i^*,n^*) \in Y$.
Now \newline
$B_{i^*} \models b_{i^*,\xi} \cap a_{i^*,n^*} = 
a_{i^*,n^*,g_\xi(i^*,n^*)}$, for every $\xi < \mu$, in particular for $\xi$
among $\varepsilon,\zeta,\varepsilon_m,\zeta_m$ 
(for $m < k$).  As $(i^*,n^*) \in Y$ we have

$$
\align
B_{i^*} \models (b_\zeta - b_\varepsilon) \cap a_{i^*,n^*} &\ge b_\zeta
\cap a_{i^*,n^*} - b_\varepsilon \cap a_{i^*,n^*} \\
  &= a_{i^*,n^*,g_\zeta(i^*,n^*)} - a_{i^*,n^*,g_\varepsilon(i^*,n^*)}
> 0
\endalign
$$
\medskip

\noindent
(as $g_\zeta(i^*,n^*) > g_\varepsilon(i^*,n^*)$ as
$(i^*,n^*) \in Y$) and similarly \newline
$B_{i^*} \models (b_{\varepsilon_m} - b_{\zeta_m}) \cap a_{i^*,n^*} = 0$.

Hence

$$
B_{i^*} \models ``b_{i^*,\zeta} - b_{i^*,\varepsilon} -
\dsize \bigcup_{m < k} (b_{i^*,\varepsilon_m} - b_{i^*,\zeta_m}) \ne 0".
$$
\medskip

\noindent
As this holds for every $(i^*,n^*) \in Y$ and $Y \in D^*$, by the assumptions
on  $D^*$ we have

$$
\{ i^* < \kappa: B_{i^*} \models ``b_{i^*,\zeta} - b_{i^*,\varepsilon} -
\dsize \bigcup_{m < k} (b_{i^*,\varepsilon_m} - b_{i^*,\zeta_m}) \ne 0"\}
\in D^+
$$
\medskip

\noindent
hence in $B,b_\zeta - b_\varepsilon \notin J$ as required.
\bigskip

\noindent
\underbar{$(\beta)' \Rightarrow (\gamma)$} \newline
Let $B_i$ be the interval Boolean algebra for $\gamma_i$, an ordinal
$< \lambda_i$. 

To prove clause $(\gamma)$ we assume that our regular $\mu$ is 
$< \text{ Depth}^+_h(\dsize \prod_{i < \kappa} B_i/D)$, and we have to find
$n_i < \omega,\lambda_{i,n} < \lambda_i$ for $i < \kappa,n < n_i$ and
$D^*$ as in
the conclusion of clause $(\gamma)$.  So there are $f_\alpha \in \dsize \prod
_{i < \kappa}B_i$ for $\alpha < \mu$ and an ideal $J$ of the Boolean algebra
$B =: \dsize \prod
_{i < \kappa} B_i/D$ such that $f_\alpha/D < f_\beta/D \text{ mod } J$ for
$\alpha < \beta$.

Remember $\mu > 2^\kappa$.  Let $f_\alpha(i) = \dsize
\bigcup_{\ell < n(\alpha,i)}[j_{\alpha,i,2 \ell},j_{\alpha,i,2 \ell +1})$
where $j_{\alpha,i,\ell} < j_{\alpha,i,\ell +1} \le \gamma_i$ for $\ell < 2n
(\alpha,i)$.  As $\mu = \text{ cf}(\mu) > 2^\kappa$, without loss of
generality $n_{\alpha,i} = n_i$.  By \cite[6.6D]{Sh:430} (better see 
\cite[6.1]{Sh:513} or \cite[7.0]{Sh:620}) we can find $A \subseteq 
A^* =: \{(i,\ell):i < \kappa,
\ell < 2n_i\}$ and $\langle \gamma^*_{i,\ell}:i < \kappa,\ell < 2n_i
\rangle$ such that $(i,\ell) \in A \Rightarrow \gamma^*_{i,\ell}$ is a limit
ordinal of cofinality $> 2^\kappa$ and
\medskip
\roster
\item "{$(*)$}"  for every $f \in \dsize \prod_{(i,\ell) \in A} \gamma^*
_{i,\ell}$ and $\alpha < \mu$ there is $\beta \in (\alpha,\mu)$ such that:

$$
\gather
(i,\ell) \in A^* \backslash A \Rightarrow j_{\beta,i,\ell} = \gamma^*
_{i,\ell} \\
(i,\ell) \in A \Rightarrow f(i,\ell) < j_{\beta,i,\ell} < \gamma^*_{i,\ell}
\endgather
$$
\endroster
\medskip

For $(i,\ell) \in A^*$ define $\beta^*_{i,\ell}$ by \newline
$\beta^*_{i,\ell} =: \sup \{\gamma^*_{i,m}:(i,m) \in A^*$ and 
$\gamma^*_{i,m} < \gamma^*_{i,\ell} \text{ and } m < 2n_i$ \nl

$\qquad \qquad \qquad \quad \text{ (actually } m < \ell 
\text{ suffice)}\}$. \newline  
Now $\beta^*_{i,\ell} < \gamma^*_{i,\ell}$ as the supremum is on a finite
set, and the case $0=\beta^*_{i,\ell} = \gamma^*_{i,\ell}$ does not occur if
$(i,\ell) \in A$.  Let

$$
\align
Y = \biggl\{ \alpha < \mu:&\text{ if } (i,\ell) \in A^* \backslash A
\text{ then } j_{\alpha,i,\ell} = \gamma^*_{i,\ell} \\
  &\text{and if } (i,\ell) \in A \text{ then } \beta^*_{i,\ell} < 
j_{\alpha,\ell,i} < \gamma^*_{\ell,i} \biggr\}.
\endalign
$$
\medskip

\noindent
Clearly $\{f_\alpha:\alpha \in Y\}$ satisfy $(*)$, so without loss 
of generality $Y = \mu$.
\medskip

\noindent
Clearly
\medskip
\roster
\item "{$(*)_1$}"  $\langle \gamma^*_{i,\ell}:\ell < 2n_i \rangle$ is non-
decreasing (for each $i$).
\endroster
\medskip

\noindent
Let $u_i = \{ \ell < 2n_i:(\forall m < \ell)[\gamma^*_{i,m} < \gamma^*
_{i,\ell}]\}$.
\mn
For $i < \kappa,\ell < n_i$ define \newline
$b_{i,\ell} =: f_\alpha(i) \cap
[\beta^*_{i,\ell},\gamma^*_{i,\ell}) \in B_i$.  \newline
Let $w_i =: \{ \ell \in u_i:
\text{for every (equivalently some) } \alpha < \mu \text{ we have}$ \newline
$B_i \models ``[\beta^*_{i,\ell},\gamma^*_{i,\ell}) \cap f_\alpha(i)$ is
$\ne 0$ and $\ne [\beta^*_{i,\ell},\gamma^*_{i,\ell})"\}$. \newline
So
\mr
\item "{$(*)_2$}"   $f_\alpha(i) \backslash 
\dsize \bigcup_{\ell \in w_i} b_{i,\ell}$ does not
depend on $\alpha$, call it $c_i(\in B_I)$.
\ermn
Let for $\ell \in w_i$

$$
\align
u_{i,\ell} =: \biggl\{ n < n_i:&[j_{\alpha,i,2n},j_{\alpha,i,2n+1})
\text{ is not disjoint to } [\beta^*_{i,\ell},\gamma^*_{i,\ell}) \\
  &\text{ for some (equivalently every) } \alpha < \mu \biggr\}.
\endalign
$$

$$
\align
A_0 = \biggl\{ (i,\ell):&i < \kappa,\ell \in w_i \text{ and for some }
n \in u_{i,\ell} \text{ we have, for some} \\
  &(\equiv \text{ every) } \alpha < \mu \text{ that } j_{\alpha,i,2n} \le 
  \beta^*_{i,\ell} < j_{\alpha,i,2n+1} < \gamma^*_{i,\ell} \biggr\}.
\endalign
$$

$$
\align
A_1 = \biggl\{(i,\ell):&i < \kappa,\ell \in w_i \text{ and for some }
n \in u_{i,\ell} \text{ we have, for some}\\
  &(\equiv \text{ every) } \alpha < \mu \text{ that } j_{\alpha,i,2n} < 
\gamma^*_{i,\ell} \le j_{\alpha,i,2n+1} \biggr\}.
\endalign
$$
\mn
Let

$$
b^0_i =: \bigcup \biggl\{ 
[\beta^*_{i,\ell},\gamma^*_{i,\ell}):\ell \in w_i \text{ and }
(i,\ell) \in A_0 \biggr\} \in B_i
$$

$$
b^1_i =: \bigcup \biggl\{ 
[\beta^*_{i,\ell},\gamma^*_{i,\ell}):\ell \in w_i \text{ and }
(i,\ell) \in A_1 \biggr\} \in B_i
$$

$$
c^1_i = b^0_i \cap b^1_i,c^2_i = b^0_i \cap (1-b^1_i),
c^3_i = b(1-b^0_i) \cap b^1_i,c^4_i = (1-b^0_i) \cap (1-b^1_i)
$$

$$
b_0 =: \langle b^0_i:i < \kappa \rangle/D \in B
$$

$$
b_1 =: \langle b^1_i:i < \kappa \rangle/D \in B
$$

$$
c_t = \langle c^t_i:i < \kappa \rangle/D \in B, c = \langle c_i:i < \kappa
\rangle/D \in B.
$$
\mn
Let $J_2 = \{b \in B:\langle (f_\alpha/D) \cap b:\alpha < \mu \rangle$ is
eventually constant modulo $J,(\exists \alpha < \mu)(\forall \beta)[\alpha
\le \beta < \mu \rightarrow (f_\alpha/D) \cap b - (f_\beta/D) \cap b 
\in J]\}$.  So $B \models c \le f_\alpha/D$.

Clearly $J_1$ is an ideal of $B$ extending $J$ and $1_B \notin B$.  Also if
$x \in J^+_1$ then for some closed unbounded $E \subseteq \mu$ we have:
$\langle (f_\alpha/D) \cap x:\alpha \in E \rangle$ is strictly increasing
modulo $J$.
\sn
Hence by easy manipulations without loss of generality:
\medskip
\roster
\item "{$(*)_3(a)$}"  if $c = c_t$ and $c \in J^+_1$ then
$\langle (f_\alpha/D) \cap c:\alpha < \mu \rangle$ is strictly increasing
modulo $J$
\sn
\item "{$(b)$}"  for at least one $t,c_t \in J^+_1$.
\ermn
By $(*)$ we can find $0 < \alpha_0 < \alpha_1 < \alpha_2 < \mu$
such that:
\medskip
\roster
\item "{$(*)_4$}"  if $i < \kappa,\ell < 2n_i,\dsize \bigwedge_{\alpha < \mu}
\gamma^*_{i,\ell} > j_{\alpha,i,\ell}$ and $k < 2$ then \newline
$\sup\{j_{\alpha_k,i,\ell_1}:j_{\alpha_k,i,\ell_1} < \gamma^*_{i,\ell}$
and $\ell_1 < 2n_i \} < j_{\alpha_{k+1},i,\ell}$.
\endroster
\medskip

\noindent
Now if in $(*)_3, c_4 \in J^+_1$ occurs then

$$
\align
B_i \models ``f_{\alpha_0}(i) \cap f_{\alpha_1}(i) \cap c^4_i 
= c_i &= \dsize \bigcup_{\ell \in w_i}(f_{\alpha_0}(i) \cap f_{\alpha_1}(i))
\cap [\beta^*_{i,\ell},\gamma^*_{i,\ell}) \\
  &= \dsize \bigcup_{\ell \in w_i} 0_{B_i} = 0_{B_i}"
\endalign
$$
\medskip

\noindent
(as for each $\ell \in w_i$ the intersection is the intersection of two unions
of intervals which are pairwise disjoint) 
whereas we know $(f_{\alpha_0}/D) \cap (f_{\alpha_1}/D) \cap c^4_i - c_i 
=_J f_{\alpha_0}/D \cap c_4 - c \notin J$; contradiction.
\medskip

Next if in $(*)_3, c_2 \in J^+_1$ holds then
$$
\align
B_i \models ``(f_{\alpha_1}(i) \cap c^4_i = c_i) - (f_{\alpha_0}(i) \cap
  c^4_i - c_i) &=\dsize \bigcup_{\ell \in w_i}[f_{\alpha_1}(i) \cap 
  [\beta^*_{i,\ell},\gamma^*_{i,\ell}) - 
  f_{\alpha_0}(i) \cap [\beta^*_{i,\ell},\gamma^*_{i,j}]] \\
  &= \dsize \bigcup_{\ell \in w_i} 0_{B_i} = 0_{B_i}"
\endalign
$$
\medskip

\noindent
(as for each $\ell \in w_i$ the term is the difference of two unions of 
intervals but the first is included in the right most interval of the second)
and we have a contradiction. \newline
Now if in $(*)_3, c_1 \in J^+$ holds then

$$
\align
B_i \models ``(f_{\alpha_2}(i) \cap c^1_i - c_i) &- (f_{\alpha_1}(i) \cap
c^1_i - c_i) \cup (f_{\alpha_0}(i) \cap c^1_i - c_i) \\
  &= \dsize \bigcup_{\ell \in w_i}[(f_{\alpha_2}(i) - f_{\alpha_1}(i) \cup 
f_{\alpha_0}(i)) \cap [\beta^*_{i,\ell},\gamma^*_{i,\ell})] \\
  &= \dsize \bigcup_{\ell \in w_i} 0_{B_i} = 0_{B_i}"
\endalign
$$
\medskip

\noindent
and we get a similar contradiction. \newline
So
\medskip
\roster
\item "{$(*)_5$}"  in $(*)_3, c_2 \in J^+_1$.
\endroster
\medskip

\noindent
Without loss of generality
\medskip
\roster
\item "{$(*)_6$}"  for $\alpha < \mu,i < \kappa$ and $\ell < 2n_i$ such
that $(i,\ell) \in A$ we have \newline
$\sup\{j_{2\alpha,i,\ell_1}:\ell_1 < 2n_i \text{ and }
j_{2 \alpha,i,\ell_1} < \gamma^*_{i,\ell}\}
< j_{2 \alpha + 1,i,\ell}$.
\endroster
\medskip

\noindent
Now for every $\alpha < \mu$ define $f'_\alpha \in \dsize \prod_{i < \kappa}
B_i$ by

$$
f'_\alpha(i) = \dsize \bigcup_{\ell \in v_i}[\beta^*_{i,\ell},\text{max}
\{j_{\alpha,i,2n+1}:n \in u_{i,\ell} \}).
$$
\medskip

\noindent
Clearly

$$
B_i \models ``f_{2 \alpha}(i) \cap c^2_i - c_i \le 
f'_\alpha(i) \le f_{2 \alpha + 1}(i) \cap c^2_i - c_i".
$$
\medskip

\noindent
Let $v_i = \{\ell \in w_i:(i,\ell) \in A_0,(i,\ell) \notin A_1\}$, so
$c^2_i = \cup\{[\beta^*_{i,\ell},\gamma^*_{i,\ell}):\ell \in v_i\}$. As
$\ell \in v_i \Rightarrow (i,\ell) \in A_0$ necessarily
\mr
\item "{$(*)_7$}"  if $\ell \in v_i$ then $\ell$ is odd and
$j_{\alpha,i,\ell-1} = \beta^*_{i,\ell} < j_{\alpha,i,2 \ell+1} < 
\gamma^*_{i,\ell}$.
\ermn
Let $Y^* =: \dsize \bigcup_{i < \kappa} (\{i\} \times v_i)$ and we
shall define now a family $D_0$ of subsets of $Y^*$.
\newline
For $Y \subseteq Y^*$, and for $\alpha < \mu$ define $f_{\alpha,Y} \in
\dsize \prod_{i < \kappa} B_i$ by

$$
f_{\alpha,Y}(i) = \cup\{[j_{\alpha,i,2 \ell},j_{\alpha,i,2 \ell + 1}):
\ell \in v_i \text{ and } (i,\ell) \notin Y\}.
$$
\medskip

\noindent
For $g \in G =: \dsize \prod_{(i,\ell) \in Y^*} [\beta^*_{i,\ell},\gamma^*
_{i,\ell})$ define $f_g \in \dsize \prod_{i < \kappa} B_i$ by
$f_g(i) = \dsize \bigcup_{\ell \in v_i}[\beta^*_{i,\ell},g((i,\ell)))$, now
\medskip
\roster
\item "{$(*)_8$}"  for every $\alpha < \mu$ for some $g = g^*_\alpha \in G$
we have $f'_\alpha = f_g$ \newline
[why? by the previous analysis in particular $(*)_7$].
\endroster
\medskip

\noindent
Let

$$
D_0 = \{ Y \subseteq Y^*:\text{ for some } Z \in D \text{ we have }
\dsize \bigcup_{i \in Z} (\{i\} \times v_i) \subseteq Y\};
$$
\medskip

\noindent
it is a filter on $Y^*$.
\medskip
\roster
\item "{$(*)_9$}"  if $g_1,g_2 \in G$ then
{\roster
\itemitem{ (a) }  $g_1 \le_{D_0} g_2 \Leftrightarrow B \models
(f_{g_1}/D) \cap c_2 \le (f_{g_2}/D) \cap c_2$
\sn
\itemitem{ (b) }  $g_1 <_{D_0} g_2 \Leftrightarrow B \models
(f_{g_1}/D) \cap c_2 < (f_{g_2}/D) \cap c_2$
\endroster}
\smallskip
\noindent
\item "{$(*)_{10}$}"  for every $g' \in G$ for some $\alpha(g') < \mu$
we have $g' < g^*_{\alpha(g')}$ (see $(*)_8$) \newline
[why?  by $(*)$].
\endroster
\medskip

\noindent
Now
\medskip
\roster
\item "{$\bigotimes$}"  cf$(\dsize \prod_{(i,\ell) \in Y^*} \gamma^*_{i,\ell}
/D_0) \ge \mu$.
\newline
[why?  if not, we can find $G^* \subseteq G = \dsize \prod_{(i,\ell) \in Y^*}
[\beta^*_{i,\ell},\gamma^*_{i,\ell})$ of cardinality $< \mu$, cofinal in
$\dsize \prod_{(i,\ell) \in Y^*} \gamma^*_{i,\ell}/D_0$.
For each $g \in G^*$ for some $\alpha(g) < \mu$ we have $g < g^*_{\alpha(g)}$,
hence $\alpha \in [\alpha(g),\mu) \Rightarrow g <_{D_0} g^*_\alpha$, let
$\alpha(*) = \sup\{\alpha(g):g \in G\}$ so $\alpha(*) < \mu$ so
$\dsize \bigwedge_{g \in G} g <_{D_0} g^*_{\alpha(*)}$; contradiction, so
$\bigotimes$ holds].
\endroster
\medskip

\noindent
So for some ultrafilter $D^*$ on $Y^*$ extending $D_0,\mu \le \text{ tcf}
\left( \dsize \prod_{(i,\ell) \in Y^*} \gamma^*_{i,\ell}/D^* \right)$, hence
$\mu \le \text{ tcf } \dsize \prod_{(i,\ell) \in Y^*} \text{ cf}
(\gamma^*_{i,\ell})/D^*$ and by \cite[II,1.3]{Sh:g} for some \newline
$\lambda'_{i,\ell} = \text{ cf}(\lambda'_{i,\ell}) \le \text{ cf}
(\gamma^*_{i,\ell}) \le \gamma_i < \lambda_i$ we have
$\mu = \text{ tcf} \left( 
\dsize \prod_{(i,\ell) \in Y^*} \lambda'_{i,\ell}/D^* \right)$ as
required (we could, instead of relying on this quotation, analyze more).
\bigskip

\noindent
\underbar{$(\beta)^- \Rightarrow (\beta)'$} \newline

Let $B_{i,\gamma}$ be the interval Boolean algebra on $\gamma$ for
$\gamma < \lambda_i,i < \kappa$, and we let $B^*_{i,\gamma}$ be generated by
$\{ a^{i,\gamma}_j:j < \gamma\}$ freely except $a^{i,\gamma}_{j_1} \le
a^{i,\gamma}_{j_2}$ for $j_1 < j_2 < \gamma_i$. \newline
So without loss of generality $B_i$ is the disjoint sum of $\{ B^*_{i,\gamma}:
\gamma < \lambda\}$.  Let \newline
$e_{i,\gamma} = 1_{B_{i,\gamma}}$ (so $\langle e_{i,\gamma}:\gamma <
\lambda_i \rangle$ is a maximal antichain of $B_i,B_i \restriction \{x \in
B_i:x \le e_{i,\gamma}\}$ is isomorphic to $B_{i,\gamma}$ and $B_i$ is
generated by $\{x:(\exists \gamma < \lambda_i)(x \le e_{i,\gamma})\}$.  
Let $\langle f_\alpha:\alpha < \mu \rangle$ and an ideal $J$ of $B$ exemplify
clause $(\beta)^-$.

Let $I_i$ be the ideal of $B_i$ generated by $\{ e_{i,\gamma}:\gamma <
\lambda_i\}$, so it is a \ub{maximal ideal}; let $I$ be such that
$(B,I) = \dsize \prod_{i < \kappa}(B_i,I_i)/D$ so clearly $|B/I| =
|2^\kappa/D| \le 2^\kappa < \text{ cf}(\mu)$ (actually $|B/I| = 2$), 
so without loss of generality $\alpha < \beta
\le \mu \Rightarrow f_\alpha/D = f_\beta/D \text{ mod }I$.  We can use
$\langle f_{1 + \alpha}/D - f_0/D:\alpha < \mu \rangle$, so without loss
of generality $f_\alpha/D \in I$, hence without loss of generality $f_\alpha
(i) \in I_i$ for $\alpha < \mu,i < \kappa$.

Let $f_\alpha(i) = \tau_{\alpha,i}(\ldots,e_{i,\gamma(\alpha,i,\varepsilon)},
a^{i,\gamma(\alpha,i,\varepsilon)}
_{j(\alpha,i,\varepsilon)},\dots)_{\varepsilon < n_{\alpha,i}}$ where
$n_{\alpha,i} < \omega$ and $\tau_{\alpha,i}$ is a Boolean term.  As $\mu$ is
regular $> 2^\kappa$, without loss of generality $\tau_{\alpha,i} = \tau_i$
and $n_{\alpha,i} = n_i$.  Let $\gamma^0_{\alpha,i,\varepsilon} =
\gamma(\alpha,i,\varepsilon)$ and $\gamma^1_{\alpha,i,\varepsilon} =
j(\alpha,i,\ell)$.

By \cite[6.6D]{Sh:430} (or better \cite[6.1]{Sh:513}) we can find a subset
$A$ of \newline
$A^* = \{(i,\varepsilon,\ell):i < \kappa \text{ and } \varepsilon <
n_i \text{ and } \ell < 2\}$ and \newline
$\langle \gamma^*_{i,n,\ell}:i < \kappa \text{ and } n < n_i \text{ and }
\ell < 2 \rangle$ such that:
\medskip
\roster
\item "{$(*)(A)$}"  $(i,\varepsilon,\ell) \in A \Rightarrow \text{ cf}
(\gamma^*_{i,\varepsilon,\ell}) > 2^\kappa$
\sn
\item "{${}(B)$}"  for every $g \in \dsize \prod_{(i,n,\ell) \in A}
\gamma^*_{i,n,\ell}$ for arbitrarily large $\alpha < \mu$ we have
\smallskip
$$
(i,n,\ell) \in A^* \backslash A \Rightarrow \gamma^\ell_{\alpha,i,n} =
\gamma^*_{i,n,\ell}
$$

$$
(i,n,\ell) \in A \Rightarrow g(i,n,\ell) < \gamma^\ell_{\alpha,i,n} <
\gamma^*_{i,n,\ell}.
$$
\endroster
\medskip

\noindent
Let

$$
\beta^*_{i,n,\ell} = \sup\{ \gamma^*_{i,n',\ell'}:n' < n,\ell' < 2
\text{ and } \gamma^*_{i,n',\ell'} < \gamma^*_{i,n,\ell} \}
$$
\medskip

\noindent
(can use $n < n_i$).  Without loss of generality

$$
(i,n,\ell) \in A^* \and \alpha < \mu \Rightarrow \gamma^\ell_{\alpha,i,n} \in
(\beta^*_{i,n,\ell},\gamma^*_{i,n,\ell})
$$

$$
(i,n,\ell) \in A^* \backslash A \and \alpha < \mu \Rightarrow 
\gamma^\ell_{\alpha,i,n} = \gamma^*_{i,n,\ell}.
$$
\medskip

\noindent
Also without loss of generality
\medskip
\roster
\item "{$(*)$}"  for $\alpha < \mu$ and $(i,n,\ell) \in A$ we have
$$
\align
\gamma^\ell_{2 \alpha + 1,i,n} > \sup \biggl\{
\gamma^{\ell'}_{2 \alpha,i,n'}:&i < \kappa,\ell' < 2,n' < n_i, \\
  &\text{and } \gamma^{\ell'}_{2 \alpha,i,n'} < \gamma^*_{i,n,\ell} \biggr\}.
\endalign
$$
\endroster
\medskip

\noindent
Let $\triangle_i = \{ \gamma^*_{i,n,0}:n < n_i \text{ and } (i,n,0) \in A^*
\backslash A\}$ and

$$
B'_i = B_i \restriction \sum\{e_{i,\gamma}:\gamma \in \triangle_i\}.
$$
\medskip

\noindent
We define $f'_\alpha \in \dsize \prod_{i < \kappa} B'_i$ by
$f'_\alpha(i) = f_{2 \alpha + 1}(i) \cap (\dsize \bigcup_{\gamma \in
\triangle_i} e_{i,\gamma}) \in B'_i \subseteq B_i$. \newline
Now easily $f'_\alpha/D \le f_{2 \alpha + 1}/D$ and (in B)
$f_{2 \alpha}/D - f'_\alpha/D \le f_{2 \alpha}/D - f'_{2 \alpha + 1}/D \in
J$, hence modulo $J$, also $\langle f'_\alpha:\alpha < \lambda \rangle$ is
increasing.  So 
$\langle B'_i:i < \kappa \rangle,\langle f'_\alpha:\alpha < \mu \rangle$
form a witness, too.  But $B'_i$ is isomorphic to the interval Boolean algebra
of the ordinal $\gamma_i = \dsize \sum_{\gamma \in \triangle_i} \gamma <
\lambda$, so we are almost done.  Well, $\gamma_i$ is an ordinal, not
necessarily a cardinal, but in the proof of 
$(\beta) \Rightarrow (\gamma)$ we allow the $\gamma_i$ to be ordinals
$< \lambda_i$ and $(\gamma) \Rightarrow (\beta)$ was proved too.  
\hfill$\square_{\scite{3.4}}$
\newpage

\head {\S4 On Existence of independent sets (for stable theories)} \endhead  \resetall  
\bigskip

The following is motivated by questions of Bays \cite{Bay} which 
continues some investigations of \cite{Sh:a} (better see \cite{Sh:c}) 
dealing with questions on Pr$_T(\mu)$,Pr$^*_T$ for stable $T$ (see 
Definition \scite{4.2} below).  We connect this to pcf, using 
\cite[3.17]{Sh:430} and also \cite[6.12]{Sh:513}).  We assume basic knowledge
on non-forking (see \cite[Ch.III,I]{Sh:c}) and we say some things on the
combinatorics but the rest of the paper does not depend on this section.
\bigskip

\proclaim{\stag{4.1} Claim}  Assume $\lambda > \theta \ge \kappa$ are regular
uncountable.  \underbar{Then} the following are equivalent:
\medskip
\roster
\item "{$(B)$}"  If $\mu < \lambda$ and $a_\alpha \in [\mu]^{< \kappa}$ for
$\alpha < \lambda$ \underbar{then} for some \newline
$A \in [\lambda]^\lambda$ we have
$\dsize \bigcup_{\alpha \in A} a_\alpha$ has cardinality $< \theta$
\sn
\item "{$(C)$}"  if $\delta = \text{ cf}(\delta) < \kappa$ and $\eta_\alpha
\in {}^\delta \lambda$ for $\alpha < \lambda$ and \newline
$|\{ \eta_\alpha \restriction i:\alpha < \lambda,i < \delta\}| < \lambda$ 
\underbar{then} for some $A \in [\lambda]^\lambda$ \newline
the set $\{ \eta_\alpha \restriction i:\alpha \in A,i <
\delta \}$ has cardinality $< \theta$.
\endroster
\endproclaim
\bigskip

\demo{Proof}  \underbar{$(B) \Rightarrow (C)$}.  Immediate.
\enddemo
\bigskip

\noindent
\underbar{$\neg(B) \Rightarrow \neg(C)$}
\bigskip

\noindent
\underbar{Case 1}:  
For some $\mu \in (\theta,\lambda),\text{cf}(\mu) < \kappa$ and $\text{pp}
(\mu) \ge \lambda$.  Without loss of generality $\mu$ is minimal. \newline
So
\medskip
\roster
\item "{$(*)$}"  ${\frak a} \subseteq \text{ Reg } \cap \mu \backslash
\theta,|{\frak a}| < \kappa,\sup({\frak a}) < \mu \Rightarrow \text{ max pcf}
({\frak a}) < \mu$.
\endroster
\bigskip

\noindent
\underbar{Subcase 1a}:  $\lambda < \text{ pp}^+(\mu)$. \newline

So by \cite[Ch.VII,1.6(2)]{Sh:g}, (if cf$(\mu) > \aleph_0$) and 
\cite[6.5;]{Sh:430} (if cf$(\mu) = \aleph_0$) we can find $\langle 
\lambda_\alpha:\alpha < \text{ cf}(\mu) \rangle$, a strictly increasing 
sequence of regulars from $(\theta,\mu)$ with limit
$\mu$ and ideal $J$ on cf$(\mu)$ satisfying $J^{bd}_{\text{cf}(\mu)} 
\subseteq J$ such that $\lambda = \text{ tcf}\left( 
\dsize \prod_{\alpha < \text{ cf}(\mu)}
\lambda_\alpha /J \right)$ and $\text{max pcf}\{ \lambda_\beta
:\beta < \alpha\} < \lambda_\alpha$.  By \cite[II,3.5]{Sh:g}, 
there is $\langle f_\zeta:\zeta < \lambda \rangle$ which is $<_J$-increasing 
cofinal in
$\dsize \prod_{\alpha < \text{ cf}(\mu)} \lambda_\alpha/J$ with
$|\{f_\zeta \restriction \alpha:\zeta < \lambda\}| < \lambda_\alpha$. \newline
Easily $\langle f_\zeta:\zeta < \lambda \rangle$ exemplifies
$\neg(C):\text{if } A \in [\lambda]^\lambda$ and 
$B =: \dsize \bigcup_{\zeta \in A}
\text{ Range}(f_\zeta)$ has cardinality $< \mu$ let $g \in \dsize \prod_\alpha
\lambda_\alpha$ be: $g(\alpha) = \sup(\lambda_\alpha \cap B)$ if 
$< \lambda_\alpha$, zero otherwise and let $\alpha_0 = \text{ Min}\{\alpha <
\text{ cf}(\mu):\lambda_\alpha > |B|\}$.  So $\alpha_0 < \text{ cf}(\mu)$ and
$\zeta \in A \Rightarrow f_\zeta \restriction [\alpha_0,
\text{cf}(\mu)) < g$, contradiction to $``<_J$-cofinal".
\bigskip

\noindent
\underbar{Subcase 1b}:  cf$(\mu) > \aleph_0$ and 
$\text{pp}^+(\mu) = \text{ pp}(\mu) = \lambda$.  Use
\cite[\S6]{Sh:513} and finish as above.
\bigskip

\noindent
\underbar{Subcase 1c}:  $\text{cf}(\mu) = \aleph_0$ and $\lambda = 
\text{ pp}^+(\mu) = \text{ pp}(\mu) = \lambda$.

Let ${\frak a},\langle {\frak b}_\tau:\tau \in R \rangle,\langle f_\tau:
\tau \in R \rangle$ be as in \cite[6.x]{Sh:513}, so $|{\frak b}_\tau| =
\aleph_0$.  Let $\eta_\tau$ be an $\omega$-sequence of ordinals enumerating
$\text{Rang}(f_\tau)$ for $\tau \in R$, now $\{ \eta_\tau:\tau \in R\}$ is
as required.
\bigskip

\noindent
\underbar{Case 2}:  Not Case 1.

So by \cite[Ch.II,5.4]{Sh:g}, we have $\theta \le \mu < \lambda \Rightarrow
\text{ cov}(\mu,\theta,\kappa,\aleph_1) < \lambda$. \newline
As we are assuming $\neg(B)$, we can find $\mu_0 < \lambda,a_\alpha \in 
[\mu_0]^{< \kappa}$
for $\alpha < \lambda$ such that $A \in [\lambda]^\lambda \Rightarrow 
| \dsize \bigcup_{\alpha \in A} a_\alpha| \ge \theta$, but by the 
previous sentence we can find
$\mu_1 < \lambda$ and $\{ b_\beta:\beta < \mu_1\} \subseteq 
[\mu_0]^{< \theta}$ such that: every $a \in [\mu_0]^{< \kappa}$ is
included in the union of $\le \aleph_0$ sets from $\{ b_\beta:\beta < \mu_1
\}$.  So we can find $c_\alpha \in [\mu_1]^{\aleph_0}$ for $\alpha < \lambda$
such that $a_\alpha \subseteq \dsize \bigcup_{\beta \in c_\alpha} b_\beta$.
Now for $A \in [\lambda]^\lambda$, if $| \dsize \bigcup_{\alpha \in A}
c_\alpha| < \theta$ then

$$
\align
| \bigcup \{ a_\alpha:\alpha \in A\}| &\subseteq| \bigcup
\{ \dsize \bigcup_{\beta \in c_\alpha} b_\beta:\alpha \in A\}| \\
  &= |\bigcup \{ b_\beta:\beta \in \dsize \bigcup_{\alpha \in A} c_\alpha\}| <
\text{ min}\{\sigma:\sigma = \text{ cf}(\sigma) > |b_\beta| \text{ for }
\beta < \mu_1\} \\
  &+ | \dsize \bigcup_{\alpha \in A} c_\alpha|^+ \le \theta + \theta = \theta
\endalign
$$
\medskip

\noindent
contradicting the choice of
$\langle a_\alpha:\alpha < \lambda \rangle$.
\medskip

So 
\medskip
\roster
\item "{$(*)$}"  $c_\alpha \in [\mu_1]^{\le \aleph_0}$, for $\alpha < \lambda,
\mu_1 < \lambda$ and \newline
$A \in [\lambda]^\lambda \Rightarrow | \dsize \bigcup_{\alpha \in A} 
c_\alpha| \ge \theta$.
\endroster
\medskip

Let $\eta_\alpha$ be an $\omega$-sequence enumerating $c_\alpha$, so
$\langle \eta_\alpha:\alpha < \lambda \rangle$ is a counterexample to
clause $(C)$. \hfill$\square_{\scite{4.1}}$
\bigskip 

\noindent
We concentrate below on $\lambda,\theta,\kappa$ regular
(others can be reduced to it).
\definition{\stag{4.2} Definition}  Let $\bold T$ be a 
complete first order theory; which
is stable (${\frak C}$ the monster model of $T$ and $A,B,\ldots$ denote 
subsets of ${\frak C}^{\text{eq}}$ of cardinality 
$<\|{\frak C}^{\text{eq}}\|$). \newline
1) $\text{Pr}_{\bold T}(\lambda,\chi,\theta)$ means:
\medskip
\roster
\item "{$(*)$}"  if $A \subseteq {\frak C}^{\text{eq}},|A| = \lambda$ then we
can find $A' \subseteq A,|A'| = \chi$ and $B',|B'| < \theta$ such that
$A'$ is independent over $B'$ \nl
(i.e. $a \in A' \Rightarrow \text{ tp}(a,B' \cup (A' \backslash \{a\}))$ does 
not fork over $B'$).
\endroster
\medskip

\noindent
2) $\text{Pr}^*_{\bold T}(\lambda,\mu,\chi,\theta)$ means:
\medskip
\roster
\item "{$(**)$}"  if $A \subseteq {\frak C}^{\text{eq}}$ is independent over
$B$ where $|A| = \lambda$ and $|B| < \mu,B \subseteq {\frak C}^{\text{eq}}$ 
\underbar{then} there are $A' \subseteq A,|A'| 
= \chi$ and $B' \subseteq B$ satisfying $|B'| < \theta$ such that 
$\text{tp}(A',B)$ does not fork over $B'$ (hence $A'$ is independent over
$B'$).
\endroster
\medskip

\noindent
3) $\text{Pr}^*_{\bold T}(\lambda,\chi,\theta)$ means
$\text{Pr}^*_{\bold T}(\lambda,\lambda,\chi,\theta)$.  
\enddefinition
\bigskip

\demo{\stag{4.3} Fact}  Assume $\lambda$ is regular 
$> \theta \ge \kappa_r(\bold T)$ \underbar{then}
\mr
\item  if $\chi = \lambda$ then 
$\text{Pr}_{\bold T}(\lambda,\chi,\theta) \Leftrightarrow 
\text{Pr}^*_{\bold T}(\lambda,\lambda,\chi,\theta)$
\sn
\item  if $\lambda \ge \chi \ge \mu \ge \theta$ then 
$\text{Pr}_{\bold T}(\lambda,\chi,\theta) \Rightarrow \text{Pr}^*_{\bold T}
(\lambda,\mu,\chi,\theta)$.
\endroster
\enddemo
\bigskip

\demo{Proof}  1) The direction $\Leftarrow$ is by the proof in 
\cite[III]{Sh:a}. \newline
[Let $A,B$ be given (the $B$ is not really necessary), such that 
$\lambda = |A| > |B| + \kappa_r(\bold T)$
so let $A = \{ a_i:i < \lambda\}$; define \newline
$A_i =: \{a_j:j < i\}, S = \{ i < \lambda:\text{cf}(i) \ge 
\kappa_r(\bold T)\}$, so by the definition of $\kappa_r(\bold T)$ for 
$\alpha \in S$ there is $j_\alpha < \alpha$ such that 
tp$(a_\alpha,A_\alpha \cup B)$ does not fork over
$A_{j_\alpha} \cup B$ so for some $j^*$ the set 
$S' = \{ \delta \in S:j_\delta =
j^*\}$ is stationary, now apply the right side with $\{a_\delta:\delta \in
S'\},A_{j^*} \cup B$, here standing for $A,B$ there].

The other direction $\Rightarrow$ follows by part (2). \nl
2) This is easy, too, by the non-forking calculus
\cite[III,Th.0.1 + (0)-(4),pgs.82-84]{Sh:a} or just read the proof. 
Let $A \subseteq {\frak C}^{\text{eq}}$ be independent over $B$, where
$|A| = \lambda$ and $|B| < \mu$.  As we are assuming
Pr$_{\bold T}(\lambda,\chi,\theta)$ there is $A' \subseteq A,|A'| = \chi$
and $B',|B'| < \theta$ such that $A'$ is independent over $B'$.   So for
every finite $\bar c \subseteq B$ for some $A_{\bar c} \subseteq A'$ of
cardinality $< \kappa(\bold T) \, (\le \kappa_r(\bold T))$ we have:
$A' \backslash A_{\bar c}$ is independent over $B' \cup \bar c$.  So
$A^* = \bigcup \{A_{\bar c}:\bar c \subseteq B \text{ finite}\}$ has
cardinality $< \kappa_r(\bold T) + |B|^+ \le \chi$ so necessarily 
$A' \backslash A^*$ has cardinality $\chi$ and it is independent 
over $\cup\{\bar c:\bar c \subseteq B \text{ finite}\} \cup B' = B \cup B'$.]
\nl
${{}}$  \hfill$\square_{\scite{4.3}}$
\enddemo
\bigskip

\demo{\stag{4.4} Discussion}  So it 
suffices to prove the equivalence \newline
$\text{Pr}^*_{\bold T}(\lambda,\mu,\chi,\theta) \Leftrightarrow \text{Pr}
(\lambda,\mu,\chi,\theta,\kappa)$ with $\kappa = \kappa_r(\bold T)$, where
\enddemo
\bigskip

\definition{\stag{4.5} Definition}  Assume
\medskip
\roster
\item "{$(*)$}"  $\lambda \ge \text{ max}\{ \mu,\chi\} \ge \text{ min}
\{ \mu,\chi\} \ge \theta \ge \kappa > \aleph_0$ and $\mu > \theta$ and
for simplicity $\lambda,\theta,\kappa$ are regular if not said otherwise 
(as the general case can be reduced to this case).
\endroster
\medskip

\noindent
1) $\text{Pr}(\lambda,\mu,\chi,\theta,\kappa)$ is defined as follows:
if $a_\alpha \in [\mu]^{< \kappa}$ for $\alpha < \lambda$ and
$| \dsize \bigcup_{\alpha < \lambda} a_\alpha| < \mu$
\underbar{then} there is $Y \in [\lambda]^\chi$ such that
$|\dsize \bigcup_{\alpha \in Y} a_\alpha| < \theta$; \newline
2) $\text{Pr}^{\text{tr}}
(\lambda,\mu,\chi,\theta,\kappa)$ is defined similarly but for some tree
$T$ each $a_\alpha$ is a branch nl
of $T$. \newline
3) We write Pr$(\lambda,\le \mu,\chi,\theta,\kappa)$ for Pr$(\lambda,\mu^+,
\chi,\theta,\kappa)$ and similarly for Pr$^{\text{tr}}$ and Pr$^*_{\bold T}$.
\enddefinition
\bigskip

\demo{\stag{4.6} Fact}  Assume $\lambda,\mu = \chi,\chi,\theta,\kappa =
\kappa_i(\bold T)$ satisfies $(*)$ of Definition \scite{4.5}.  Then \nl
1)  $\text{Pr}(\lambda,\mu,\chi,\theta,
\kappa_r(\bold T)) \Rightarrow \text{Pr}^*_{\bold T}(\lambda,\mu,\chi,\theta) 
\Rightarrow \text{Pr}^{\text{tr}} (\lambda,\mu,\chi,\theta,
\kappa_r(\bold T))$. \nl
2) Pr$(\lambda,\chi,\chi,\theta,\kappa_r(\bold T)) \Rightarrow
\text{ Pr}_{\bold T}(\lambda,\chi,\theta) \Rightarrow \text{ Pr}^{\text{tr}}
(\lambda,\chi,\chi,\theta,\kappa_r(\bold T))$. \nl
3) We have obvious monotonicity properties.
\enddemo
\bigskip

\demo{Proof}  Straight.  \nl
1) First we prove the first implication so assume Pr$(\lambda,\mu,\chi,
\theta,\kappa_i(\bold T))$, let $\kappa = \kappa_i(\bold T)$, hence $(*)$
of \scite{4.5} holds and we shall prove Pr$^*_{\bold T}(\lambda,\mu,\chi,
\theta)$.  So (see Definition \scite{4.2}(2)) we have $A \subseteq
{\frak C}^{\text{eq}}$ is independent over $B \subseteq {\frak C}^{\text{eq}},
|A|= \lambda$ and $|B| < \mu$.  Let $A = \{a_\alpha:\alpha < \lambda\}$, with
no repetitions of course and $B = \{b_j:j < j(*)\}$ so $i(*) < \mu$. For each
$\alpha < \lambda$, there is a subset $u_\alpha$ of $j(*)$ of cardinality
$< \kappa_r(\bold T) = \kappa$ such that tp$(a_\alpha,B)$ does not fork over
$\{b_j:j \in u_\alpha\}$.  So $u_i \in [\mu]^{< \kappa}$ and
$| \dbcu_{\alpha < \lambda} u_\alpha| \le |j(*)| < \mu$ hence as we are
assuming Pr$(\lambda,\mu,\chi,\theta,\kappa)$, there is $Y \in [\lambda]^\chi$
such that $|\dbcu_{\alpha \in Y} u_\alpha| < \theta$.  Let $B' = \{b_j:j \in
\dbcu_{\alpha \in Y} u_\alpha\},A' = \{a_\alpha:\alpha \in Y\}$ so $B'
\subseteq B,|B'| < \theta$ and $A' \subseteq A,|A'| = \chi$ and by the
nonforking calculus, tp$(A',B)$ does not fork over $B'$ (even $\{a_\alpha:
\alpha \in Y\}$ is independent over $(B,B')$).

Second, we prove the second implication, so we assume 
Pr$^*_{\bold T}(\lambda,\mu,\chi,\theta)$ and we shal prove
Pr$^{\text{tr}}(\lambda,\mu,\chi,\theta,\kappa_r(\bold T))$.  Let 
$\kappa = \kappa_r(\bold T)$.
\mn
Let $T$ be a tree and for $\alpha < \lambda,a_\alpha$ a branch, $|a| <
\kappa,|\dbcu_{\alpha < \lambda} a_\alpha| < \mu$.  Without loss of generality
$T = \dbcu_{\alpha < \lambda} a_\alpha,\lambda = \dbcu_{\zeta < \kappa}
A_\zeta$, where $A_\zeta = \{\alpha:\text{otp}(a_\alpha) = \zeta\}$.
Without loss of generality $T \subseteq {}^{\kappa >}\mu,T = \dbcu_{\zeta
< \kappa} T_\zeta$ such that

$$
\eta \in T_\zeta \backslash \{<>\} \Rightarrow \eta(0) = \zeta,
$$
\mn
$T_\zeta$ can be replaced by $\{\eta \restriction C_\zeta:\eta \in T_\zeta\}$
where $0 \in C_\zeta$, otp $C = 1 + \text{ cf}(\zeta),\sup(C) = \zeta$.  So
\wilog \,

$$
T = \cup \{T_\sigma:\sigma \in \text{ Reg } \cap \kappa_T\}
$$

$$
<> \ne \eta \in T_\sigma \Rightarrow \eta(0) = \sigma.
$$
\mn
Without loss of generality $\lambda = \cup \{A_\sigma:\sigma \in
\text{ Reg } \cap \kappa\}$ and $A \dbcu_{\alpha \in A_\sigma} a_\alpha =
T_\sigma$.  It is enough to take care of one $\sigma$ (otherwise little
more work).  So \wilog:

$$
\alpha < \lambda \Rightarrow \text{ otp}(a_\alpha) = \sigma.
$$
\mn
As $\sigma = \text{ cf}(\sigma) < \kappa$ there are $A_i \subseteq
{\frak C}^{\text{eq}}$ such that $\langle A_i:i \le \sigma \rangle$
increase continuously and $p \in S(A_\sigma) \dsize \bigwedge_{i < \sigma}
p \restriction A_{i+1}$ forks over $A_i$ say $\varphi(x,c_i) \in p
\restriction A_{i+1}$ forks over $A_i$ and $A_i = \{c_j:j < i\}$.

By nonforking calculus we can find $\langle f_\eta:\eta \in T \rangle,
f_\eta$ elementary mapping

$$
\text{Dom}(f_\eta) = A_{\ell g(\eta)}
$$
\mn
$\langle f_\eta:\eta \in T \rangle$ nonforking tree, that is

$$
\nu \triangleleft \eta \Rightarrow \text{ tp}(\text{Rang}(f_\eta),\cup
\{\text{Rang}(f_\rho):\rho \in T,\rho \restriction \ell g(\nu) +1
\ntriangleleft \eta\})
$$
\mn
does not fork over $A_\nu$. \nl
For $\alpha < \lambda$, let $g_\alpha = \cup\{f_\nu:\nu \in a_\alpha\},
A_\alpha = \dbcu_{nu \in a_\alpha} \text{ Rang}(f_\nu) = g_\alpha(A_\sigma)$
and $p_\alpha = g_\alpha(p)$. \nl
Let $b_\alpha \in {\frak C}$ realize

$$
A_\eta = \text{ Rang}(f_\eta)
$$

$$
\text{tp}(b_\alpha,\dbcu_{\eta \in T} A_\eta \cup \{b_\beta:\beta \ne
\alpha\}) \text{ does not fork over } A_\alpha.
$$
\mn
Now we apply Pr$^*_{\bold T}(\lambda,\mu,\chi,\theta)$ on

$$
A = \{b_\alpha:\alpha < \lambda\}
$$

$$
B = \dbcu_{\eta \in T} A_\eta.
$$
\mn
So there are $A' \subseteq A,|A'| = \chi$ and $B' \subseteq B,|B'| <
\theta$, tp$(A',B)$ does not fork over $B'$, hence (for some $y \in
[\lambda]^\chi$) we have $A' = \{a_\alpha:\alpha \in Y\}$ independent over
$B'$.  So there is $T' \subseteq T$ subtree such that $|T'| = |B'| + \sigma
< \theta$ such that $B' \subseteq \dbcu_{\rho \in T'} A_\rho$.  Throwing
``few" $(< |B'|^+ + \kappa_r(\bold T))$ members of $A'$ that is of $Y$ we
get $A'$ independent over $B'$. \nl
Easily $Y$ is as required.
\enddemo
\bigskip

\noindent
\underbar{Discussion}  So if Pr,Pr$^{tr}$ are equivalent, $\kappa =
\kappa_r(\bold T)$ then Pr$^*_T$ is equivalent to them (for the suitable
cardinal parameter, so we would like to prove such equivalence).  Now 
\scite{4.1} gives the equivalence when $\theta = \kappa_r(\bold T),
\lambda = \chi =
\text{ cf}(\lambda)$ and ``for every $\mu < \lambda$".  We give below more
general cases; e.g. if $\lambda$ is a successor of regular or
$\{ \delta < \lambda:\text{cf}(\delta) =
\theta^*\} \in I(\lambda)$ or ...
\bigskip

\demo{\stag{4.8} Fact}  Assume $\lambda,\mu,\chi,\theta,\kappa$ as in $(*)$ of
Definition \scite{4.5} and $\mu^* \in [\theta,\mu)$ and \newline
cf$(\mu^*) < \kappa$. \newline
0)  $\text{Pr}(\lambda,\mu,\chi,\theta,\kappa) \Rightarrow \text{ Pr}
^{\text{tr}}(\lambda,\mu,\chi,\theta,\kappa)$. \newline
[Why?  Straight]. \newline
1)  If $\kappa < \lambda$ and $\mu < \lambda$ and cf$(\mu) \ge \kappa$,
\underbar{then} Pr$(\lambda,\le \mu,\chi,\theta,\kappa) \Leftrightarrow 
(\forall \mu_1 < \mu)\text{Pr}(\lambda,\le \mu_1,\chi,\theta,\kappa)$;
similarly for Pr$^{\text{tr}}$. \newline
2)  If $\text{pp}(\mu^*) > \lambda$ \underbar{then}
$\neg \text{ Pr}^{\text{tr}}(\lambda,\mu,\chi,\theta,\kappa)$ (by
\cite[1.5A]{Sh:355}, see \cite[6.10]{Sh:513}). \newline
3)  If $\text{pp}(\mu^*) \ge \lambda$ and
\medskip
\roster
\item "{$(a)$}"  $\{ \delta < \lambda:\text{cf}(\delta) = \theta\} \in 
I[\lambda]$ or just
\sn
\item "{$(a)^-$}"   for some $S \in I[\lambda],(\forall \delta \in S),
\text{cf}(\delta) = \theta$ and
\sn
\item "{$(a)_S$}"   for every closed $e \subseteq \lambda$ of order
type $\chi,e \cap S \ne \emptyset$).
\endroster
\medskip

\noindent
\underbar{Then} $\neg \text{ Pr}^{\text{tr}}(\lambda,\mu,\chi,\theta,\kappa)$.
\newline
[Why? As in \cite[Ch.VIII,6.4]{Sh:g} based on \cite[Ch.II,5.4]{Sh:g} 
better still \cite[Ch.II,3.5]{Sh:g}]. \newline
4)  If $\lambda$ is a successor of regular and $\theta^+ < \lambda$, 
\underbar{then} the assumption (b) of part (3) holds
(see \cite[Ch.VIII,6.1]{Sh:g} based on \cite[\S4]{Sh:351}). \newline
5)  If $\mu < \lambda$ and cov$(\mu,\theta,\kappa,\aleph_1) < \lambda$
(equivalently \newline
$(\forall \tau)[\theta < \tau \le \mu \and \text{ cf}(\tau)
< \kappa \rightarrow \text{ pp}_{\aleph_1\text{-complete}}(\tau) < \lambda)$,
\underbar{then} $\neg \text{Pr}(\lambda,\mu^+,\chi,\theta,\kappa)$ implies 
that for some $\mu_1 \in (\mu,\lambda)$ we have
$\neg \text{Pr}(\lambda,\mu_1,\chi,\theta,\aleph_1)$ (as in Case 2 in the 
proof of \scite{4.1}). \newline
6)  $\text{Pr}(\lambda,\mu,\chi,\theta,\aleph_1) \Leftrightarrow 
\text{ Pr}^{\text{tr}}(\lambda,\mu,\chi,\theta,\aleph_1)$. \newline
7)  Pr$(\lambda,\mu,\lambda,\theta,\kappa)$ \underbar{iff} for every $\tau
\in [\theta,\mu)$ we have: Pr$(\lambda,\le \tau,\lambda,\tau,\kappa)$;
similarly for Pr$^{\text{tr}}$. \newline
8)  Pr$(\lambda,\le \mu,\lambda,\theta,\kappa)$ \underbar{iff} 
Pr$^{\text{tr}}(\lambda,\le \mu,\lambda,\theta,\kappa)$ (by \scite{4.1}).
\enddemo
\bigskip

\proclaim{\stag{4.9} Claim}  Under GCH we get equivalence:
$\text{Pr}(\lambda,\mu,\chi,\theta,\kappa) \Leftrightarrow \text{ Pr}
^{\text{tr}}(\lambda,\mu,\chi,\theta,\kappa)$.
\endproclaim
\bigskip

\demo{Proof}  $\text{Pr} \Rightarrow \text{ Pr}^{\text{tr}}$ is trivial;
so let us prove $\neg \text{ Pr} \Rightarrow \neg \text{ Pr}^{\text{tr}}$, 
so assume
\newline
$\{ a_\alpha:\alpha < \lambda\} \subseteq [\mu]^{< \kappa}$ exemplifies
$\neg \text{Pr}(\lambda,\mu,\chi,\theta,\kappa)$.  Without loss of generality
\newline
$|a_\alpha| = \kappa^* < \kappa$.  By \scite{4.8}(1) without 
loss of generality $\lambda > \mu$, so necessarily
\medskip
\roster
\item "{$(c)$}"  $\lambda = \mu^+,\mu > \kappa^* \ge \text{ cf}(\mu)$ or
\sn
\item "{$(d)$}"  $\lambda = \mu^+,\kappa = \lambda$.
\endroster
\medskip

\noindent
In Case (a) let $T$ be the set of sequences of
bounded subsets of $\mu$ each of cardinality $\le \kappa^*$ of length
$< \text{ Min}\{ \text{cf}(\mu),\kappa^*\}$.  For each $\alpha < \lambda$ let
$\bar b^d = \langle b_{\alpha,\varepsilon}:\varepsilon < \text{ cf}(\mu) 
\rangle$ be a sequence, every initial segment is in $T$ and $a_\alpha = 
\dsize \bigcup_{\varepsilon < \text{ cf}(\mu)} b_{\alpha,\varepsilon}$, so
\newline
$t_\alpha = \{ \bar b^\alpha \restriction \zeta:\zeta < \text{ cf}(\mu)\}$
is a cf$(\mu)$-branch of $T$, and it should be clear.
\enddemo
\bigskip

\remark{\stag{4.10} Remark}  We can get independence result: by 
instances of Chang's Conjecture (so the consistency 
strength seems somewhat more than huge 
cardinals, see Foreman \cite{For}, Levinski Magidor Shelah \cite{LMSh:198}).
\endremark
\newpage

\head {\S5 Cardinal invariants for general regular cardinals: restriction
on the depth} \endhead  \resetall
\bigskip

Cummings and Shelah \cite{CuSh:541} prove that there are no non-trivial
restrictions on some cardinal invariants like ${\frak b}$ and ${\frak d}$, 
even for all regular cardinals simultaneously; i.e., on functions like 
$\langle {\frak b}_\lambda:\lambda \in \text{ Reg} \rangle$.  But not 
everything is
independent of ZFC.  Consider the cardinal invariants ${\frak d}
{\frak p}^{\ell +}_\lambda$, defined below.
\definition{\stag{5.1} Definition}  1) For an ideal $J$ on a regular cardinal 
$\lambda$ let 
\sn
if $\lambda > \aleph_0$

$$
\align
{\frak d}{\frak p}^{1+}_\lambda = 
\text{ Min} \biggl\{ \theta:&\text{ there is no sequence } \langle C_\alpha:
\alpha < \theta \rangle \text{ such that:} \\
  &(a) \quad C_\alpha \text{ is a club of } \lambda, \\
  &(b) \quad \alpha < \beta \Rightarrow |C_\alpha \backslash C_\beta| <
\lambda, \\
  &(c) \quad C_{\alpha + 1} \subseteq \text{ acc}(C_\alpha) \biggr\}.
\endalign
$$
\mn
where acc$(C)$ is the set of accumulation points of $C$. \nl
If $\lambda \ge \aleph_0$

$$
\align
{\frak d}{\frak p}^{2+}_{\lambda,J} = 
\text{ Min} \biggl\{ \mu:&\text{ there are no }
f_\alpha \in {}^\lambda \lambda \text{ for} \\
  &\,\alpha < \mu 
\text{ such that } \alpha < \beta < \mu \Rightarrow f_\alpha <_J 
f_\beta \biggr\}.
\endalign
$$
\medskip

\noindent
if $\lambda \ge \aleph_0$

$$
\align
{\frak d}{\frak p}^{3+}_{\lambda,J} = 
\text{ Min} \biggl\{ \mu:&\text{ there is no
sequence } \langle A_\alpha:\alpha < \mu \rangle \text{ such that:} \\
  &\,A_\alpha \in J^+ \text{ and} \\
  &\,\alpha < \beta < \mu \Rightarrow [A_\beta
\backslash A_\alpha \in J^+ \and A_\alpha \backslash A_\beta \in J] \biggr\}.
\endalign
$$
\medskip

\noindent
If $J = J^{bd}_\lambda$, we may omit it.  We can replace $J$ by its 
dual filter. \newline
2) For a regular cardinal $\lambda$ let

$$
{\frak d}_\lambda = \text{ Min} \biggl\{|F|:F \subseteq {}^\lambda \lambda
\text{ and } (\forall g \in {}^\lambda \lambda)(\exists f \in F)(g
<_{J^{bd}_\lambda} f) \biggr\}
$$

\noindent
(equivalently $g < f$)

$$
{\frak b}_\lambda = \text{ Min} \biggl\{|F|:F \subseteq {}^\lambda \lambda
\text{ and } \neg(\exists g \in {}^\lambda \lambda)(\forall f \in F)
[f <_{J^{bd}_\lambda} g] \biggr\}.
$$

\noindent
We shall prove here that in the ``neighborhood" of singular cardinals there
are some connections between the ${\frak d}{\frak p}^{\ell+}_\lambda$'s
(hence by monotonicity, also with the ${\frak b}_\lambda$'s).

We first note connections for ``one $\lambda$".
\enddefinition
\bigskip

\demo{\stag{5.2} Fact}  1) If $\lambda = \text{ cf}(\lambda) > \aleph_0$ then

$$
{\frak b}_\lambda < {\frak d}{\frak p}^{1+}_\lambda \le {\frak d}{\frak p}
^{2+}_\lambda \le {\frak d}{\frak p}^{3+}_\lambda.
$$

\noindent
2)  ${\frak b}_{\aleph_0} < {\frak d}{\frak p}^{2+}_{\aleph_0} =
{\frak d}{\frak p}^{3+}_{\aleph_0}$. \nl
3) In the definition of ${\frak d}{\frak p}^{1+}_\lambda,C_{\alpha +1}
\subseteq \text{ acc}(C_\alpha) \text{ mod } J^{\text{bd}}_\lambda$ suffice.
\enddemo
\bigskip

\demo{Proof}  1) \underbar{First inequality}:  
${\frak b}_\lambda < {\frak d}{\frak p}^{1+}_\lambda$.

We choose by induction on $\alpha < {\frak b}_\lambda$, a club $C_\alpha$
of $\lambda$ such that \newline
$\beta < \alpha \Rightarrow |C_\beta \backslash
C_\alpha| < \lambda$ and $\beta < \alpha \Rightarrow C_{\beta + 1} \subseteq
\text{ acc}(C_\beta)$.

For $\alpha = 0$ let $C_\alpha = \lambda$, for $\alpha = \beta + 1$ let
$C_\alpha = \text{ acc}(C_\beta)$, and for $\alpha$ limit let, for each
$\beta < \alpha,f_\beta \in {}^\lambda \lambda$ be defined by $f_\beta(i) = 
\text{ Min}(C_\alpha \backslash (i+1))$.  So $\{f_\beta:\beta < \alpha\}$ 
is a subset of ${}^\lambda \lambda$ of cardinality $\le |\alpha| < {\frak b}
_\lambda$, so there is $g_\alpha \in {}^\lambda \lambda$ such that
$\beta < \alpha \Rightarrow f_\beta <_{J^{bd}_\lambda} g_\alpha$.  

Lastly,
let $C_\alpha = \{ \delta < \lambda:\delta \text{ a limit ordinal such that }
(\forall \zeta < \delta)[g_\alpha(\zeta) < \delta]\}$, now $C_\alpha$ is as
required.

So $\langle C_\alpha:\alpha < {\frak b}_\lambda \rangle$ exemplifies
${\frak b}_\lambda < {\frak d}{\frak p}^{1+}_\lambda$.
\bigskip

\noindent
\underbar{Second inequality}:
${\frak d}{\frak p}^{1+}_\lambda \le {\frak d}{\frak p}^{2+}_\lambda$

Assume $\mu < {\frak d}{\frak p}^{1+}_\lambda$.  Let $\langle C_\alpha:
\alpha < \mu \rangle$ exemplify it, and let us define for $\alpha < \mu$
the function $f_\alpha \in {}^\lambda \lambda$ by: $f_\alpha(\zeta)$ is the
$(\zeta + 1)$-th member of $C_\alpha$; clearly $f_\alpha \in
{}^\lambda \lambda$ and $f_\alpha$ is strictly increasing.  Also, if $\beta <
\alpha$ then $C_\beta \backslash C_\alpha$ is a bounded subset of $\lambda$,
say by $\delta_1$, and there is $\delta_2 \in
(\delta_1,\lambda)$ such that otp$(\delta_2 \cap C_\beta) = \delta_2$.  So
for every $\zeta \in [\delta_2,\lambda)$ clearly $f_\beta(\zeta) =$ the
$(\zeta + 1)$-th member of $C_\beta =$ the $(\zeta +1)$-th member of
$C_\beta \backslash \delta_1 \le$ the $(\zeta + 1)$-th member of $C_\alpha$.
So $\beta < \alpha \Rightarrow f_\beta \le_{J^{bd}_\lambda} f_\alpha$.  
Lastly,
for $\alpha < \mu,C_{\alpha + 1} \subseteq \text{ acc}(C_\alpha)$
hence $f_\alpha(\zeta) =$ the $(\zeta + 1)$-th member of $C_\alpha <$ 
the $(\zeta + \omega)$-th member of $C_\alpha \le$ the $(\zeta + 1)$-th
member of acc$(C_\alpha) \le$ the $(\zeta + 1)$-th member of 
$C_{\alpha +1}$.  So $\beta < \alpha \Rightarrow f_\beta <_{J^{bd}_\lambda}
f_{\beta + 1} \le_{J^{bd}_\lambda} f_\alpha$, so $\langle f_\alpha:\alpha <
\lambda \rangle$ exemplifies $\mu < {\frak d}{\frak p}^{2+}_\lambda$.
\bigskip

\noindent
\underbar{Third inequality}: 
${\frak d}{\frak p}^{2+}_\lambda \le {\frak d}{\frak p}^{3+}_\lambda$

Assume $\mu < {\frak d}{\frak p}^{2+}_\lambda$ and let $\langle f_\alpha:
\alpha < \mu \rangle$ exemplify this.

Let $c:\lambda \times \lambda \rightarrow \lambda$ be one to one and let

$$
A_\alpha = \{c(\zeta,\xi):\zeta < \lambda \text{ and } \xi < f_\alpha
(\zeta)\}.a
$$

\noindent
Now $\langle A_\alpha:\alpha < \mu \rangle$ exemplifies 
$\mu < {\frak d}{\frak p}^{3+}_\lambda$. \newline
\medskip

\noindent
2), 3)  Easy. \hfill$\square_{\scite{5.2}}$
\enddemo
\bigskip

\proclaim{\stag{5.3} Observation}  Suppose 
$\lambda = \text{ cf}(\lambda) > \aleph_0$. \newline
1) If $\langle f_\alpha:\alpha \le \gamma^* \rangle$ is $<_{J^{bd}
_\lambda}$-increasing \underbar{then} we can find a sequence 
$\langle C_\alpha:\alpha < \gamma^* \rangle$ of clubs of $\lambda$, such 
that $\alpha < \beta \Rightarrow
|C_\alpha \backslash C_\beta| < \lambda$ and $C_{\alpha + 1} \subseteq
\text{ acc}(C_\alpha) \text{ mod } J^{bd}_\lambda$. \newline
2)  ${\frak d}{\frak p}^{1+}_\lambda = {\frak d}{\frak p}^{+2}_\lambda$ or
for some $\mu,{\frak d}{\frak p}^{1+}_\lambda = 
\mu^+,{\frak d}{\frak p}^{2+}_\lambda = \mu^{++}$ (moreover though there 
is in $({}^\lambda \lambda,<_{J^{bd}_\lambda}$) an increasing sequence of 
length $\mu^+$, there is none of length $\mu^+ + 1$).
\endproclaim
\bigskip

\demo{Proof}  1) Let

$$
\align
C^* = \biggl\{ \delta < \lambda:&\,\delta \text{ a limit ordinal and }
(\forall \beta < \delta)f_{\gamma^*}(\beta) < \delta \\
  &\text{ and } \omega^\delta = \delta 
\text{ (ordinal exponentiation)} \biggr\};
\endalign
$$

\noindent
this is a club of $\lambda$.

For each $\alpha < \gamma^*$ let

$$
\align
C_\alpha = \biggl\{ \delta + \omega^{f_\alpha(\delta)} \cdot \beta:&\,\delta
\in C^* \text{ and } \beta < f_\alpha(\delta) \\
  &\text{ and } f_\alpha(\delta) < f_{\gamma^*}(\delta) \biggr\}.
\endalign
$$

\noindent
2) Follows. \hfill$\square_{\scite{5.3}}$
\enddemo
\bigskip

\noindent
Now we come to our main concern.
\proclaim{\stag{5.4} Theorem}  Assume
\medskip
\roster
\item "{$(a)$}"  $\kappa$ is regular uncountable, $\ell \in \{1,2,3\}$
\sn
\item "{$(b)$}"  $\langle \mu_i:i < \kappa \rangle$ is (strictly) increasing
continuous with limit $\mu$, \newline
$\lambda_i = \mu^+_i,\lambda = \mu^+$
\sn
\item "{$(c)$}"  $2^\kappa < \mu$ and $\mu^\kappa_i < \mu$
\sn
\item "{$(d)$}"  $D$ a normal filter on $\kappa$
\sn
\item "{$(e)$}"  $\theta_i < {\frak d}{\frak p}^{\ell+}_{\lambda_i}$ and
$\theta = \text{ tcf}(\dsize \prod_{i < \kappa} \theta_i/D)$ or just \newline
$\theta <$ Depth$^+(\dsize \prod_{i < \kappa} \theta_i/D)$.
\endroster
\medskip

\noindent
\underbar{Then} $\theta < {\frak d}{\frak p}^{\ell +}_\lambda$.
\endproclaim
\bigskip

\demo{Proof}  By \scite{5.11}, \scite{5.12}, \scite{5.5} below for 
$\ell = 1,2,3$ respectively (the conditions there are easily checked).
\hfill$\square_{\scite{5.4}}$
\enddemo
\bigskip

\remark{\stag{5.4A} Remark}  1) Concerning assumption (e), 
e.g. if $2^{\mu_i} = \mu^{+5}_i$ and $2^\mu = \mu^{+5}$, then necessarily 
$\mu^{+ \ell} =
\text{ tcf}(\dsize \prod_{i < \kappa} \mu^{+ \ell}_i/D)$ for $\ell =
1,\dotsc,5$ and so $\dsize \bigwedge_{i < \kappa}{\frak d}{\frak p}
^{\ell +}_{\lambda_i} = 2^{\mu_i} \Rightarrow {\frak d}{\frak p}_\lambda 
= 2^\mu$ and we can use $\mu_i = (2^\kappa)^{+i},\lambda_i = \mu^+_i,
\theta_i = \mu^{+5}_i,\theta = \mu^{+5}$.

So this theorem really says that the function $\lambda \mapsto {\frak d}
{\frak p}_\lambda$ has more than the cardinality exponentiation restrictions.
\newline
2) Note that Theorem \scite{5.4} is trivial if $\dsize \prod_{i < \kappa}
\lambda_i = 2^\mu = \lambda$, so (see \cite[V]{Sh:g}) it is natural to assume
$E =: \{D':D' \text{ a normal filter on } \kappa\}$ is nice, but this will not
be used. \newline
3) Note that the proof of \scite{5.12} (i.e. the case $\ell = 2$) does not
depend on the longer proof of \scite{5.5}, whereas the proof of \scite{5.11}
does. \nl
4) Recall that for an $\aleph_1$-complete filter $D$ say on $\kappa$ and 
$f \in {}^\kappa\text{Ord}$ we define $\|f\|_D$ by $\|f\|_D = \cup \{\|g\|_D
+1:g \in {}^\kappa\text{Ord}$ and $g <_D f\}$. \nl
5) Below we shall use the assumption
\medskip
\roster
\item "{$(*)$}"  $\|\lambda\|_{D+A} = \lambda$ for every $A \in D^+$. \nl
This is not strong assumption as
{\roster
\itemitem{ $(a)$ }  if SCH holds, then the only case of interests if $\langle
\chi_i:i < \kappa \rangle$ is increasing continuous with limit $\chi$ and
$\| \langle \chi^+_i:i < \kappa \rangle\|_D = \chi^+$ for any normal filter 
$D$ on $\kappa$; so our statements are degenerated and say nothing,
\sn
\itemitem{ $(b)$ }  if SCH fails, there are nice filters for which this
phenomenum is ``popular" see \cite[V,1.13,3.10]{Sh:g} (see more in 
\sciteu{10.17}).
\endroster}
\endroster
\endremark
\bigskip

\proclaim{\stag{5.5} Theorem}  Assume
\medskip
\roster
\item "{$(a)$}"  $D$ is an $\aleph_1$-complete filter on $\kappa$ 
\sn
\item "{$(b)$}"  $\langle \lambda_i:i < \kappa \rangle$ is a sequence of
regular cardinals $> (2^\kappa)^+$
\sn
\item "{$(c)$}"  $\| \langle \lambda_i:i < \kappa \rangle\|_{D+A} = \lambda$
for $A \in D^+$
\sn
\item "{$(d)$}"  $\mu_i < {\frak d}{\frak p}^{3+}_{\lambda_i}$
\sn
\item "{$(e)$}"  $\mu = \text{ tcf}(\Pi \mu_i/D)$ or at least
\sn
\item "{$(e^-)$}"  $\mu < \text{ Depth}^+
(\Pi \mu_i,<_D)$ and $\mu > 2^\kappa$.
\endroster
\medskip

\noindent
\underbar{Then} $\mu < {\frak d}{\frak p}^{3+}_\lambda$.
\endproclaim
\bigskip

\remark{Remark}  Why not assume just $\|f\|_D = \lambda$ for $f =:
\langle \lambda_i:i < \kappa \rangle$?  Note that cla$^\alpha_I(f,A)$, 
see below, does not make much sense.
\endremark
\bn
A related theorem
\definition{\stag{5.10A} Definition} 

$$
\align
{\frak a}_\lambda = \text{ Min} \biggl\{\mu:&\text{there is no }
{\Cal P} \subseteq [\lambda]^\lambda \text{ of cardinality} \\
  &\mu \text{ such that } A \ne B \in {\Cal P} \Rightarrow |A \cap B| <
\lambda \biggr\}.
\endalign
$$ 
\enddefinition
\bigskip

\proclaim{\stag{5.10B} Theorem}  Assume
\medskip
\roster
\item "{$(a)$}"  $D$ is an $\aleph_1$-complete filter on $\kappa$ 
\sn
\item "{$(b)$}"  $\langle \lambda_i:i < \kappa \rangle$ is a sequence of
regular cardinals $> (2^\kappa)^+$
\sn
\item "{$(c)$}"  $\| \langle \lambda_i:i < \kappa \rangle\|_{D+A} = \lambda$
for $A \in D^+$
\sn
\item "{$(d)$}"  $\mu_i < {\frak a}_{\lambda_i}$
\sn
\item "{$(e)$}"  $\mu = \text{ tcf}(\Pi \mu_i/D)$ or at least
\sn
\item "{$(e^-)$}"  $\mu < \text{ Depth}^+
(\Pi \mu_i,<_D)$ and $\mu > 2^\kappa$.
\ermn
\ub{Then} $\mu < {\frak a}_\lambda$. \nl
We shall prove later.
\endproclaim
\bigskip

\noindent
We delay the proof of \scite{5.5}.
\demo{\stag{5.6} Fact}  Assuming \scite{5.5}(a), for any 
$f \in {}^\kappa(\text{Ord} \backslash (2^\kappa)^+)$ we have: 
$T_D(f)$ is smaller or equal to the cardinality of $\|f\|_D$ remembering
(\sciteu{10.5}(4) above and)

$$
T_D(f) = \text{ sup} \biggl\{ |F|:F \subseteq \dsize \prod_{i < \kappa} f(i) 
\text{ and } f \ne g \in F \Rightarrow f \ne_D g \biggr\}.
$$
\enddemo
\bigskip

\demo{Proof}  Why?  Let $F$ be as in the definition of $T_D(f)$, note:
$f_i \ne_D f_j \and f_i \le_D f_j \Rightarrow f_i <_D f_j$.  Note that
$i < \kappa \Rightarrow f(i) \ge (2^\kappa)^+$, necessarily $|F| > 2^\kappa$.
Now for each ordinal $\alpha$ let $F^{[\alpha]} =: \{f \in F:\|f\|_D = 
\alpha\}$.  Clearly $F^{[\alpha]}$ has at most
$2^\kappa$ members (otherwise $f_i \in F^{[\alpha]}$ for
$i < (2^\kappa)^+$ are pairwise distinct so for some $i < j,f_i <_D f_j$
\cite[\S2]{Sh:111} or simply use Erd\"os-Rado on 
$c(i,j) = \text{ min}\{ \zeta < \kappa:f_i(\zeta) > f_j(\zeta)\}$). \newline
So $\|f\|_D \ge \text{ sup}\{ \|g\|_D:
g \in F\} \ge \text{ otp}\{\alpha:F^{[\alpha]} \ne \emptyset\} \ge |\{
F^{[\alpha]}:F^{[\alpha]} \ne \emptyset\}| \ge |F|/2^\kappa = |F|$.  So
$\|f\|_D \ge T_D(f)$.
\hfill$\square_{\scite{5.6}}$
\enddemo
\bigskip

\definition{\stag{5.7} Definition}  For $f \in {}^\kappa\text{Ord}$ (natural
to add $0 \notin \text{ Rang}(f))$ and $D$ an
$\aleph_1$-complete filter on $\kappa$ let
$\dsize \prod^*_{i < \kappa} f(i) = \{g:\text{Dom}(g) = \kappa,f(i) > 0
\Rightarrow g(i) < f(i)$ and $f(i) = 0 \Rightarrow g(i) = 0\}$ and
\sn
1)  cla$(f,D) = \biggl\{(g,A):g \in 
\dsize \prod^*_{i < \kappa}f(i) \text{ and }
A \in D^+ \biggr\}$ \newline
\smallskip
\noindent
$\quad$ cla$^\alpha(f,D) = \{(g,A) \in \text{ cla}(f,D):\|g\|_{D+A} =
\alpha\}$. 
\sn
Here ``cla" abbreviates ``class".
\smallskip
\noindent
2)  For $(g,A) \in \text{ cla}(f,D)$ let

$$
J_D(g,A) = \{ B \subseteq \kappa:\text{ if } B \in (D+A)^+ \text{ then }
\|g\|_{(D+A)+B} > \|g\|_{D+A}\}.
$$

\smallskip
\noindent
3) We say $(g',A') \approx (g'',A'')$ if (both are in cla$(f,D)$ and)
$A' = A'' \text{ mod } D$ and $J_D(g',A') = J_D(g'',A'')$ and $g' = g''
\text{ mod } J_D(g',A')$. \newline
4)  For $I$ an ideal on $\kappa$ disjoint to $D$ we let

$$
I*D = \{A \subseteq \kappa:\text{for some } X \in D \text{ we have }
A \cap X \in I\},
$$
\medskip

\noindent
(usually we have $\{\kappa \backslash A:A \in D\} \subseteq I$ so
$I*D = I$), and let

$$
\text{cla}_I(f,D) = \{ (g,A):g \in \dsize \prod^*_{i < \kappa}f(i)
\text{ and } A \in (I*D)^+ \}.
$$
\medskip

\noindent
5) On cla$_I(f,D)$ we define a relation $\approx_I$ \newline
$(g_1,A_1) \approx_I (g_2,A_2)$ if:
\medskip
\roster
\item "{$(a)$}"  $A_1 = A_2$ mod $D$ and
\sn
\item "{$(b)$}"  there is $B_0 \in I$ such that: if $B_0 \subseteq B \in I$
then \newline
$\|g_1\|_{(D+A_1)+(\kappa \backslash B)} = \|g_2\|_{(D+A_2)+(\kappa \backslash
B)}$ and \newline
$J_{(D+A_1)+(\kappa \backslash B)}(g_1,A_1) = J_{(D+A_1)+(\kappa \backslash
B)}(g_2,A_2)$.
\endroster

$$
\align
J_{D,I}(g_1,A_1) = \{A \subseteq \kappa:&\text{for some } B_0 \in I
\text{ if } B_0 \subseteq B \in I \tag "{$6)$}" \\
  &\text{we have } A \in J_{(D+A_1)+(\kappa \backslash B_0)}(g_1,A_1) \}.
\endalign
$$
\medskip

\noindent
7) Let com$(D)$ be \ub{the maximal} $\theta$ such that $D$ is 
$\theta$-complete.
\enddefinition
\bigskip

\demo{\stag{5.8} Fact}  For $f \in {}^\kappa\text{Ord}$ and $D$ an 
$\aleph_1$-complete filter on $\kappa$ and $A \in D^+$: \newline
0) If $f_1 \le f_2$ then cla$(f_1,D) \subseteq \text{ cla}(f_2,D)$ and for
$g',g'' \in \dsize \prod^*_i f_1(i),A \subseteq \kappa$ we have
$(g',A) \approx (g'',A)$ in cla$(f_1,D)$ iff $(g',A) \approx (g'',A)$ in
cla$(f_2,D)$ (so we shall be careless about this). \newline
1) $J_D(g,A)$ is an ideal on $\kappa,\text{com}(D)$-complete, and normal if
$D$ is normal.  \newline
2) $A$ does not belong to $J_D(g,A)$, and it includes $\{B \subseteq \kappa:
B = \emptyset \text{ mod } (D+A)\}$.  If $B \in J^+_D(g,A)$ then $A \cap B
\in D^+$ and $\|g\|_{D+(A \cap B)} = \|g\|_{D+A}$. \newline
3) $\approx$ is an equivalence relation on cla$(f,D)$, similarly $\approx_I$
on cla$_I(f,D)$. \newline
4) Assume
\medskip
\roster
\item "{$(i)$}"   $(g,A) \in \text{ cla}^\alpha(f,D),g' \in 
\dsize \prod_{i < \kappa} f(i)$ and
\sn
\item "{$(ii)$}"  (a) $\quad g' = g \text{ mod}(D+A)$ or \newline
                  (b) $\,\,$ for some $B \in J_D(g,A)$ we have 
                             $\alpha \in B \Rightarrow g'(\alpha)
                             > \|g\|_D$ \nl
		
                	     $\qquad$ (or just 
                             $\|g\|_{D+A} < \|g'\|_{D+B}$) and \nl
 
$\qquad \qquad g' \restriction (\kappa \backslash B) = g \restriction (\kappa
\backslash B) \text{ mod } D$.
\endroster
\medskip

\noindent
\underbar{Then} $(g',A) \approx (g,A)$. \newline
5) For each $\alpha$, in cla$^\alpha(f,D)/\approx$ there are at 
most $2^\kappa$ classes. \nl
6) For $f \in {}^\kappa(\text{Ord})$, in cla$(f,D)/ \approx$
there are at most $2^\kappa + \sup_{A \in D^+} \|f\|_{D+A}$ classes.

\enddemo
\bigskip

\demo{Proof}  0) Easy. \nl
1) Straight (e.g. ideal as for $B \subseteq \kappa$ we have
\newline
$\|g\|_D = \text{ Min}\{\|g\|_{D+A},\|g\|_{D+(\kappa-A)}\}$, where we
stipulate $\|g\|_{{\Cal P}(\kappa)} = \infty$ see \cite{Sh:71}). \newline
2) Check. \newline
3) Check. \newline
4) Check. \newline
5) We can work also in cla$^\alpha(f+1,D)$ (this change gives more
elements and by (0) it preserves $\approx$).  Assume $\alpha$ is a 
counterexample (note that ``$\le 2^{2^\kappa}$" is totally immediate).  
Let $\chi$ be large enough; choose $N \prec ({\Cal H}(\chi),
\in,<^*_\chi)$ of cardinality $2^\kappa$ such that $\{f,D,\kappa,\alpha\} \in
N$ and ${}^\kappa N \subseteq N$.  So necessarily there is $(g,A) \in
\text{ cla}^\alpha(f,D)$ such that $(g,A)/ \approx \notin N$, by the 
definition of cla$^\alpha$, clearly $\|g\|_D + A = \alpha$.  
Let $B =: \{i < \kappa:g(i) \notin N\}$.
\enddemo
\bigskip

\noindent
\underbar{Case 1}:  $B \in J_D(g,A)$. 

Let $g' \in \dsize \prod_{i < \kappa} (f(i)+1)$ be defined by: 
$g'(i) = g(i)$ if
$i \in \kappa \backslash B$ and $g'(i) = f(i)$ if $i \in B$.  By part (4) we
have $(g',A) \approx (g,A)$ and by the choice of $N$ we have $(g',A) \in N$
as $A \in {\Cal P}(\kappa) \subseteq N,g' \in N$ (as Rang$(g') \subseteq N 
\and {}^\kappa N \subseteq N)$.  So as $({\Cal H}(\chi),\in,<^*_\chi) \models
(\exists x \in \text{ cla}^\alpha(f,A))(x \approx (g',A))$ there is $(g'',A) 
\in N$ such that $(g'',A) \approx (g,A)$ as required. 
\bigskip

\noindent
\underbar{Case 2}:  $B \notin J_D(g,A)$.

Let $g' \in {}^\kappa\text{Ord}$ be: $g'(i) = \text{ Min}(N \cap (f(i)+1)
\backslash g(i)) \le f(i)$ if $i \in B,g'(i) = g(i)$ if $i \notin B$
(note: $f(i) \in N,g(i) < f(i)$ so $g'$ is well defined).

Clearly $g' \in N$, (as Rang $(g') \subseteq N$ and 
${}^\kappa N \subseteq N$), and

$$
\align
({\Cal H}(\chi),\in,<^*_\chi) \models &(\exists x)(x \in \dsize 
\prod^*_{i < \kappa}
f(i) \wedge (\forall i \in \kappa \backslash B)(x(i) = g'(i)) \wedge \\
  &(\forall i \in B)(x(i) < g'(i)) \wedge \|x\|
_{D+(A \cap B)} = \alpha)
\endalign
$$

\noindent
(why?  because $x=g$ is like that, last equality as $B \notin J_D(g,A)$).
So there is such $x$ in $N$, call it $g''$.  So $g'' \in \dsize \prod_{i <
\kappa}(f(i)+1)$ and $\|g''\|_{D+(A \cap B)} = \alpha$ and for \newline
$i \in B,g''(i)
\in g'(i) \cap N$ hence $g''(i) < g(i)$ by the definition of $g'(i)$.
\newline
So $g'' < g \text{ mod } D+(A \cap B)$, but this contradicts
$\|g''\|_{D+(A \cap B)} = \alpha = \|g\|_{D+(A \cap B)}$, the last equality
as $B \notin J_D(g,A)$. \nl
6) Immediate from (5).   \hfill$\square_{\scite{5.8}}$
\bigskip

\demo{\stag{5.9} Fact}  Assume $f \in {}^\kappa\text{Ord}$ and $D$ an 
$\aleph_1$-complete filter on $\kappa$ and $I$ an com$(D)$-complete ideal 
on $\kappa$. \newline
1) If $(g,A) \in \text{ cla}_I(f,D)$ then $J_{D,I}(g,A)$ is an
ideal on $\kappa$, which is com$(D)$-complete and normal if $D,I$
are normal. \newline
If $B \in (J_{D,I}(g,A))^+$ then $\|g\|_{D+(A \cap B)} = \|g\|_{D+A}$, and
$(D + (A \cap B)) \cap I = \emptyset$. \newline
2) $\approx_I$ is an equivalence relation on cla$(f,D)$. \newline
3) If $(g,A) \in \text{ cla}(f,D)$ and $g' \in \dsize \prod^*_{i < \kappa}
f(i)$ and $g' = g \text{ mod } J_{D,I}(g,A)$ \underbar{then} for some
$A'$ we have $(g',A') \approx_I (g,A')$ so $(g',A') \in \text{ cla}(f,D)$ and
$\|g'\|_{D+A'} = \|g\|_{D+A'}$ (in fact $A' = \{i \in A:g'(i) = g(i)\}$ is
O.K.).
\enddemo
\bigskip

\demo{Proof}  Easy.
\enddemo
\bigskip

\demo{\stag{5.9A} Fact}  Let $\kappa,f,D$ be as in \scite{5.9}. \nl
1) If $f_\zeta \in {}^\kappa\text{Ord}$, for $\zeta \le \delta$, 
cf$(\delta) > \kappa$ and for each $i$ the sequence $\langle 
f_\zeta(i):\zeta \le \delta \rangle$ is increasing continuous {then} 
$\|f_\delta\|_D = \underset {\zeta < \delta}\to \sup \|f_\zeta\|_D$. \newline
2) If $\delta = \|f\|_D$, cf$(\delta) > 2^\kappa$ \ub{then}
$\{i:\text{cf}(f(i)) \le 2^\kappa\} \in J_D(f,\kappa)$. \newline
3) If $\|f\|_D = \delta,A \in J^+_D(f,\kappa)$ \ub{then} 
$\dsize \prod^*_{i<\kappa} f(i)/(D+A)$ is not (cf$(\delta))^+$-directed. \nl
4) If $\|f\|_D = \delta$ and $A \in J^+_D(f,\kappa)$ \underbar{then}
cf$(\delta) \le \text{ cf}(\dsize \prod^*_{i<\kappa}f(i)/(D+A))$. \newline
5) If $\|f\|_D = \delta$ and $A \subseteq \kappa,(\forall i \in A)
\text{cf}(f(i)) > \kappa$ and \newline
max pcf$\{f(i):i \in A\} < \text{ cf}(\delta)$ \newline
(or just cf$(\delta) > \text{ max}\{\text{cf } \dsize \prod^*_{i <\kappa}
f(i)/D':D' \text{ an ultrafilter extending } D+A\}$) \ub{then} 
$A \in J_D(f,\kappa)$. \newline
6) If $\|f\|_D = \delta$, cf$(\delta) > 2^\kappa$, \ub{then} 
$\dsize \prod^*_{i <\kappa}f(i)/J_D(f,\kappa)$
is cf$(\delta)$-directed. \newline
7) If $\|f\|_D = \delta$, cf$(\delta) > 2^\kappa$, \ub{then} for some
$A \in J^+_D(f,\kappa)$ we have \nl
$\dsize \prod^*_{i<\kappa}f(i)/J_D(f,\kappa)+(kappa \backslash A)$ 
has true cofinality cf$(\delta)$. \newline
8) Assume $\|f\|_D = \lambda = \text{ cf}(\lambda) > 2^\kappa$.  \newline
Then $(\forall A \in D^+)(\|f\|_{D+A} = \lambda)$ \ub{implies} 
tcf$(\dsize \prod^*_{i < \kappa} f(i)/D) = \lambda$. \nl
9) If $\|f\|_D = \delta$, cf$(\delta) > 2^\kappa$ \ub{then} 
tcf $\dsize \prod^*_{i < \kappa} f(i)/J_D(f,\kappa) = \text{ cf}(\delta)$. 
\enddemo
\bigskip

\demo{Proof}  1) Let $g <_D f_\delta$, so $A = \{i < \kappa:g(i) <
f_\delta(i)\} \in D$, now for each $i \in A$ we have $g(i) < f_\delta(i)
\Rightarrow (\exists \alpha < \delta)(g(i) < f_\alpha(i)) \Rightarrow$ there
is $\alpha_i < \delta$ such that $(\forall \alpha)[\alpha_i \le \alpha \le
\delta \Rightarrow g(i) < f_{\alpha_i}(i)]$.  Hence $\alpha(*) =: \sup \{
\alpha_i:i \in A\} < \delta$ as cf$(\delta) > \kappa$, so 
$g <_D f_{\alpha(*)}$ hence $\|g\|_D < \|f_{\alpha(*)}\|_D$; this suffices
for one inequality, the other is trivial. \newline
2) Let $A = \{i:\text{cf}(i) \le 2^\kappa\}$, and assume toward contradiction
that $A \in J^+_D(f,\kappa)$.  For each $i \in A$ let $C_i \subseteq f(i)$ be
unbounded of order type cf$(f(i)) \le 2^\kappa$.  \newline
Let $F = \{g \in \dsize \prod^*_{i < \kappa}(f(i)+1)$: if $i \in A$ then 
$g(i) \in C_i$, if $i \in \kappa \backslash A$ then $g(i) = f(i)\}$.  
So $|F| \le 2^\kappa$ and:
\medskip
\roster
\item "{$(*)$}"  if $g <_{D+A} f$ then for some $g' \in F,g <_{D+A} g'$ \nl
hence $\delta = \|f\|_{D+A} = \sup\{\|g\|_{D+A}:g \in F\}$ but the supremum
is on $\le |F| < \text{ cf}(\delta)$ ordinals each $< \delta$ because
$g' \in F \Rightarrow g' <_{D+A} f$ as $\|f\|_D = \delta \Rightarrow f \ne_D
0_\kappa$, contradiction to cf$(\delta) > 2^\kappa$.
\endroster
\medskip

\noindent
3) Assume this fails, so $\|f\|_D = \delta,A \in J^+_D(f,\kappa)$ anad
$\dsize \prod^*_{i < \kappa} f(i)/(D+A)$ is (cf$(\delta))^+$-directed.  Let
$C \subseteq \delta$ be unbounded of order type cf$(\delta)$; as
$\|f\|_{D+A} = \delta$ (because $A \in J^+_D(f,A))$ for each $\alpha \in C$
there is $f_\alpha <_{D+A} f$ such that $\|f_\alpha\|_{D+A} \ge \alpha$ (even
$= \alpha$ by the definition of $\|-\|_{D+A}$).  As $\dsize \prod^*_{i<\kappa}
f(i)/(D+A)$ is (cf$(\delta))^+$-directed there is $f' <_{D+A} f$ such that
$\alpha \in C \Rightarrow f_\alpha <_{D+A} f'$.  By the first inequality
$\|f'_{D+A}\| < \|f\|_{D+A} = \delta$, and by the second inequality
$\alpha \in C \Rightarrow \alpha \le \|f_\alpha\|_{D+A} \le \|f'\|_{D+A}$
hence $\delta = \sup(C) \le \|f'\|_{D+A}$, a contradiction. \newline
4) Same proof as part (2). \newline
5) By part (4) and \cite[Ch.II,3.1]{Sh:g}. \newline
6) Follows. \newline
7) If not, by part (2) without loss of generality $\forall i[\text{cf}(f(i)) 
> 2^\kappa]$; let $C \subseteq \delta$ be unbounded, otp$(C) = \text{ cf}
(\delta)$.  For
each $\alpha \in C$ and $A \in J^+_D(f,\kappa)$ choose $f_{\alpha,A} <_D f$
such that $\|f_{\alpha,A}\|_{D+A} = \alpha$.  Let $f_\alpha$ be 
$f_\alpha(i) = \sup\{f_{\alpha,A}(i):A \in J^+_D(f,\kappa)\}$.
As $\bigl( \dsize \prod^*_{i < \kappa} f_\alpha(i),<_{J_D(f,\kappa)} \bigr)$ 
is cf$(\delta)$-directed (see part (8)), by the assumption 
toward contradiction and the pcf theorem we have $\dsize \prod^*_{i < \kappa}
f(i)/J_D(f,\kappa)$ is
(cf$(\delta))^+$-directed.  Hence we can find $f^* < f$ such that
$\alpha \in C \Rightarrow f_\alpha <_{J_D(f,\kappa)} f^*$.  Let 
$\beta = \sup\{\|f^*\|_{D+B}:B \in J^+_D(f,A)\}$, it
is $< \delta$ as cf$(\delta) > 2^\kappa$; hence there is $\alpha,
\beta < \alpha \in C$, so by the choice of $f^*$ we have 
$f_\alpha <_{J_D(f,\kappa)}f^*$, and let $A =: \{i < \kappa:
f_\alpha(i) < f^*(i)\}$ so $A \in J^+_D(f,\kappa)$, so $f_{\alpha,A} \le
f_\alpha <_{D+A} f^*$ hence $\alpha \le \|f_{\alpha,A}\|_{D+A} \le
\|f_\alpha\|_{D+A} \le \|f^*\|_{D+A} \le \beta$ contradicting the choice of
$\alpha$. \nl
8) For every $\alpha < \lambda$ we can choose $f_\alpha <_D f$ such that
$\|f_\alpha\|_D = \alpha$.  Let $a_\alpha = \{\|f_\alpha\|_{D+A}:A \in D^+\}$,
as $A \in D^+ \Rightarrow \alpha \le \|f_\alpha\|_D \le
\|f_\alpha\|_{D+A} < \|f\|_{D+A} = \lambda$,
clearly $a_\alpha$ is a subset of $\lambda \backslash \alpha$, 
and its cardinality is $\le 2^\kappa < \lambda$.  So we can find an 
unbounded $E \subseteq \lambda$
such that $\alpha < \beta \in E \Rightarrow \sup(a_\alpha) < \beta$.  So
if $\alpha < \beta,\alpha \in E,\beta \in E$, let $A = \{i < \kappa:f_\alpha
(i) \ge f_\beta(i)\}$, and if $A \in D^+$, \ub{then} $\|f_\beta\|_{D+A} \le
\|f_\alpha\|_{D+A} \le \sup a_\alpha < \beta$, contradiction.  Hence
$A = \emptyset \text{ mod } D$ that is $f_\alpha <_D f_\beta$.  Also if
$g <_D f$, then $a =: \{\|g\|_{D+A}:A \in D^+\}$ is again a subset of 
$\lambda$ of cardinality $\le 2^\kappa$ hence for some $\beta < \lambda$,
sup$(a) < \beta$, so as above $g <_D f_\beta$.  Together $\langle f_\alpha:
\alpha \in E \rangle$ exemplify $\lambda = \text{ tcf}(\Pi f(i),<_D)$. \nl
9) Similar proof (to part (8)), using parts (6), (7).
\hfill$\square_{\scite{5.9A}}$
\enddemo
\bigskip

\remark{\stag{5.10} Remark}  We think Claims \scite{5.8}, \scite{5.9},
\scite{5.9A} (and Definition \scite{5.7}) can be applied to the problems 
from \cite{Sh:497}.
\endremark
\bigskip

\demo{Proof of \scite{5.5}}   Fix $f \in {}^\kappa\text{Ord}$ as
$f(i) = \lambda_i$ and let $\approx,\approx_I$ be as in Definition
\scite{5.7}.  \nl
For each $i < \kappa$ let $\bar X^i = \langle
X^i_\alpha:\alpha < \mu_i \rangle$ be a sequence of members of $[\lambda_i]
^{\lambda_i}$ such that

$$
\alpha < \beta < \mu_i \Rightarrow X^i_\alpha \backslash X^i_\beta \in
J^{bd}_{\lambda_i} \and X^i_\beta \backslash X^i_\alpha \notin
J^{bd}_{\lambda_i}.
$$
\medskip

\noindent
(it exists by assumption $(d)$). \nl
Let $\bar g^* = \langle g^*_\zeta:\zeta < \mu \rangle$ be a $<_D$-increasing
sequence of members of $\dsize \prod_{i < \kappa} \mu_i$, it exists by
assumption $(e)$ or $(e)^-$.  \newline
Let $I =: \{B \subseteq \kappa:\text{if } B \in D^+ \text{ then }
\|f\|_{D+B} > \lambda\}$, it is a com$(D)$-complete ideal on 
$\kappa$ disjoint to $D$, i.e., $I = J_D(\bar \lambda,\kappa) \supseteq 
\{ \kappa \backslash A:A \in D\}$, and $\approx_I,\approx$ are equal by
assumption $(c)$.  
For any sequence $\bar X = \langle X_i:i < \kappa \rangle \in \dsize
\prod_{i < \kappa}[\lambda_i]^{\lambda_i}$, let 

$$
Y[\bar X] =: \biggl\{ \|h\|_{D+A}:h \in \dsize \prod_{i < \kappa}
X_i \text{ and } A \in I^+ \biggr\}
$$
\mn
and

$$
\align
{\Cal Y}[\bar X] =: \biggl\{ (h,A)/\approx:&\,h \in \dsize \prod_{i < \kappa}
X_i \text{ and } (h,A) \in \text{ cla}^\alpha_I(\bar \lambda,D) \\
  &\text{ for some } \alpha < \lambda \biggr\}
\endalign
$$
\medskip

\noindent
Note: $Y[\bar X] \subseteq \lambda$ and ${\Cal Y}[\bar X] \subseteq
{\Cal Y}^* =: \dsize \bigcup_{\alpha < \lambda} \text{ cla}^\alpha
(\bar \lambda,D)/ \approx$. \newline
Note that by \scite{5.8}(6)
\mr
\item "{$\boxtimes$}"  $\dsize \bigcup_{\alpha < \lambda} 
\text{ cla}^\alpha(f,D)/\approx$ has cardinalty $\le \lambda$.
\sn
\item "{$(*)_0$}"  for $\bar X \in \dsize \prod_{i < \kappa}[\lambda_i]
^{\lambda_i}$, the mapping $(g,A)/\approx_I \, \mapsto \|g\|_{D+A}$ is from
${\Cal Y}[\bar X]$ onto $Y[\bar X]$ with every $\alpha \in Y[\bar X]$ 
having at most $2^\kappa$ preimages \newline
[why? by \scite{5.8}(5)]
\sn
\item "{$(*)_1$}"  if $\bar X \in \dsize \prod_{i < \kappa}[\lambda_i]
^{\lambda_i}$ then ${\Cal Y}[\bar X]$ has cardinality $\lambda$. \newline
[why?  by the definition of $\|-\|_D$ for every $\alpha < \lambda$ for some
$g \in \dsize \prod_{i < \kappa} \lambda_i/D$ we have $\|g\|_D = \alpha$;
as sup$(X_i) = \lambda_i > g(i)$ we can find $g' \in \dsize \prod_{i < \kappa}
(X_i \backslash g(i))$ such that $g \le g' < \langle \lambda_i:
i < \kappa \rangle$, so $\alpha = \|g\|_D \le
\|g'\|_D < \|\langle \lambda_i:i < \kappa \rangle\|_D = \lambda$.
Clearly for some $\alpha'$ and $A,(g',A) \in \text{ cla}^{\alpha'}(f,A)$, 
so $A \in I^+ \subseteq D^+$, and $\alpha \le \alpha' = \|g'\|_{D+A} <
\|f\|_{D+A} = \lambda$ (as $A \in I^+$).  So $\alpha' \in Y[\bar X]$ hence
$Y[\bar X] \nsubseteq \alpha$; as $\alpha < \lambda$ was arbitrary,
$Y[\bar X]$ has cardinality $\ge \lambda$, by $\boxtimes$ equality holds 
hence (by $(*)_0$) also ${\Cal Y}[\bar X]$ has cardinality $\lambda$.]
\sn
\item "{$(*)_2$}"  if $\bar X',\bar X'' \in 
\dsize \prod_{i < \kappa}[\lambda_i]^{\lambda_i}$, and \newline
$\{i < \kappa:X'_i \subseteq X''_i \text{ mod } J^{bd}_{\lambda_i}\}
\in D$ \underbar{then}
{\roster 
\itemitem{ $(a)$ }  $Y[\bar X'] \subseteq Y[\bar X''] \text{ mod } 
J^{bd}_\lambda$
\sn
\itemitem{ $(b)$ }  ${\Cal Y}[\bar X'] \backslash 
{\Cal Y}[\bar X'']$ has cardinality $< \lambda$ \newline
[why? define $g \in \dsize \prod_{i < \kappa} \lambda_i$ by $g(i) =
\text{ sup}(X'_i \backslash X''_i)$ if \newline
$i \in A^* =: \{ i < \kappa:
X'_i \subseteq X''_i$ mod $J^{bd}_{\lambda_i}\}$ and $g(i) = 0$ otherwise.
Let $\alpha(*) = \text{ sup}\{\|g\|_{D+A} + 1:A \in I^+\}$, as $\lambda$ is
regular $> 2^\kappa$ clearly $\alpha(*) < \lambda$ (see assumption (c) or
definition of $I$).
Assume $\beta \in Y[\bar X'] \backslash \alpha(*)$ and we shall prove that
$\beta \in Y[\bar X'']$, moreover, ${\Cal Y}[\bar X'] \cap 
(\text{cla}^\beta (\bar X,D)/\approx_I) \subseteq {\Cal Y}[\bar X'']$, 
this clearly suffices for both clauses.  We can find 
$f^* \in \dsize \prod_{i < \kappa} (X'_i \cap X''_i) \cup \{0\}$) such that
$\|f^*\|_D > \beta$. \newline
So let a member of ${\Cal Y}[\bar X'] \cap (\text{cla}^\beta(\bar \lambda,
D)/\approx)$ have the form $(h,A)/\approx_I$, where $A \in I^+,h \in \dsize
\prod_{i < \kappa} X'_i$ and $\beta = \|h\|_{D+A}$ and let \nl
$A_1 =: \{i < \kappa:h(i) \le g(i)\}$.  We know \nl 
$\beta = \|h\|_{D+A} = \text{ Min}
\{\|h\|_{D+(A \cap A_1)},\|h\|_{D+(A \backslash A_1)}\}$ (if
$A \cap A_1 = \emptyset \text{ mod } D$, then $\|h\|_{D+A \cap A_1}$ can be
considered $\infty$). \newline
If $\beta = \|h\|_{D+(A \cap A_1)}$ then note $h \le_{D+(A \cap A_1)} g$
hence \nl
$\beta = \|h\|_{D+(A \cap A_1)} \le \|g\|_{D+(A \cap A_1)} <
\alpha(*)$, contradiction to an assumption on $\beta$.  So
$\beta = \|h\|_{D+(A \backslash A_1)}$ and $A \cap A_1 \in J_{D,I}(h,A)$, and
define $h' \in \dsize \prod_{i < \kappa} f(i)$ by: $h'(i)$ is $h(i)$ if
$i \in A \backslash A_1$ and $h'(i)$ is $f^*(i)$ if $i \in \kappa \backslash 
(A \backslash A_1)$.  So $h' \in \dsize
\prod_{i < \kappa} f(i)$ and $h' =_{D+(A \backslash A_1)} h$ hence
$\|h'\|_{D+(A \backslash A_1)} = \|h\|_{D+(A \backslash A_1)} = \beta$,
and clearly $\beta = \|h'\|_{D+(A \backslash A_1)} \in Y[\bar X'']$, as 
required for clause (a), moreover $(h,A) \approx (h',A)$ so $(h',A)/\approx
\in {\Cal Y}[\bar X'']$ as required for clause (b).]
\endroster}
\sn
\item "{$(*)_3$}"  If $\bar X',\bar X'' \in \dsize \prod_{i < \lambda}
[\lambda_i]^{\lambda_i}$ and \newline
$\{i < \kappa:X''_i \nsubseteq X'_i \text{ mod } J^{bd}_\lambda\} \in
D$ \underbar{then} \newline
${\Cal Y}[\bar X''] \backslash {\Cal Y}[\bar X']$ has cardinality $\lambda$.
\newline
[why?  let $\alpha < \lambda$, it is enough to find $\beta \in [\alpha,
\lambda]$ such that

$$
({\Cal Y}[\bar X''] \backslash {\Cal Y}[\bar X']) \cap (\text{cla}^\beta
(f,D) / \approx) \ne \emptyset.
$$

\noindent
We can find $g \in \dsize \prod_{i < \kappa} \lambda_i$ such that
$\|g\|_D = \alpha$.  Define $g' \in \dsize \prod_{i < \kappa} X''_i$ by:
$g'(i)$ is Min$(X''_i \backslash X'_i \backslash g(i))$ when well defined,
Min$(X''_i)$ otherwise.  By assumption $g \le_D g'$ and, of course,
$g' \in \dsize \prod_{i < \kappa} X''_i \subseteq \dsize \prod_{i < \kappa}
\lambda_i$, so $\|g'\|_D \ge \alpha$.  So \newline
$\bigl( (g',\kappa)/\approx \bigr) \in {\Cal Y}[\bar X'']$ but trivially 
$((g',\kappa) / \approx) \notin {\Cal Y}[\bar X']$, so we are done].
\endroster
\medskip

\noindent
Together $(*)_0 - (*)_3$ give that $\langle {\Cal Y}[\langle 
X^i_{g^*_\zeta(i)}]:i < \kappa \rangle:\zeta < \mu \rangle$ is a sequence 
of subsets of ${\Cal Y}^*$ of length $\mu$ (see $(*)_1$), $|{\Cal Y}^*| 
= \lambda$,
which is increasing modulo $[{\Cal Y}^*]^{< \lambda}$ (by $(*)_2$), and in
fact, strictly increasing (by $(*)_3$).  So modulo changing names we have
finished. \hfill$\square_{\scite{5.5}}$
\enddemo
\bigskip

\demo{Proof of \scite{5.10B}}  Similar to the proof of \scite{5.5}
\enddemo
\bigskip

\proclaim{\stag{5.11} Theorem}  Assume
\medskip
\roster
\item "{$(a)$}"  $D$ an $\aleph_1$-complete filter on $\kappa$
\sn
\item "{$(b)$}"  $\bar \lambda = \langle \lambda_i:i < \kappa \rangle$ is
a sequence of regular cardinals $> 2^\kappa$
\sn
\item "{$(c)$}"  $\lambda = \|\bar \lambda\|_{D+A}$ for $A \in D^+$
\sn
\item "{$(d)$}"  $\mu_i < {\frak d}{\frak p}^{1+}_{\lambda_i}$
\sn
\item "{$(e)$}"  $\mu < \text{ Depth}^+(\dsize \prod_{i < \kappa} \mu_i,<_D)$.
\endroster
\medskip

\noindent
\underbar{Then} $\mu < {\frak d}{\frak p}^{1+}_\lambda$.
\endproclaim
\bigskip

\demo{Proof}  Let Club$(\lambda) = \{C:C \text{ a club of } \lambda\}$ so
Club$(\lambda) \subseteq [\lambda]^\lambda$ for $\lambda = \text{ cf}(\lambda)
> \aleph_0$.

For any sequence $\bar C \in \dsize \prod_{i < \kappa}$ Club$(\lambda_i)$ let
${\Cal C}(\bar C)$ be the set acc$(c \ell(Y(\bar C))$ where 
$Y[\bar C] =: \{\|g\|_D:g \in \dsize \prod_{i < \kappa} C_i\} 
(\subseteq \lambda)$; i.e. ${\Cal C}(\bar C) = \{ \delta < \lambda:\delta = 
\sup(\delta \cap Y[\bar C])\}$.  Clearly
\medskip
\roster
\item "{$(*)_1$}"  for $\bar C \in \dsize \prod_{i < \kappa}
\text{ Club}(\lambda_i) \text{ we have } {\Cal C}(\bar C) \in
\text{ Club}(\lambda)$ \newline
[the question is why it is unbounded, and this holds as $\| \bar \lambda\|_D =
\lambda$ by its definition]
\sn
\item "{$(*)_2$}"  if $\bar C',\bar C'' \in \dsize \prod_{i < \lambda}$ 
Club$(\lambda_i),g^* \in \Pi \lambda_i$ and $C''_i = C'_i \backslash g^*(i)$
\underbar{then} \newline
${\Cal C}(\bar C') = {\Cal C}(\bar C'')$ mod $J^{bd}_\lambda$ \newline
[why? let $\alpha(*) = \sup \{\|g^*\|_{D+A}:A \in D^+$ and $\|g^*\|_{D+A} <
\lambda\}+1$, so as $2^\kappa < \lambda = \text{ cf}(\lambda)$ clearly
$\alpha(*) < \lambda$.  We shall show ${\Cal C}(\bar C') \backslash
\alpha(*) = {\Cal C}(\bar C'') \backslash \alpha(*)$; for this it suffices to
prove $Y(\bar C') \backslash \alpha(*) = Y(\bar C'') \backslash
\alpha(*)$.  If $\alpha \in Y(\bar C') \backslash \alpha(*)$ let 
$\alpha = \|h\|_D$ where $h \in \dsize \prod_i C'_i$, and let
$A = \{i < \kappa:h(i) < g^*(i)\}$, so if $A \in (J_D(\bar \lambda,\kappa))^+$
then $\alpha \le \|h\|_{D+A} < \lambda$ and $\|h\|_{D+A} \le 
\|g^*\|_{D+A} < \alpha(*)$ but $\alpha \ge \alpha(*)$ contradiction.  
So $A \in J_D(\bar \lambda,\kappa)$ hence $A \notin D^+$ by clause (c), 
so $g^* \le_D h$.  Now clearly there is $h' =_D h$ with 
$h' \in \dsize \prod_{i < \kappa} C''_i$, so 
$\alpha = \|h\|_D = \|h'\|_D \in {\Cal C}(\bar C'')$.  
The other inclusion is easier.]
\endroster
\medskip

\noindent

\medskip
\roster
\item "{$(*)_3$}"  if $\bar C',\bar C'' \in \dsize \prod_{i < \kappa}
\text{ Club}(\lambda_i)$ and 
$\{i < \kappa:C''_i \subseteq \text{ acc}(C'_i)\} \in D$ \underbar{then}
\newline
${\Cal C}(\bar C'') \subseteq \text{ acc}({\Cal C}(\bar C'))$ \newline
[why?  let $\beta \in {\Cal C}[\bar C'']$ but $\beta \notin \text{ acc}
({\Cal C}(\bar C'))$ and we shall get a contradiction.  Clearly $\beta >
\sup({\Cal C}(\bar C') \cap \beta)$ (as $\beta \notin \text{ acc}({\Cal C}
(\bar C'))$.  As ${\Cal C}[\bar C'']$ is acc$(c\ell Y[\bar C''])$, clearly 
there is $\alpha \in Y[\bar C'']$ such that $\beta > \alpha > \sup({\Cal C}
(\bar C') \cap \beta) \ge \sup(Y(\bar C') \cap \beta)$, but 
$Y[\bar C''] = \{\|g\|_D:g \in \dsize \prod_{i < \kappa} C''_i\}$, 
so there is $g \in \dsize \prod_{i < \kappa} C''_i$ such that 
$\|g\|_D = \alpha$.  As $\{i:C''_i \subseteq \text{ acc}(C'_i)\} \in D$,
clearly

$$
B =: \{i < \kappa:g(i) \in \text{ acc}(C'_i)\} \in D.
$$ 
\medskip

\noindent
So if $h \in \dsize \prod_{i < \lambda} \lambda_i, h <_D g$ then we can
find $h' \in \dsize \prod_{i < \kappa} C'_i$ such that $h <_D h' <_D g$
(just $h'(i) = \text{ Min}(C'_i \backslash (h(i)+1)$ noting $B \in D$) hence
\newline
$\alpha = \|g\|_D = \text{ sup}\{\|h\|_D:h(i) \in g(i) \cap C'_i$ 
when $i \in B,
h(i) = \text{ Min}(C'_i) \text{ otherwise}\}$ and in this set there is no
last element and it is included in $Y(\bar C')$, so necessarily
$\alpha \in {\Cal C}(\bar C')$, contradicting the choice of $\alpha:\beta >
\alpha > \sup({\Cal C}(\bar C') \cap \beta)$.]
\medskip
\noindent
\item "{$(*)_4$}"  if $\bar C',\bar C'' \in \dsize \prod_{i < \kappa}
\text{ Club}(\lambda_i)$ and 
$\{i:C''_i \subseteq \text{ acc}(C'_i) \text{ mod } J^{bd}_{\lambda_i}\}
\in D$ \underbar{then} \newline
${\Cal C}(\bar C'') \subseteq \text{ acc}({\Cal C}(\bar C')) \text{ mod }
J^{bd}_\lambda$ \newline
[why?  by $(*)_2 + (*)_3$, i.e., define $C'''_i$ to be $C''_i \backslash 
g(i)$ where \nl
$g(i) =: \sup (C''_i \backslash \text{ acc}(C'_i)) +1)$ 
when $C''_i \subseteq \text{ acc}(C'_i)$ and the empty set otherwise.  Now by 
$(*)_2$ we know \newline
${\Cal C}(\bar C'') = {\Cal C}(\bar C''') \text{ mod } J^{bd}_\lambda$ and 
by $(*)_3$ we know \newline
${\Cal C}(\bar C''') \subseteq \text{ acc}({\Cal C}(\bar C'))$.]
\endroster
\medskip

\noindent
Now we can prove the conclusion of \scite{5.11}.  
Let $\langle C^i_\alpha:\alpha <
\mu_i \rangle$ witness $\mu_i < {\frak d}{\frak p}^{1+}_{\lambda_i}$ and
$\langle g_\alpha:\alpha < \mu \rangle$ witness $\mu < \text{ Depth}^+
(\dsize \prod_{i < \kappa} \lambda_i,<_D)$.  Let $C_\alpha =: {\Cal C}
(\langle C^i_{g_\alpha(i)}:i < \kappa \rangle)$ for $\alpha < \mu$.  So
$\langle C_\alpha:\alpha < \mu \rangle$ witness $\mu < {\frak d}{\frak p}
^{1+}_\lambda$. \hfill$\square_{\scite{5.11}}$
\enddemo
\bigskip

\proclaim{\stag{5.12} Theorem}  Assume
\medskip
\roster
\item "{$(a)$}"  $\kappa$ is regular uncountable
\sn
\item "{$(b)$}"  $\bar \lambda = \langle \lambda_i:i < \kappa \rangle$ is
a sequence of regular cardinals $> \kappa$
\sn
\item "{$(c)$}"  $D$ is a normal filter on $\kappa$ (or just
$\aleph_1$-complete)
\sn
\item "{$(d)$}"  $\lambda = \| \bar \lambda\|_D = \text{ tcf}(\dsize
\prod_{i < \kappa} \lambda_i/D),\lambda$ regular
\sn
\item "{$(e)$}"  $\mu_i < {\frak d}{\frak p}^{2+}_{\lambda_i}$
\sn
\item "{$(f)$}"  $\mu < \text{ Depth}^+(\dsize \prod_{i < \kappa} \mu_i,
<_D)$.
\endroster
\medskip

\noindent
\underbar{Then} $\mu < {\frak d}{\frak p}^{2+}_\lambda$.
\endproclaim
\bigskip

\demo{Proof}  Let $\langle f^i_\alpha:\alpha < \mu_i \rangle$ exemplify
$\mu_i < {\frak d}{\frak p}^{+2}_{\lambda_i}$, let
$\langle g_\alpha:\alpha < \mu \rangle$ exemplify $\mu < \text{ Depth}^+
(\dsize \prod_{i < \kappa} \mu_i,<_D)$, and let $\langle h_\zeta:\zeta < 
\lambda \rangle$ exemplify $\lambda = \text{ tcf}(\dsize \prod_{i < \kappa}
\lambda_i,<_D)$.

Now for each $\alpha < \mu$ we define $f_\alpha \in {}^\lambda \lambda$
as follows:

$$
f_\alpha(\zeta) = \|\langle f^i_{g_\alpha(i)}(h_\zeta(i)):i < \kappa \rangle
\|_D.
$$

\noindent
Clearly $f_\alpha(\zeta)$ is an ordinal and as $f^i_{g_\alpha(i)} \in
{}^{(\lambda_i)} \lambda_i$ clearly
$\langle f^i_{g_\alpha(i)}(h_\zeta(i)):i < \kappa \rangle <_D \langle
\lambda_i:i < \kappa \rangle$ hence $f_\alpha(\zeta) < \| \bar \lambda\|_D
= \lambda$, so really $f_\alpha(\zeta) < \lambda$, so
\medskip
\roster
\item "{$(*)_1$}"  $f_\alpha \in {}^\lambda \lambda$.
\endroster
\medskip

\noindent
The main point is to prove $\beta < \alpha < \mu \Rightarrow f_\beta
<_{J^{bd}_\lambda} f_\alpha$.

Suppose $\beta < \alpha < \mu$, then $g_\beta <_D g_\alpha$ hence
$A =: \{i < \kappa:g_\beta(i) < g_\alpha(i)\} \in D$ so $i \in A \Rightarrow
f^i_{g_\beta(i)} <_{J^{bd}_{\lambda_i}} f^i_{g_\alpha(i)}$.  We can define
$h \in \dsize \prod_{i < \kappa} \lambda_i$ by: \newline
$h(i)$ is
sup$\{\zeta + 1:f^i_{g_\beta(i)}(\zeta) \ge f^i_{g_\alpha(i)}(\zeta)\}$ if
$i \in A$ and $h(i)$ is zero otherwise.

But $\langle h_\zeta:\zeta < \lambda \rangle$ is $<_D$-increasing cofinal
in $(\dsize \prod_{i < \kappa} \lambda_i,<_D)$ hence there is $\zeta(*) <
\lambda$ such that $h <_D h_{\zeta(*)}$. \newline
So it suffices to prove:

$$
\zeta(*) \le \zeta < \lambda \Rightarrow f_\beta(\zeta) < f_\alpha(\zeta).
$$

\noindent
So let $\zeta \in [\zeta(*),\lambda)$, so

$$
B =: \{ i < \kappa:h(i) < h_{\zeta(*)}(i) \le h_\zeta(i) \text{ and }
i \in A\}
$$

\noindent
belongs to $D$ and by the definition of $A$ and $B$ and $h$ we have

$$
i \in B \Rightarrow f^i_{g_\beta(i)}(h_\zeta(i)) < f^i_{g_\alpha(i)}
(h_\zeta(i)).
$$

\noindent
So

$$
\langle f^i_{g_\beta(i)}(h_\zeta(i)):i < \kappa \rangle <_D
\langle f^i_{g_\alpha(i)}(h_\zeta(i)):i < \kappa \rangle
$$

\noindent
hence (by the definition of $\|-\|_D$)

$$
\|\langle f^i_{g_\beta(i)}(h_\zeta(i)):i < \kappa \rangle\|_D < \| \langle
f^i_{g_\alpha(i)}(h_\zeta(i)):i < \kappa \rangle\|_D
$$

\noindent
which means

$$
f_\beta(\zeta) < f_\alpha(\zeta).
$$

\noindent
As this holds for every $\zeta \in [\zeta(*),\lambda)$ clearly

$$
f_\beta <_{J^{bd}_\lambda} f_\alpha.
$$

\noindent
So $\langle f_\alpha:\alpha < \mu \rangle$ is $<_{J^{bd}_\lambda}$-increasing,
so we have finished. \hfill$\square_{\scite{5.12}}$
\enddemo
\bn

\ub{\stag{5.11a} Discussion}:  Now assumption (c) in \scite{5.11} 
(and in \scite{5.5}) is not so serious: 
once we quote \cite[V]{Sh:g} (to satisfy the assumption in the usual case we
are given $\lambda = \text{ cf}(\lambda),\mu < \lambda \le \mu^\kappa$,
cf$(\mu) = \kappa,(\forall \alpha < \mu)(|\alpha|^\kappa < \mu)$ and we like
to find $\langle \lambda_i:i < \kappa \rangle$, and normal $D$ such that
$\|\langle \lambda_i:i < \kappa \rangle\|_{D+A} = \lambda$.  E.g.
\cite[Ch.V]{Sh:g}) if SCH fails above $2^{2^\theta},
\theta$ regular uncountable, $D$ a normal filter on $\theta,\|f\|_D \ge
\lambda = \text{ cf}(\lambda) > 2^{2^\theta}$, (so if ${\Cal E} =$ family 
of normal filters on $\theta$, so ${\Cal E}$ is nice and 
$\text{rk}^3_E(f) \ge \|f\|_D \ge \lambda$), so $g_\kappa$ from
\cite[Ch.V,3.10,p.244]{Sh:g} is as required. \nl
Still we may note
\demo{\stag{5.11A} Fact}  Assume
\medskip
\roster
\item "{$(a)$}"  $D$ an $\aleph_1$-complete filter on $\kappa$
\sn
\item "{$(b)$}"  $f^* \in {}^\kappa\text{Ord}$ and cf$(f^*(i)) > 2^\kappa$
for $i < \kappa$.
\endroster
\medskip

\noindent
\underbar{Then} for any 
$\bar C = \langle C_i:i < \kappa \rangle,C_i$ a club of $f^*(i)$ and
$\alpha < \|f^*\|_D$ we can find $f \in 
\dsize \prod_{i < \kappa} C_i$ such that:
\mr
\item "{$(\alpha)$}"  $A \in (J_D(f^*,\kappa))^+ \Rightarrow 
\alpha < \|f\|_{D+A} = \|f\|_D < \|f^*\|_D$
\sn
\item "{$(\beta)$}"  $A \in J_D(f^*,\kappa) \cap D^+ \Rightarrow 
\|f\|_{D+A} \ge \|f^*\|_D$
\endroster
\enddemo
\bigskip

\demo{Proof}  We choose by induction on $\zeta \le \kappa^+$, a function
$f_\zeta$ and $\langle f_{\zeta,A}:A \in (J_D(f^*,\kappa))^+ \rangle$ 
such that:
\medskip
\roster
\item "{$(a)$}"  $f_\zeta \in \dsize \prod_{i < \kappa} C_i$
\sn
\item "{$(b)$}"  $\varepsilon < \zeta \Rightarrow \dsize \bigwedge_i
f_\varepsilon(i) < f_\zeta(i)$
\sn
\item "{$(c)$}"  for $\zeta$ limit $f_\zeta(i) = 
\underset {\varepsilon < \zeta}\to \sup f_\varepsilon(i)$
\sn
\item "{$(d)$}"  for $A \in (J_D(f^*,A))^+$, letting $\alpha_{\zeta,A} =:
\|f_\zeta\|_{D+A}$ we have \newline
$f_{\zeta,A} \in \dsize \prod_{i < \kappa}
f^*(i),\|f_{\zeta,A}\|_D > \alpha_{\zeta,A}$ and \nl
$f_{\zeta,A}(i) \ge f_\zeta(i)$ for $i < \kappa$
\sn
\item "{$(e)$}"  $f_{\zeta,A}(i) < f_{\zeta + 1}(i)$ for $i < \kappa,A \in
(J_D(f^*,A))^+$
\sn
\item "{$(f)$}"  $\|f_0\|_D \ge \alpha$ and \nl
$A \in J_D(f^*,\kappa) \Rightarrow \|f_0\|_{D+A} \ge \|f^*\|_D$.
\endroster
\medskip

\noindent
There is no problem to carry the definition: for defining $f_0$ for each
$A \in J_D(f^*,\kappa)$ choose $g_A <_{D+A} f^*$ such that
$\|g_A\|_{D+A} \ge \|f\|_D$ (possible as $\|f^*\|_{D+A} > \|f^*\|_D$ by the
assumption on $A$).  Let $g^* < f^*$ be such that $\|g^*\|_D \ge \alpha$,
(possible as $\alpha < \|f^*\|_D$) and let $f_0 \in \dsize \prod_{i < \kappa}
f^*(i)$ be defined by $f_0(i) = \text{ Min}(C_i \backslash \sup
\{g^*(i),g_A(i):A \in (J_D(f^*,\kappa))\})$.  For $\zeta$ limit 
there is no problem to define 
$f_\zeta$; and also for $\zeta$ successor.  If $f_\zeta$ is defined, we should
choose $f_{\zeta,A}$.  For clause (d) note that
$\|f^*\|_{D+A} = \|f^*\|_D$ as $A \in (J_D(f^*,A))^+$ and use the definition
of $\|f\|_D$.  We use, of course,
$\dsize \bigwedge_i$ cf$(f^*(i)) > 2^\kappa$.

Now $f_{\kappa^+}$ is as required.  Note: $f <_D f_{\kappa^+} \Rightarrow
\dsize \bigvee_{\zeta < \kappa^+} f <_D f_\zeta$, and for \newline
$A \in (J_D(f^*,\kappa))^+,\|f_{\kappa^+}\|_{D+A} = 
\underset {\zeta < \kappa^+}\to \sup \|f_\zeta\|_{D+A} 
= \underset {\zeta < \kappa^+}\to \sup \alpha_{\zeta,A} 
\le \underset {\zeta < \kappa^+}\to \sup \|f_{\zeta+1}\|_D =
\|f_{\kappa^+}\|_D$. \newline
${}$ \hfill$\square_{\scite{5.11A}}$
\enddemo
\bigskip

\demo{\stag{5.11B} Conclusion}  1) In \scite{5.11} we can weaken 
assumption (c) to
\mr
\item {$(c)^-$}"  $\|\langle \lambda_i:i < \kappa \rangle\|_D = \lambda$.
\ermn
2) In \scite{5.5} we can weaken assumption (c) to (c)a$^-$.
\enddemo
\bigskip

\demo{Proof}  1) In the proof of \scite{5.11}, choose $g^{**} \in
\dsize \prod_{i <\kappa} \lambda_i$ satisfying (exists by \scite{5.11A}):
\mr
\item "{$(*)_0$}"  $A \in J_D(\bar \lambda,\kappa) \cap D^+ \Rightarrow 
\|g^{**}\|_{D+A} \ge \lambda$ (which is $\|\bar \lambda\|_D$).
\ermn
We redefine $Y[\bar C]$ as $\{\|g\|_D:g \in \dsize \prod_{i < \kappa} C_i$
but $g(i) \ge g^{**}(i)$ for $i < \kappa\}$.  The only change is during the
proof of $(*)_2$ there.  Now if $\alpha \in Y[\bar C'] \backslash
\alpha(*)$ then there is $h \in \dsize \prod_{i < \kappa} \lambda_i$ such
that $[i < \kappa \Rightarrow h(i) \ge g^{**}(i)]$ and $\|h\|_D = \alpha$ 
and let $A = \{i < \kappa:h(i) < g^*(i)\}$.  Now if $A \in (J_D(\bar \lambda,
\kappa))^+$ we get a contradiction as there and if $A = \emptyset$ mod $D$
we finish as there.  So we are left with the case $A \in J_D(\bar \lambda,
\kappa) \cap D^+,\| \bar \lambda\|_{D+A} > \|\bar \lambda\|_D \le \lambda$ 
hence $\|g^{**}\|_{D+A} \le \lambda$ hence $\|h\|_{D+A} \le \lambda > \alpha$ 
hence necessarily $\|h\|_{D+(\kappa \backslash A)} = \alpha$ (as
$\|h\|_D = \text{ Min}\{\|h\|_{D+A},\|h\|_{D + (\kappa \backslash A)}\}$.
Now choose $h' \in \dsize \prod_{i < \kappa} \lambda_i$ by $h' \restriction 
(\kappa \backslash A) = h \restriction (\kappa \backslash A)$ and 
$[i \in A \Rightarrow h'(i) = \text{ Min}(C''_i
\backslash h(i))]$ so $h' \in \dsize \prod_{i < \kappa} C''_i,h \le h' <
\bar \lambda,\lambda \le \|h\|_{D+A} \le \|h'\|_{D+A} \le \|h'\|_{D+A}$ and
so \nl
$\|h'\|_D = \text{ Min}\{\|h'\|_{D+A},\|h'\|_{D + (\kappa \backslash A)}\} 
= \alpha$. \nl
So we are done. \nl
2) Let $g^{**}$  be as in the proof of part (1).  We let there

$$
\align
Y[\bar X] =: &\biggl\{ \|h\|_{D+A}:h \in \dsize \prod_{i < \kappa}
(\lambda_i \backslash g^{**}(i)) \text{ and } A \in I^+ \biggr\} \\
  &\text{remembering } I = J_D(\bar \lambda,\kappa).
\endalign
$$

$$
\align
{\Cal Y}[\bar X] =: \biggl\{ (h,A)/ \approx_I:&h \in \dsize \prod_{i < \kappa}
(X_i \backslash g^{**}(i)) \text{ and} \\
  &(h,A) \in \text{ cla}^\alpha_I(\lambda,D) 
\text{ for some } \alpha < \lambda \biggr\}
\endalign
$$
\mn
and we can restrict ourselves to sequences $\bar X$ such that $X_i \cap
g^{**}(i) = \emptyset$.  In the proof of $(*)_2$ make $g \ge g^{**}$.
\hfill$\square_{\sciteu{10.19}}$\sciteuphantom{10.19}
\enddemo
\bigskip

\proclaim{\stag{5.16} Claim}  Assume
\medskip
\roster
\item "{$(a)$}"  $J$ is a filter on $\kappa$
\sn
\item "{$(b)$}"  $\lambda$ a regular cardinal, $\lambda_i > 2^\kappa,
\theta > 2^\kappa$
\sn
\item "{$(c)$}"  $\dsize \prod_{i < \kappa} \lambda_i/J$ is $\lambda$-like
i.e.
{\roster
\itemitem{ $(i)$ }  $\lambda = \text{ tcf } \Pi \lambda_i/J$
\sn
\itemitem{ $(ii)$ }  $T_J(\langle \lambda_i:i < \kappa \rangle) = \lambda$
(follows from (i) + (iii) actually) and
\sn
\itemitem{ $(iii)$ }  if $\mu_i < \lambda_i$ then $T_J(\langle \mu_i:
i < \kappa \rangle) < \lambda$
\endroster}
\item "{$(d)$}"  $\kappa < \theta = \text{ cf}(\theta) < \lambda_i$ for
$i < \kappa$
\sn
\item "{$(e)$}"  $i < \kappa \Rightarrow S^{\lambda_i}_\theta = \{ \delta
< \lambda_i:\text{cf}(\delta) = \theta\} \in I[\lambda_i]$ (see below)
\sn
\item "{$(f)$}"  $(\forall \alpha < \theta)[|\alpha|^\kappa < \theta]$.
\endroster
\medskip

\noindent
\underbar{Then} $S^\lambda_\theta = \{ \delta < \lambda:\text{cf}(\delta) =
\theta\} \in I[\lambda]$.
\endproclaim
\bigskip

\remark{Remark}  Remember that for $\lambda$ regular uncountable

$$
\align
I[\lambda] = \biggl\{ A \subseteq \lambda:&\text{ for some club } E 
\text{ of } \lambda \bar{\Cal P} = \langle {\Cal P}_\alpha:\alpha < \lambda
\rangle, \\
  &\,{\Cal P}_\alpha \subseteq {\Cal P}(\alpha),|{\Cal P}| < \lambda \\
  &\text{ for every } \delta \in A \cap E,\text{cf}(\delta) < \delta
\text{ and for some closed} \\ 
  &\text{ unbounded subset } a \text{ of } \delta \text{ of order type} \\
  &< \delta,(\forall \alpha < \delta)(\exists \beta < \delta)(a \cap \alpha 
  \in {\Cal P}_\beta) \biggr\}.
\endalign
$$
\endremark
\bigskip

\demo{Proof}  Clearly each $\lambda_i$ is a regular cardinal and
$\lambda = \text{ tcf}(\dsize \prod_{i < \kappa} \lambda_i/J)$, so let \nl
$\bar f = \langle f_\alpha:\alpha < \lambda \rangle$ be a $<_J$-increasing
sequence of members of $\dsize \prod_{i < \kappa} \lambda_i$, which is
cofinal in $\dsize \prod_{i < \kappa} \lambda_i/J$.  So without loss of
generality if $\bar f \restriction \delta$ has a $<_J$-eub $f'$ then
$f_\delta =_J f'$.

For each $i < \kappa$ (see the references above) we can find $\bar e^i = 
\langle e^i_\alpha:\alpha < \lambda_i \rangle$ and $E_i$ such that:
\medskip
\roster
\widestnumber\item{(iii)}
\item "{$(i)$}"  $E_i$ is a club of $\lambda_i$
\sn
\item "{$(ii)$}"  $e^i_\alpha \subseteq \alpha$ and otp$(e^i_\alpha) \le
\theta$
\sn
\item "{$(iii)$}"  if $\beta \in e^i_\alpha$ \ub{then} $e^i_\beta = e^i_\alpha
\cap \beta$
\sn
\item "{$(iv)$}"  if $\delta \in E_i$ and cf$(\delta) = \theta$, \ub{then}
$\delta = \sup(e^i_\delta)$.
\endroster
\medskip

\noindent
Choose $\bar N = \langle N_i:i < \lambda \rangle$ such that $N_i \prec
({\Cal H}(\chi),\in,<^*_\chi)$ where, e.g., $\chi = \beth_8(\lambda)^+$, \nl
$\|N_i\| < \lambda,N_i$ is increasing continuous, $\bar N \restriction (i+1) 
\in N_{i+1},N_i \cap \lambda$ is an ordinal, and \newline
$\{ \bar f,J,\lambda,\langle \lambda_i:i < \kappa \rangle,\langle \bar e^i:
i < \kappa \rangle \} \in N_0$.  Let $E = \{ \delta < \lambda:
N_\delta \cap \lambda = \delta\}$, so it suffices to prove
\medskip
\roster
\item "{$(*)$}"  if $\delta \in E \cap S^\lambda_\theta$ then there is 
$a$ such that:
{\roster
\itemitem{ $(i)$ }  $a \subseteq \delta$
\sn
\itemitem{ $(ii)$ }  $\delta = \sup(a)$
\sn
\itemitem{ $(iii)$ }  $|a| < \lambda_0$
\sn
\itemitem{ $(iv)$ }  $\alpha < \delta \Rightarrow a \cap N_\alpha \in
N_\delta$.
\endroster}
\endroster
\medskip

By clause (b) in the assumption necessarily $\bar f \restriction \delta$ has
$a <_J$-eub (\cite[Ch.II,\S1]{Sh:g}) so necessarily $f_\delta$ is an $<_J$-eub 
of $\bar f \restriction
\delta$ hence $A^* = \{i < \kappa:\text{cf}(f_\delta(i)) = \theta\} =
\kappa \text{ mod } J$.  By clause (f) of the assumption for each 
$i \in A^*,e^i_{f_\delta(i)}$ is well
defined, and let $e^i_{f_\delta(i)} = \{ \alpha^i_\zeta:\zeta < \theta\}$
with $\alpha^i_\zeta$ increasing with $\zeta$.  For each $\zeta < \theta$
we have $\langle \alpha^i_\zeta:i < \kappa \rangle <_J f_\delta$ hence for
some $\gamma(\zeta) < \delta$ we have $\langle \alpha^i_\zeta:i < \kappa 
\rangle <_J f_{\gamma(\zeta)}$, but $T_D(f_{\gamma(\zeta)}) < \lambda$ and
$\gamma(\zeta) \in N_{\gamma(\zeta)+1}$ hence $f_{\gamma(\zeta)} \in
N_{\gamma(\zeta)+1}$ hence for some $g_\zeta <_J f_{\gamma(\zeta)}$ we have:
$g_\zeta \in N_{\gamma(\zeta)+1}$ and $A_\zeta = \{i < \kappa:g_\zeta(i) =
\alpha^i_\zeta\} \ne \emptyset \text{ mod } J$.  As $\theta = \text{ cf}
(\theta) > 2^\kappa$ for some $A \subseteq \kappa$ we have
$B =: \{ \zeta < \theta:A_\zeta = A\}$ is unbounded in $\theta$.
\medskip

\noindent
Now for $\zeta < \theta$ let

$$
\align
a_\zeta = \biggl\{ \text{ Min}\{ \gamma < \lambda:&\neg(f_\gamma
\le_{J+(\kappa \backslash A)} g)\}: \\
  &g \in \dsize \prod_{i < \kappa} \{\alpha^i_\varepsilon:\varepsilon <
\zeta\} = \dsize \prod_{i<\kappa} e^i_{(\alpha^i_\zeta)} \biggr\}.
\endalign
$$
\mn
Clearly $\zeta < \xi < \theta \Rightarrow a_\zeta \subseteq a_\xi$.  Also for
$\zeta < \theta,a_\zeta$ is definable from $\bar f$ and $g_\zeta \restriction
A$, hence belongs to $N_{\gamma(\zeta)+1}$, but its cardinality is
$\le \theta + 2^\kappa < \lambda$ hence it is a subset of 
$N_{\gamma(\zeta)+1}$.  Moreover, also $\langle a_\xi:\xi < \zeta \rangle$
is definable from $\bar f$ and $\left< \langle\{\alpha^i_\varepsilon:
\varepsilon < \xi \}:i < A \rangle:\xi \le \zeta \right>$ hence from $\bar f$
and $g_\zeta \restriction A$ and $\langle \bar e^i:i < \kappa \rangle$, all
of which belong to $N_0 \prec N_{\gamma(\zeta)+1}$, hence $\zeta \in B
\Rightarrow \langle a_\xi:\xi \le \zeta \rangle \in N_{\gamma(\zeta)+1}$, is
a bounded subset of $\delta$.  Now
\mr
\item "{$(*)$}"  $\dbcu_{\xi < \theta} a_\xi$ is unbounded in $\delta$ \nl
[why?  let $\beta < \delta$ so for some $\zeta < \theta$ we have: 
$$
f_\beta(i) < f_\delta(i) \Rightarrow f_\beta(i) < \alpha^i_\zeta <
f_\delta(i)
$$
so
$$
\text{Min}\{\gamma:\neg(f_\gamma \le_{J+(\kappa \backslash A)}\langle
\alpha^i_\zeta:i < \kappa \rangle)\} \in (\beta,\delta) \cap a_{\zeta +1}.
$$
\ermn
Let $w = \{\zeta < \theta:a_\zeta \text{ is bounded in }
a_{\zeta +1}\}$

$$
a'_\zeta = \bigl\{ \text{Min}\{\gamma \in a_{\xi +1}:\gamma \text{ is an upper
bound of } a_\xi\}:\xi < \zeta \bigr\}.
$$
\mn
So $\cup \{a'_\zeta:\zeta < \theta\}$ is as required.
\hfill$\square_{\scite{5.16}}$
\enddemo
\bigskip

\remark{\stag{5.16A} Remark}  1) If we want to weaken clause (c) in claim
\sciteu{10.20} retaining only (i) there (and omitting (ii) + (iii)), it 
is enough if we add:
\medskip
\roster
\item "{$(g)$}"  for each $i < \kappa$ and $\delta \in S^{\lambda_i}_\theta,
\{\gamma < \delta:\text{cf}(\gamma) > \kappa$ and $\gamma \in e^i_\delta\}$
is a stationary subset of $\delta$.
\ermn
2) In part (1) of this remark, we can apply cf$(\gamma) > \kappa$ by
cf$(\gamma) = \sigma$, \ub{if} $D$ is $\sigma^+$-complete or at least not
$\sigma$-incomplete. \nl
3) This is particularly interesting if $\lambda = \mu^+ = \text{ pp}(\mu)$.
\endremark
\newpage

\head {\S6 The class of cardinal ultraproducts modulo $D$} \endhead\resetall
\bigskip

We presently concentrate on ultrafilters (for filters: two versions).  This
continues \cite[\S3]{Sh:506}, see history there and in \cite{CK}, \cite{Sh:g}.
\bigskip

\demo{\stag{6.4} Fact}  Assume
\medskip
\roster
\item "{$(a)$}"  $D$ is an ultrafilter on $\kappa$ and $\theta = \text{ reg}
(D)$
\sn
\item "{$(b)$}"  $\mu = \text{ cf}(\mu)$ and $\alpha < \mu \Rightarrow
|\alpha|^{< \text{ reg}(D)} < \mu$
\sn
\item "{$(c)$}"  $\bar n = \langle n_i:i < \kappa \rangle,0 < n_i < \omega,
A^* = \dsize \bigcup_{i < \kappa}(\{i\} \times n_i)$
\sn
\item "{$(d)$}"  $Y \in D \Rightarrow \mu \le \text{ max pcf}\{\lambda_{i,n}:
i \in Y,n < n_i\}$
\sn
\item "{$(e)$}"  for each $i < \kappa,n < n$ we have $\lambda_{(i,n)}$ 
is regular $> \kappa$ strictly increasing
with $n$, stipulating $\lambda_{(i,n_i)} = \mu$. 
\endroster
\medskip

\noindent
\underbar{Then} for some $\langle m_i:i < \kappa \rangle \in \dsize
\prod_{i < \kappa}(n_i +1)$ we have:
\medskip
\roster
\item "{$(\alpha)$}"  $\mu \le \text{ tcf}(\dsize \prod_{i < \kappa}
\lambda_{(i,m_i)} /D)$
\sn
\item "{$(\beta)$}"  $\mu > \text{ max pcf}\{\lambda_{(i,n)}:i < \kappa
\text{ and } n < m_i\}$.
\endroster
\enddemo
\bigskip

\demo{Proof}  We try to choose by induction on $\zeta < \text{ reg}(D),
B_\zeta$ and $\langle n^\zeta_i:i < \kappa \rangle$ such that:
\medskip
\roster
\widestnumber\item{(iii)}
\item "{$(i)$}"  $B_\zeta \in D$
\sn
\item "{$(ii)$}"  $n^\zeta_i < n_i \text{ non-decreasing in } \zeta$
\sn
\item "{$(iii)$}"  $B_\zeta = \{i:n^\zeta_i < n^{\zeta +1}_i\}$ and
\sn
\item "{$(iv)$}"   max pcf$\{\lambda_{(i,n)}:i < \kappa \text{ and }
n \le n^\zeta_i\} < \mu$.
\endroster
\medskip

\noindent
If we succeed, then $\{B_\zeta:\zeta < \text{ reg}(D)\}$ exemplifies $D$ is
reg$(D)$-regular, contradiction.  During the induction we choose
$B_\zeta$ in step $\zeta + 1$.  For $\zeta = 0$ try $n^\zeta_i = 0$, if
this fails then $m_i = 0$ (for $i < \kappa$) is as required.  For $\zeta$
limit let $n^\zeta_i = n^\xi_i$ for every $\xi < \zeta$ large enough, this
is O.K. as \newline
max pcf$\{\lambda_{(i,n)}:i < \kappa \text{ and } n < n^\zeta_i\} \le
\dsize \prod_{\xi < \zeta}$ max pcf$\{\lambda_{(i,n)}:i < \kappa \text{ and }
n \le n^\xi_i\} < \mu$ by assumption (b).  Lastly, for $\zeta = \xi + 1,\{i <
\kappa:n^\xi_i < n_i\} \in D$ (otherwise contradiction as
$\lambda_{(i,n_i)} = \mu$ and assumption (d)), and if 
$\mu \le \text{ tcf}(\dsize \prod_{i < \kappa} \lambda_{n^\xi_i + 1}/D)$ 
we are done with $m_i = 
n^\xi_i + 1$, if not there is $B_\xi \in D$ such that 
max pcf$\{\lambda_{n^\xi_i + 1}:i \in B\} < \mu$ and let

$$
n^\zeta_i = \cases n^\xi_i + 1 \quad &\text{ \underbar{if} } \quad
i \in B_\xi,n^\xi_i < n_i \\
n^\xi_i \quad &\text{ \underbar{if} } \quad \text{otherwise}.
\endcases
$$
${}$\hfill$\square_{\scite{6.4}}$
\enddemo
\bigskip

\proclaim{\stag{6.5} Lemma}    Assume
\medskip
\roster
\widestnumber\item{(iii)}
\item "{$(i)$}"  $D$ is an ultrafilter on $\kappa$
\sn
\item "{$(ii)$}"  $\mu = \text{ cf}(\mu)$ and $\alpha < \mu \Rightarrow
|\alpha|^{< \text{ reg}(D)} < \mu$
\sn
\item "{$(iii)$}"  at least one of the following occurs:
{\roster
\itemitem{ $(\alpha)$ }  $\alpha < \mu \Rightarrow |\alpha|^{\text{reg}(D)}
< \mu$
\sn
\itemitem{ $(\beta)$ }  $D$ is closed under decreasing sequences of length
$\theta$. 
\endroster}
\endroster
\medskip

\noindent
\underbar{Then} there is a minimal $g/D$ such that: \newline
$\mu = \text{ tcf }\left( \dsize \prod_{i < \kappa} g(i)/D \right)$ and
$\dsize \bigwedge_{i < \kappa} \text{ cf}(g(i)) > \kappa$. \nl
We shall prove it somewhat later.
\endproclaim
\bigskip

\remark{\stag{6.5A} Remark}  1) Note that necessarily (in \scite{6.5})

$$
\{i < \kappa:g(i) \text{ a regular cardinal}\} \in D.
$$
\medskip

\noindent
2) $g$ is also $<_D$-minimal under: $\mu \le \text{ tcf }
\left(\dsize \prod_{i < \kappa}g(i)/D \right) \and \{i:\text{cf}(g(i)) >
\kappa\} \in D$. \newline
[Why?  assume $g' <_D g_\beta,\mu \le \text{ tcf}\left( 
\dsize \prod_{i < \kappa} g'(i)/D \right)$, and
$X = \{i:\text{cf}(g(i)) \le \kappa\} = \emptyset \text{ mod } D$; clearly 
$\mu \le \text{ tcf}\left(\dsize \prod_{i < \kappa} \text{ cf}(g'(i))/D 
\right)$.  If Lim$_D\text{ cf}(g'(i))$ is singular, by 
\cite[II,1.5A]{Sh:g} for some $\langle \lambda_i:i < \kappa \rangle$, we have
$\mu = \text{ tcf}(\Pi \lambda_i/D)$ and \nl
Lim$_D \lambda_i = \text{ Lim}_D \text{ cf}(g(i))$ and 
$(\forall i)[\text{cf}(g(i)) > \kappa \rightarrow \lambda_i \ge \kappa$], so
again without loss 
of generality $\dsize \bigwedge_{i < \kappa} \lambda_i > \kappa$.  Now 
$\langle \lambda_i:i < \kappa \rangle$ contradicts the choice of $g$.  If
Lim$_D \text{cf}(g(i))$ is regular, it is $\mu$ and all is easier].
\newline
3)  If $|\kappa^\kappa/D| < \mu$ then we can omit (in \scite{6.5} and 
\scite{6.5A}(2)) the clause \newline
``$\{i:\text{cf}(g(i)) > \kappa\} \in D$".
\endremark
\bigskip

\demo{\stag{6.6} Conclusion}  If assumptions (i)-(iii) of \scite{6.5} 
hold and 
\medskip
\roster
\item "{$(iv)$}"  $\mu > 2^\kappa$
\endroster
\medskip

\noindent
\ub{then} without loss of generality each $g(i)$ is a regular cardinal and
$\left( \dsize \prod_{i < \kappa} g(i)/D,<_D\right)$ is $\mu$-like (i.e. of
cardinality $\mu$ but every proper initial segment has smaller cardinality.
\enddemo
\bigskip

\remark{\stag{6.6A} Remark}  We use $\mu > 2^\kappa$ in \scite{6.6} rather 
than $\mu > |\kappa^\kappa/D|$ as \scite{6.5A}(3) (which concerns
\scite{6.5}, \scite{6.5A}(3)) as the proof of \scite{6.6} uses \scite{1.3}.
\endremark
\bigskip

\demo{Proof of \scite{6.6}}  If $D$ is $\aleph_1$-complete this is 
trivial, so assume not hence reg$(D) > \aleph_0$.

Let $g \in {}^\kappa(\mu + 1)$ be as in \scite{6.5}, so without loss 
of generality as in \scite{6.5A}(2), and remember \scite{6.5A}(1) so
without loss of generality each $g(i)$ is a regular cardinal.  
Clearly $\dsize \prod_{i < \kappa}g(i)$ has cardinality $\ge \mu$.  
Assume first $\mu = \chi^+$. \newline
Let $g' < \dsize \prod_{i < \kappa} g(i)$, then by \scite{6.5A}(3) and 
choice of $g$

$$
\sup\{\text{tcf } \Pi \lambda_i/D:\lambda_i \le g'(i) \text{ for }
i < \kappa\} \le \chi.
$$
\medskip

\noindent
But as reg$(D) > \aleph_0$ by clause (ii) of the assumption we have 
$\alpha < \mu \Rightarrow |\alpha|^{\aleph_0} < \mu$ so \scite{2.5} 
applies and $| \dsize \prod_{i < \kappa} g'(i)/D| \le \chi$, so really 
$\dsize \prod_{i < \kappa} g(i)/D$ is $\mu$-like.

If $\mu$ is not a successor, then it is weakly 
inaccessible and $\mu = \sup(Z)$, where \newline
$Z = \{\chi^+:\kappa^\kappa/D < \chi^{\aleph_0} = \chi < \mu\}$, so for each
$\chi \in Z$ we can find $g_\chi \in {}^\kappa(\mu + 1)$ such that
$\dsize \prod_{i < \kappa} g_\chi(i)/D$ is $\chi$-like so necessarily for
$\chi_1 < \chi_2$ in $Z$ we have $g_{\chi_1} <_D g_{\chi_2}$.  It is enough 
to find a $<_D$-lub for $\langle f_\chi:\chi \in Z \rangle$, and as 
$\mu > 2^\kappa$ this is immediate. \hfill$\square_{\scite{6.6}}$
\enddemo
\bigskip

\demo{Proof of \scite{6.5}}  First try to choose, by induction on 
$\alpha,f_\alpha$ such that:
\medskip
\roster
\item "{$(A)$}"  $f_\alpha \in {}^\kappa(\mu + 1)$
\sn
\item "{$(B)$}"  $\mu = \text{ tcf}\left( \dsize \prod_{i < \kappa} f_\alpha
(i)/D \right)$
\sn
\item "{$(C)$}"  $\beta < \alpha \Rightarrow f_\alpha <_D f_\beta$
\sn
\item "{$(D)$}"  each $f_\alpha(i)$ is a regular cardinal $> \kappa$. 
\endroster
\medskip

\noindent
Necessarily for some $\alpha^*$ we have: $f_\alpha$ is well defined iff
$\alpha < \alpha^*$.  Now $\alpha^*$ cannot be zero as the constant function
with value $\mu$ can serve as $f_0$.  Also if $\alpha^*$ is a successor
ordinal, say $\alpha^* = \beta + 1$, then $f_\beta$ is as required in the
desired conclusion. 

So $\alpha^*$ is a limit ordinal, and by passing to a subsequence, without 
loss of generality $\alpha^* =
\text{ cf}(\alpha^*)$ and call it $\theta$. \newline
Without loss of generality
\mr
\item "{$(E)$}"  $\mu = \text{ max pcf}\{f_\alpha(i):i < \kappa\}$.
\ermn
We now try to choose by induction on $\zeta < \text{ reg}(D)$ the objects
$\alpha_\zeta,A_\zeta,{\frak b}_\zeta$ such that:
\medskip
\roster
\item "{$(a)$}"  $\alpha_\zeta < \theta$ is strictly increasing with $\zeta$
\sn
\item "{$(b)$}"  $A_\zeta \in D$
\sn
\item "{$(c)$}"  ${\frak b}_\zeta \subseteq \{f_{\alpha_\xi}(i):\xi \le \zeta,
\text{ and } i \in A_\xi\}$
\sn
\item "{$(d)$}"  ${\frak b}_\zeta$ is increasing with $\zeta$
\sn
\item "{$(e)$}"  max pcf$({\frak b}_\zeta) < \mu$
\sn
\item "{$(f)$}"  for each $i$ the sequence \newline
$\langle f_{\alpha_\xi}(i):\xi \le
\zeta \text{ and } i \in A_\xi \text{ and } f_{\alpha_\xi}(i) \notin
{\frak b}_\zeta \rangle$ is strictly decreasing
\sn
\item "{$(g)$}"  $\alpha_0 = 0, A_0 = \kappa,{\frak b}_\zeta = \emptyset$
\sn
\item "{$(h)$}"  $\alpha_{\zeta + 1} = \alpha_\zeta + 1$ and \nl
$A_{\zeta + 1} = \{i \in A_\zeta:f_{\alpha_{\zeta + 1}}(i) < 
f_{\alpha_\zeta}(i) \text{ and } f_{\alpha_\zeta}(i) \in {\frak b}_\zeta\}$
\sn
\item "{$(i)$}"  for $\zeta$ limit, $\alpha_\zeta$ is the first $\alpha <
\theta$ which is $\ge \dsize \bigcup_{\varepsilon \le \zeta} \alpha_
\varepsilon$ such that for some $B \in D$ we have: \newline
$\mu > \text{ max pcf}\{f_{\alpha_\xi}(i):\xi < \zeta,i \in A_\xi 
\text{ and } i \in B \text{ and } f_{\alpha_\xi}(i) \le
f_\alpha(i)\}$
\sn
\item "{$(j)$}"  ${\frak b}_{\zeta + 1} = {\frak b}_\zeta$
\sn
\item "{$(k)$}"  for $\zeta$ limit $A_\zeta$ satisfies the requirements on
$B$ in clause (i) and \newline
${\frak b}_\zeta = \dbcu_{\varepsilon < \zeta} {\frak b}_\varepsilon 
\bigcup \{f_\xi(i):\xi < \zeta \text{ and } i \in A_\zeta,A_\xi \cap
A_\zeta \text{ and } f_{\alpha_\xi}(i) \le f_{\alpha_\zeta}(i)\}$
\sn
\item "{$(\ell)$}"  for $\xi \le \zeta$ we have $\{i \in A_\xi:f_{\alpha_\xi}
(i) \notin {\frak b}_\zeta\} \in D$.
\endroster
\medskip

So for some $\zeta^* \le \text{ reg}(D)$ we have $(\alpha_\zeta,A_\zeta,
{\frak b}_\zeta)$ is well defined iff $\zeta < \zeta^*$. \newline
We check the different cases and get a contradiction in each (so $\alpha^*$
had necessarily been a successor ordinal giving the desired conclusion).
\enddemo
\bigskip

\noindent
\underbar{CASE 1}:  $\zeta^* = 0$.

We choose $\alpha_0 = 0,A_0 = \kappa,{\frak b}_0 = \emptyset$; so clause (g)
holds, first part of clause (a) (i.e. $\alpha_\zeta < \theta$) holds, 
clause (b)
and clause (c) are totally trivial, clause (e) holds as max pcf$(\emptyset)
= 0$ (formally we should have written sup pcf$({\frak b}_\zeta)$), clause (f)
speaks on the empty sequence, and the other clauses are empty in this case.
\bigskip

\noindent
\underbar{CASE 2}:  $\zeta^* = \zeta + 1$.

We choose $\alpha_{\zeta^*} = \alpha_{\zeta + 1} = \alpha_\zeta + 1,
A_{\zeta^*} = \{i \in A_\zeta:f_{\alpha_{\zeta^*}+1}(i) < f_{\alpha_\zeta}
(i) \text{ and } f_{\alpha_\zeta}(i) \notin {\frak b}_\zeta\}$ and 
${\frak b}_{\zeta + 1} \supseteq {\frak b}_\zeta$ is defined by clause (j).  
Clearly $\alpha_\zeta < \alpha_{\zeta + 1} < \theta$ and
$A_{\zeta + 1} \in D$ as $A_\zeta \in D$ and $f_{\alpha_\zeta + 1} <_D
f_{\alpha_\zeta}$ and $\{i:f_{\alpha_\zeta}(i) \notin {\frak b}_\zeta\} \in
D$ by clause $(ell)$; so clause (b) holds.  Now clause (a) holds 
trivially and clauses (g) and (i) are irrelevant.   
Clause (h) holds by our choice.

For clause (f), the new cases are when $f_{\alpha_{\zeta + 1}}(i)$ appears in
the sequence, i.e., $i \in A_{\zeta + 1}$ such that 
$f_{\alpha_{\zeta + 1}}(i) \notin \dsize \bigcup_{\xi \le \zeta + 1}
{\frak b}_\xi = {\frak b}_{\zeta + 1} = {\frak b}_\zeta$ but $i \in 
A_{\zeta + 1} \Rightarrow i \in A_\zeta \and f_{\alpha_\zeta}(i) \notin
{\frak b}_\zeta$ so also $f_{\alpha_\zeta}(i)$ appears in the sequence and as
$i \in A_{\zeta + 1} \Rightarrow f_{\alpha_\zeta}(i) > 
f_{\alpha_{\zeta + 1}}(i) = f_{\alpha_{\zeta + 1}}(i)$ plus the induction
hypothesis; we are done.

As for clause $(\ell)$ for 
$\xi \le \zeta + 1$, clearly \newline
$\{i \in A_\xi:f_{\alpha_\xi}(i) \notin {\frak b}_{\zeta + 1}\} = A_\xi \cap
\{i < \kappa:f_{\alpha_\xi}(i) \notin {\frak b}_{\zeta + 1}\}$.  Now the first
belongs to $D$ by clause (b) proved above and the second belongs to $D$ as
max pcf$({\frak b}_{\zeta + 1}) < \mu$ by clause (e) proved below as
tcf$\left( \dsize \prod_{i < \kappa} f_{\alpha_\xi}(i)/D \right) = \mu$ by
clause (B).

We have chosen ${\frak b}_{\zeta + 1} = {\frak b}_\zeta$, so 
(using the induction hypothesis) clauses (c), (d), (e) trivially hold and 
also clause (j) holds by the choice of $b_{\zeta^*}$, and (k) irrelevant 
so we are done.  
\bigskip

\noindent
\underbar{CASE 3}:  $\zeta^* = \zeta$ is a limit ordinal $< \text{ reg}(D)$.

Let ${\frak b}_\zeta = \dsize \bigcup_{\xi < \zeta} {\frak b}_\xi$, 
so clause (j) holds; also clauses (c), (d) hold trivially and clause 
(e) holds by basic pcf:

$$
\text{max pcf}({\frak b}_\zeta) \le \dsize \prod_{\xi < \zeta}
\text{ max pcf}({\frak b}_\xi) < \mu
$$

\noindent
as

$$
\mu = \text{ cf}(\mu) \and (\forall \alpha < \mu)[|\alpha|^{< \text{ reg}(D)}
< \mu)] \and \zeta < \text{ reg}(D).
$$
\medskip

\noindent
Now we try to define $\alpha_\zeta$ by clause (i).
\bigskip

\noindent
\underbar{SUBCASE 3A}:  $\alpha_\zeta$ is not well defined.

Let $w_i = \{ \xi < \zeta:i \in A_\xi \text{ and } f_{\alpha_\xi}(i)
\notin {\frak b}_\zeta\}$.  Note that by the induction hypothesis (clause (f))
for each $\varepsilon < \zeta$ and $i < \kappa$ we have the sequence 
$\langle f_{\alpha_\xi}(i):\xi < \varepsilon \text{ and } i \in A_{\xi} 
\text{ and }
f_{\alpha_\xi}(i) \notin {\frak b}_\varepsilon \rangle$ is strictly 
decreasing, so as ${\frak b}_\varepsilon \subseteq {\frak b}_\zeta$ clearly
$\langle f_{\alpha_\xi}(i):\xi < \varepsilon \text{ and } \xi \in
w_i \rangle$ is strictly decreasing.  As this holds for each $\varepsilon <
\zeta$ and $\zeta$ is a limit ordinal, clearly $\langle f_{\alpha_\xi}(i):\xi
\in w_i \rangle$ is strictly decreasing hence $w_i$ is finite.

Now for each $B \in D$ we have (first inequality by the assumption of the
subcase, second by the definition of the $w_i$'s)

$$
\align
\mu &\le \text{ max pcf} \biggl\{ f_\xi(i):\xi < \zeta,i \in A_\xi 
\text{ and } i \in B \biggr\} \\
  &\le \text{ max}\biggl\{\text{max pcf}({\frak b}_\zeta),\text{ max pcf}
\{f_\xi(i):\xi \in w_i \text{ and } i \in B\} \biggr\},
\endalign
$$
\medskip

\noindent
and max pcf$({\frak b}_\zeta) < \mu$ as said above, hence necessarily
\mr
\item "{$(*)$}"  $B \in D \Rightarrow \mu \le 
\text{ max pcf}\{f_{\alpha_\xi}(i):\xi \in w_i \text{ and } i \in B\}$.
\ermn
As $w_i$ is finite and each $f_\alpha(i)$ is a regular cardinal $> \kappa$ we
have $\{i:w_i \ne \emptyset\} \in D$.

By Claim \scite{6.4} (the case there of $\{i:m_i = n_i\} \in D$ is impossible
by $(*)$ above) we can find 
$g \in \dsize \prod_{i < \kappa} w_i/D$, more exactly 
$g \in {}^\kappa\text{Ord},w_i \ne \emptyset \Rightarrow g(i) \in w_i$ 
and $B \in D$ such that:
\medskip
\roster
\item "{$(\alpha)$}"  $\mu \le \text{ tcf}\left( \dsize \prod_{i < \kappa}
g(i)/D \right)$
\sn
\item "{$(\beta)$}"  $\mu > \text{ max pcf}\{f_{\alpha_\xi}(i):\xi \in w_i
\text{ and } i \in B \text{ and } f_{\alpha_\xi}(i) < g(i)\}$.
\endroster
\medskip

\noindent
Now by the choice of $\langle f_\alpha:\alpha < \theta \rangle$ and clause
$(\alpha)$ necessarily (and \cite[Ch.II,1.5A]{Sh:g}) 
for some $\alpha < \theta$ we have $f_\alpha <_D g$.  Now let
$B^\xi_\alpha = \{i < \kappa:f_\alpha(i) \ge f_{\alpha_\xi}(i)\}$, if
$B^\xi_\alpha \in D$ then $B^* = \{i < \kappa:\xi \in w_i \text{ and }
i \in B \text{ and } g_\alpha(i) > f_{\alpha_\xi}(i)\} \supseteq
\{i < \kappa:i \in A_\xi\} \cap \{i < \kappa:f_\alpha(i) \notin 
{\frak b}_\zeta\} \cap \{i < \kappa:f_\alpha(i) \ge f_{\alpha_\xi}(i)\}$
which is the intersection of three members of $D$ hence belongs to $D$, but
$\{f_{\alpha_\xi}(i):i \in B^*\}$ is included in the set in the right side of
clause $(\beta)$ hence $\mu > \text{ max pcf}\{f_{\alpha_\xi}(i):i \in
B^*\}$ contradicting $B^* \in D$, tcf$(\dsize \prod_{i < \kappa}
f_{\alpha_\xi}(i)/D) = \mu$.  So necessarily $B^\xi_\alpha \notin D$, hence
$f_\alpha <_D f_{\alpha_\xi}$ hence $\alpha < \alpha_\xi$.  So
$\dsize \bigcup_{\xi < \zeta} \alpha_\xi < \alpha < \theta$.
Let $B' = B \cap \{i < \kappa:f_\alpha(i) < g(i)\}$ so $B' \in D$ and

$$
\align
\biggl\{ f_{\alpha_\xi}(i):&\xi < w_i \text{ and } i \in B' \text{ and }
f_{\alpha_\xi}(i) \le f_\alpha(i) \biggr\} \subseteq \\
  &\biggl\{ f_{\alpha_\xi}(i):\xi < w_i \text{ and } i \in B \text{ and }
f_{\alpha_\xi}(i) < g(i) \biggr\} \subseteq \\
  &{\frak b}_\zeta \cup \biggl\{ f_{\alpha_\xi}(i):\xi \in w_i \text{ and }
f_{\alpha_\xi}(i) < g(i) \biggr\}
\endalign
$$
\mn
hence

$$
\align
&\text{max pcf} \biggl\{ f_{\alpha_\xi}(i):\xi < \zeta,i \in A_\xi 
\text{ and } i \in B' 
\text{ and } f_{\alpha_\xi}(i) \le f_\alpha(i) \biggr\} \le \\
  &\text{max} \biggl\{\text{max pcf}({\frak b}_\zeta),\text{max pcf}
\{ f_{\alpha_\xi}(i):\xi \in w_i \text{ and } i \in B \text{ and} \\
  &\qquad \qquad \qquad \qquad \qquad \qquad \qquad \qquad
f_{\alpha_\xi}(i) < g(i)\} \biggr\} < \mu
\endalign
$$

\noindent
(the first term is $< \mu$ as clause (e) 
was proved in the beginning of Case 3, the second term is $< \mu$ by 
clause $(\beta)$).  So $\alpha$ is as required in
clause (i) so $\alpha_\zeta$ is well defined; contradiction.
\bigskip

\noindent
\underbar{CASE 3B}:  $\alpha_\zeta$ is well defined.

Let $B \in D$ exemplify it.  We choose $A_\zeta$ as $B$ and we define
${\frak b}_\zeta$ by clause (k).

Now clause (a) follows from clause (i) (which holds by the assumption of
the subcase), clause (b) holds by the choice of $B$ (and of $A_\zeta$),
clause (c) by the choice of ${\frak b}_\zeta$, clause (d) by the choice of
${\frak b}_\zeta$, clause (e) by the choice of ${\frak b}_\zeta$. Now for
clause (f) by the induction hypothesis and clause (d) we should consider
only $f_{\alpha_\xi}(i) > f_{\alpha_\zeta}(i)$ when $\xi < \zeta,i \in
A_\xi \cap A_\zeta$ and $f_{\alpha_\xi}(i),f_{\alpha_\zeta}(i) \notin
b_\zeta$, but clauses (i) + (k) (i.e. the choice of ${\frak b}_\zeta$) take
care of this clauses (g), (h), (j) are irrelevant and $(\ell)$ follows from
(e).

So we are done.
\bigskip

\noindent
\underbar{CASE 4}:  $\zeta^* = \text{ reg}(D)$.

The proof is split according to the two cases in the assumption (iii).
\bigskip

\noindent
\underbar{SUBCASE 4A}:  $\alpha < \mu \Rightarrow |\alpha|^{\text{reg}(D)}
< \mu$.

Let ${\frak b} = \cup\{{\frak b}_\xi:\xi < \zeta^*\}$ so max pcf$({\frak b})
< \mu$, hence for each $\xi < \zeta^*$ we have $A'_\xi =: \{i \in A_\xi:
f_{\alpha_\xi}(i) \notin {\frak b}\} \in D$.  Let $w_i = \{ \xi < \zeta^*:
i \in A'_\xi \text{ and } f_{\alpha_\xi}(i) \notin {\frak b}\}$.  Now
for any $\zeta < \zeta^*$ and $i < \kappa$ the sequence $\langle 
f_{\alpha_\xi}(i):\xi < \zeta \text{ and } \xi \in w_i \rangle$ is strictly
decreasing (by clause (f)) hence 
$\langle f_{\alpha_\xi}(i):\xi < \zeta^* \text{ and } \xi \in w_i \rangle$
is strictly decreasing hence $w_i$ is finite.  Also for each $\xi <
\zeta^*$ the set $A'_\xi$ belongs to $D$, so $\{A'_\xi:\xi < \zeta^*\}$
exemplifies $D$ is $|\zeta^*|$-regular, but $\zeta^* = \text{ reg}(D)$,
contradiction.
\bigskip

\noindent
\underbar{SUBCASE 4B}:  $D$ is closed under decreasing sequences of length 
reg$(D)$.

Let ${\frak b} = \dsize \bigcup_{\zeta < \zeta^*} {\frak b}_\zeta$.

In this case, for each $\xi < \zeta^*$, the sequence
$\langle \{i \in A_\xi:f_{\alpha_\xi}(i) \notin {\frak b}_\zeta\}:\zeta \in
[\xi,\zeta^*]\} \rangle$ is a decreasing sequence of length $\zeta^* =
\text{ reg}(D)$ of members of $D$ so the intersection, $A'_\xi =
\{i \in A_\xi:f_{\alpha_\zeta}(i) \notin {\frak b}\} \in D$, and we continue
as in the first subcase. \hfill$\square_{\scite{6.5}}$
\bn

\definition{\stag{6.7} Definition}  1) For an ultrafilter $D$ on $\kappa$ let
reg$'(D)$ be: reg$(D)$ if $D$ is closed under intersection of decreasing
sequences of length reg$(D)$ and (reg$(D))^+$ otherwise. \newline
2) $\text{reg}''(D)$ is: reg$(D)$ \underbar{if} $(a)^-$ below holds and
(reg$(D))^+$ otherwise
\medskip
\roster
\item "{$(a)$}"  $\text{reg}'(D) = \text{ reg}(D)$ or just
\sn
\item "{$(a)^-$}"  letting $\theta = \text{ reg}(D)$, in $\theta^\kappa/D$
there is a $<_D$-first function above the constant functions.
\endroster
\enddefinition
\bigskip

\proclaim{\stag{6.8} Theorem}  If $D$ is an ultrafilter on $\kappa$ and 
$\theta = \text{ reg}'(D)$ \underbar{then} \newline
$\mu = \mu^{< \theta} \ge |2^\kappa/D| \Rightarrow \mu \in \{\Pi \lambda_i
/D:\lambda_i \in \text{ Card}\}$.
\endproclaim
\bigskip

\demo{Proof}  Apply Lemma \scite{6.6} with $D,\kappa,\mu^+$ here standing for
$D,\kappa,\mu$ there; note that assumption (iii) there holds as the definition
of reg$'(D)(=\theta)$ was chosen appropriately.

Let $g^*/D = \langle \lambda^*_i:i < \kappa \rangle$ be as there, so
as $\left( \dsize \prod_{i < \kappa} \lambda^*_i/D \right)$ is $\mu^+$-like,
for some $f \in \dsize \prod_{i < \kappa} \lambda_i$, we have
$| \dsize \prod_{i < \kappa} f(i)/D| = \mu$ as required. 
\hfill$\square_{\scite{6.8}}$
\enddemo
\bigskip

\remark{Remark}  Can reg$'(D) \ne \text{ reg}(D)$?  This is 
equivalent to: $D$ is not closed
under intersections of decreasing sequences of length $\theta =
\text{ reg}(D)$.  So if reg$'(D) \ne \text{ reg}(D) = \theta$ then $\theta$
is regular and for some function $\bold i:\kappa \rightarrow \theta$ the
ultrafilter \nl
$D' = \{ A \subseteq \theta:\bold i^{-1}(A) \in D\}$ is an
ultrafilter on $\theta$, with reg$(D') = \theta$ so $D'$ is not regular.
\medskip

\noindent
This leads to the well known problem (Kanamori \cite{Kn}) if $D$ 
is a uniform ultrafilter on $\kappa$ with reg$(D) = \kappa$ does
$\kappa^\kappa/D$ have a first function above the constant ones.
\endremark
\bigskip

\demo{\stag{6.10} Fact}  If $\mu = \theta = \text{ reg}(D) < \text{ reg}'(D),
\mu = \dsize \sum_{i < \theta} \mu_i,\mu^\kappa_i = \mu_i < \mu_{i+1}$ and

$$
|\dsize \prod_{i < \kappa} f(i)/D| \ge \mu \text{ then }
|\dsize \prod_{i < \kappa} f(i)/D| \ge \mu^\theta = \mu^\kappa.
$$
\enddemo

\newpage
    
REFERENCES.  
\bibliographystyle{lit-plain}
\bibliography{lista,listb,listx,listf,liste}

\enddocument

\bye